\documentclass[12pt,english]{article}
\usepackage{amssymb}
\usepackage{amsfonts}
\usepackage{amsmath}
\usepackage{latexsym}
\usepackage{indentfirst}
\usepackage{fancyhdr}
\usepackage{stmaryrd}
\usepackage{mathrsfs}
\usepackage{cite}
\usepackage{babel}
\makeatletter
\date{}
\numberwithin{equation}{section}
\newtheorem{theorem}{Theorem}[section]
\newtheorem{lemma}[theorem]{Lemma}
\newtheorem{defi}[theorem]{Definition}
\newtheorem{coral}[theorem]{Corallary}
\newtheorem{pro}[theorem]{Proposition}

\newcommand{\al}{\alpha}
\newcommand{\ga}{\gamma}
\newcommand{\sgm}{\sigma}
\newcommand{\be}{\beta}
\newcommand{\G}{\Gamma}

\newcommand{\mbb}{\mathbb}
\newcommand{\mcr}{\mathscr}

\newcommand{\kn}{\mbox{ker}\:}

\newcommand{\ove}{\overline}
\newcommand{\ptl}{\partial}

\makeatletter
\newcommand{\rmnum}[1]{\romannumeral#1}
\newcommand{\Rmnum}[1]{\expandafter\@slowromancap\romannumeral#1@}
\makeatother

\begin{document}

\title{Some Representations of Nongraded\\ Divergence-Free Lie Algebras }
\author{Ling \uppercase{Chen}
\thanks{ E-mail: chenling@amss.ac.cn (L. Chen)}\\%EndAName
{\small{ Beijing International Center for Mathematical Research, Peking University,}}\\
{\small{Beijing 100871, P. R. China}}}\maketitle

\begin{abstract}

 Divergence-free Lie algebras are originated from the Lie algebras of
volume-preserving transformation groups. Xu constructed a certain
nongraded generalization,  which may not contain any toral Cartan
subalgebra. In this paper, we give a complete classification of the
generalized weight modules over these algebras with weight
multiplicities less than or equal to one.

\end{abstract}

\section{Introduction}

The Lie algebras of volume-preserving transformation groups consist
of divergence-free first-order differential operators. They are
often called the special Lie algebras of Cartan type, or
divergence-free Lie algebras (see \cite{SX}). These algebras play
important roles in geometry, Lie groups and dynamics. Various graded
generalizations of them were studied by Kac \cite{Ka1}, Osborn
\cite{O}, Dokovic and Zhao \cite{DZ3}, and Zhao \cite{Zk1}. Bergen
and Passman\cite{BP} gave a certain characterization of graded
divergence-free Lie algebras. Kac also investigated suppersymmetric
graded generalizations of them \cite{Ka2,Ka3}. Motivated from his
classification of quadratic conformal algebras corresponding to
certain Hamiltonian pairs in integrable systems (cf. \cite{X3,X4}),
Xu\cite{X2} introduced nongraded generalizations of the Lie algebras
of Cartan type and determined their simplicity. An important feature
of Xu's algebras is that they do not contain any toral Cartan
subalgebras in generic case. The aim of this work is to study
certain representations of Xu's nongraded divergence-free Lie
algebras, whose isomorphism classes were determined by Su and Xu
\cite{SX}.

The representations of the Lie algebras of Cartan type were  vastly
studied.  Shen \cite{Sg1, Sg2, Sg3} studied mixed product of graded
modules over graded Lie algebras of Cartan type, and obtained
certain irreducible modules over a field with characteristic $p$.
Penkov and Serganova \cite{PS} gave an explicit description of the
support of an arbitrary irreducible weight module of the classical
Witt algebra, as well as its subalgebra with constant divergence.
 Rao \cite{R1,
R2} investigated the irreducibility of the weight modules over the
derivation Lie algebra of the algebra of Laurent polynomials
virtually constructed by Shen \cite{Sg1}. Zhao \cite{Zy1} determined
the module structure of Shen's mixed product over Xu's nongraded Lie
algebras of Witt type (cf. \cite{X2}). Moreover, she \cite{Zy3}
obtained a composition series for a family of modules with
parameters over Xu's nongraded Hamiltonian Lie algebras (cf.
\cite{X2}). These modules are constructed from finite-dimensional
multiplicity-free irreducible modules of symplectic Lie algebras by
Shen's mixed product.

Starting from Kac's conjecture \cite{Ka4,KR} on irreducible representations over the Virasoro algebra (which can also be viewed as a central extension of the rank one classical Witt algebra), Kaplansky \cite{Kap,KS} and Santharoubane \cite{KS} gave classifications of multiplicity-one
representations over classical Virasoro algebras. Later, Mathieu \cite{M} classified the Harish-Chandra modules over the Virasoro algebras, which confirms Kac's conjecture. Moreover, Su \cite{S1,S2} generalized Kaplansky and Santharoubane's result
to multiplicity-one modules over high rank Virasoro algebras and super-Virasoro algebras. Based on their classifications of multiplicity-one
representations over generalized Virasoro algebras in \cite{SZk},
Zhao \cite{Zk2} classified the multiplicity-one representations
over graded generalized Witt algebras. Su and Zhou
\cite{SZh} then generalized Zhao's result \cite{Zk2} to generalized
weight modules over the nongraded Witt algebras
introduced by Xu \cite{X2}.

Motivated by Kaplansky and Santharoubane's works in \cite{Kap,KS},
we \cite{C} gave a complete classification on multiplicity-one
representations over graded generalized divergence-free Lie algebras
introduced by Dokovic and Zhao \cite{DZ3}. Based on our previous
result,  we give in this paper a complete classification of the
irreducible and indecomposable generalized weight modules over
the nongraded divergence-free Lie algebras introduced by Xu
\cite{X2} with weight multiplicities less than or equal to one.

Throughout this paper, we let $\mbb{F}$ be an algebraically closed
field of characteristic zero. All the vector spaces are assumed over
$\mbb{F}$. Denote by $\mbb{N}$ the set of nonnegative integers
$\{0,1,2,\cdots\}$ and by $\mbb{Z}_+$ the set of positive integers.
For $m,n\in\mbb{Z}$, we shall use the notation
\begin{equation}
\overline{m,n}=\{m,m+1,m+2,\ldots,n\} \textrm{ if } m\leq n;\quad \overline{m,n}=\emptyset \textrm{ if }m>n.
\end{equation}
For any positive integer $n$, an additive subgroup $G$ of
$\mbb{F}^{n}$ is called {\it nondegenerate} if $G$ contains an
$\mbb{F}$-basis of $\mbb{F}^{n}$. Let $l_{1}$, $l_{2}$ and $l_{3}$
be three nonnegative integers such that \begin{equation}
l=l_{1}+l_{2}+l_{3}>0.
\end{equation}
Take any nondegenerate additive subgroup $\G$ of $\mbb{F}^{l_{2}+l_{3}}$
and let $\G=\{0\}$ when $l_{2}+l_{3}=0$. Let $\mcr{A}(l_{1},l_{2},l_{3};\G)$
be the semi-group algebra $\mbb{F}[\G\times \mbb{N}^{l_1+l_2}]$ with
the basis \begin{equation}
\{x^{\al}t^{\textbf{i}}\mid\al\in\G,\textbf{i}\in \mbb{N}^{l_1+l_2}\}.
\end{equation}
 Define a commutative associative algebraic operation ``$\cdot$'' on
$\mcr{A}(l_{1},l_{2},l_{3};\G)$ by
\begin{equation}
x^{\al}t^{\textbf{i}}\cdot x^{\be}t^{\textbf{j}}=x^{\al+\be}t^{\textbf{i}+\textbf{j}} \quad \textrm{for }\al,\be\in\G, \ \textbf{i},\textbf{j}\in \mbb{N}^{l_1+l_2}.\end{equation}
 Note that $x^{0}t^{\textbf{0}}$ is the identity element, which is denoted
by $1$ for convenience. Moreover, we  write $t^{\textbf{i}}$ and
$x^{\al}$ instead of $x^{0}t^{\textbf{i}}$ and
$x^{\al}t^{\textbf{0}}$ for short. When the context is clear, we
omit the notation ``$\cdot$'' in any associative algebra product.

For $k\in\mbb{Z}_+$ and $p\in\ove{1,l}$, we denote
\begin{equation}
k_{[p]}=(0,\cdots,0,\stackrel{p}{k},0,\cdots,0)\in\mbb{N}^{l}.
\end{equation}
For convenience, $\al\in\mbb{F}^{l_{2}+l_{3}}$ is written as
\begin{equation}
\al=(\al_1,\al_2,\cdots,\al_l)\ \textrm{ with }\al_1=\al_2=\cdots=\al_{l_1}=0,
\end{equation}
and $\textbf{i}\in\mbb{N}^{l_1+l_2}$ is written as
\begin{equation}
\textbf{i}=(i_1,i_2,\cdots,i_{l})\ \textrm{ with }i_{l_1+l_2+1}=i_{l_1+l_2+2}=\cdots=i_{l}=0.
\end{equation}
Moreover, we denote
\begin{equation}
|\textbf{i}|=\sum_{p=1}^{l_1+l_2}i_p \quad \textrm{ for }\textbf{i}\in\mbb{N}^{l_1+l_2}.
\end{equation}
We define linear transformations $\{\ptl_1,\ptl_2,\ldots,\ptl_l\}$
 on $\mcr{A}(l_{1},l_{2},l_{3};\G)$ by
 \begin{equation}\label{0.1}
\ptl_p( x^{\al}t^{\textbf{i}})=\al_{p} x^{\al}t^{\textbf{i}}+i_px^{\al}t^{\textbf{i}-1_{[p]}} \ \ \textrm{ for } p\in\ove{1,l}, \; \textbf{i}\in\mbb{N}^{l_1+l_2},\; \al\in\G,
\end{equation}
where $i_p=0$ and the monomial $i_px^{\al}t^{\textbf{i}-1_{[p]}}$ is treated as zero when $\textbf{i}-1_{[p]}\not\in \mbb{N}^{l_1+l_2}$.
Then $\{\ptl_{p}\mid p\in\overline{1,l}\}$
are mutually commutative derivations of $\mcr{A}(l_{1},l_{2},l_{3};\G)$, and they are $\mbb{F}$-linearly
independent.
Let
\begin{equation}
\mcr{D}=\sum_{i=1}^{l}\mbb{F}\ptl_{i}
\end{equation}
 and
 \begin{equation}
\mcr{W}(l_{1},l_{2},l_{3};\G)=\mcr{A}(l_{1},l_{2},l_{3};\G)\mcr{D}.
\end{equation}
 Then $\mcr{W}(l_{1},l_{2},l_{3};\G)$ is the simple Lie algebra of
nongraded Witt type constructed by Xu {\cite{X2}}. Its Lie bracket is given by
\begin{equation}
[\sum_{p=1}^{l}u_{p}\ptl_{p},\sum_{q=1}^{l}v_{q}\ptl_{q}]
=\sum_{p,q=1}^{l}(u_{p}\ptl_{p}(v_{q})-v_{p}\ptl_{p}(u_{q}))\ptl_{q} \quad \textrm{ for } u_{p}, v_{p} \in\mcr{A}(l_{1},l_{2},l_{3};\G).
\end{equation}

Define the divergence on $\mcr{W}(l_{1},l_{2},l_{3};\G)$ by
\begin{equation}\label{0.3}
\textrm{div}\ptl=\sum_{i=1}^{l}\ptl_{i}(u_{i})\quad\textrm{for }\ptl=\sum_{i=1}^{l}u_{i}\ptl_{i}\in\mcr{W}(l_{1},l_{2},l_{3};\G).
\end{equation}
Let $\rho\in\G$ be any element. We define
\begin{equation}\label{0.5}
D_{p,q}(u)=x^{\rho}(\ptl_p(x^{-\rho}u)\ptl_q-\ptl_q(x^{-\rho}u)\ptl_p)
\end{equation}
for $p,q\in\overline{1,l}$ and $u\in\mcr{A}(l_1,l_2,l_3;\G)$. Then
\begin{equation}\label{0.6}
\mcr{S}(l_1,l_2,l_3;\rho,\G)=\textrm{Span}\{D_{p,q}(u)\mid p,q\in\overline{1,l},\;u\in\mcr{A}(l_1,l_2,l_3;\G)\}
\end{equation}
forms a Lie subalgebra of the Lie algebra $\mcr{W}(l_1,l_2,l_3;\G)$. The Lie algebra $\mcr{S}(l_1,l_2,l_3;\rho,\G)$ was
first introduced by Xu \cite{X2}, as a nongraded generalization of graded simple Lie algebras of Special type, and further
 studied by Su and Xu \cite{SX}. It was called a {\it divergence-free Lie algebra} \cite{SX} in the sense that
$\mcr{S}(l_1,l_2,l_3;0,\G)$ consists of the divergence-free elements of
$\mcr{W}(l_1,l_2,l_3;\G)$. When $l_2>0$, the algebra
$\mcr{S}(l_1,l_2,l_3;\rho,\G)$ does not contain any toral Cartan
subalgebra. The special case $\mcr{S}(0,0,l_3;\rho,\G)$ was
introduced and studied by Dokovic and Zhao \cite{DZ3}. In \cite{C},
we classified the multiplicity-one representations of the Lie
algebras $\mcr{S}(0,0,l_3;0,\G)$ with $l_3\geq 3$. In this paper, we
further study the generalized multiplicity-one representations of
the Lie algebras $\mcr{S}(l_{1},l_{2},l_{3},0;\G)$ with
$l_{1}+l_{2}\geq 3$ and $l_{2}+l_{3}\geq3$.

A linear transformation $T$ of a vector space $V$ is called {\it locally finite} if
\begin{equation}
\textrm{dim}(\textrm{span}_{\mbb{F}}\{T^{n}(v)\mid n\in\mbb{Z}_{+}\})<\infty,\quad \forall\; v\in V;
\end{equation}
it is called {\it  semi-simple}, if for any $v\in V$, $\textrm{span}_{\mbb{F}}\{T^{n}(v)\mid n\in\mbb{Z}_{+}\}$ has a
$\mbb{F}$-basis consisting of eigenvectors of $T$; it is called {\it locally nilpotent}, if for any $v\in V$,
 there exists $n\in\mbb{Z}_{+}$ (may depend on $v$) such that $T^{n}(v)=0$. A subspace $U$ of $\textrm{End}V$ is called
 {\it locally finite} or {\it locally nilpotent}, if each element of $U$ is locally finite or locally nilpotent, respectively.

For convenience, we define
\begin{equation}\label{4.12}
<\ptl,\be>=\ptl(\be)=\be(\ptl)=\sum_{p=l_1+1}^l a_p\be_p
\end{equation}
for $\ptl=\sum_{p=1}^l a_p\ptl_p\in\mcr{D}$ and $\be=(\be_1,\be_2,\cdots,\be_l)\in\mbb{F}^{l_2+l_3}$. Let
\begin{equation}
\mcr{D}_1=\sum_{p=1}^{l_1}\mbb{F}\ptl_{p},\ \ \ \mcr{D}_2=\sum_{p=l_1+1}^{l_1+l_2}\mbb{F}\ptl_{p},\ \ \ \mcr{D}_3=\sum_{p=l_1+l_2+1}^{l}\mbb{F}\ptl_{p}.
\end{equation}
Observe that $ad\mcr{D}_1$ is locally nilpotent and $ad\mcr{D}$ is locally finite on $\mcr{S}$. Motivated from this and Su and Zhou's work \cite{SZh},
 we are interested in considering $\mcr{S}$-modules $V$ on which $\mcr{D}_1$ and $\mcr{D}$ act locally nilpotently and locally finitely, respectively.
  For such a module $V$, it has the decomposition
\begin{equation}\label{4.1}
V=\bigoplus_{\be\in\mbb{F}^{l_2+l_3}}V_{\be},
\end{equation}
where
\begin{equation}
V_{\be}=\{v\in V\mid (\ptl-\be(\ptl))^n (v)=0\textrm{ for }\ptl\in\mcr{D}\textrm{ and some } n\in\mbb{Z}_{+}\}.
\end{equation}
 We set
\begin{equation}\label{2.2}
V_{\be}^{(n)}=\{v\in V\mid (\ptl-\be(\ptl))^{n+1} (v)=0\textrm{ for }\ptl\in\mcr{D}\},\  \be\in\mbb{F}^{l_2+l_3},\ n\geq 0.
\end{equation}
Clearly, $V_{\be}=\bigcup_{n\in\mbb{N}} V_{\be}^{(n)}$ and $
V_{\be}\not=\{0\} \Longleftrightarrow V_{\be}^{(0)}\not=\{0\}.$

A module $V$ of $\mcr{S}(l_{1},l_{2},l_{3},0;\G)$ is called a {\it
generalized weight module} if it has the decomposition (\ref{4.1}).
For $\be\in\mbb{F}^{l_2+l_3}$, if $V_{\be}\not=\{0\}$, then
$V_{\be}$ is called the {\it generalized weight space with weight
$\be$}, and $V_{\be}^{(0)}$ is called the {\it weight space with
weight $\be$}. Moreover, we call $\dim V_{\be}^{(0)}$ the {\it
weight multiplicity corresponding to weight $\be$}.

The aim of this paper is to classify all the irreducible and
indecomposable generalized weight modules of
$\mcr{S}(l_{1},l_{2},l_{3},0;\G)$ with weight multiplicities equal
to or less than one. Since our investigations are based on the
previous result \cite{C}, we  only consider the nongraded
divergence-free Lie algebra $\mcr{S}(l_{1},l_{2},l_{3},0;\G)$ with
$l_2+l_3\geq 3$. When $l_1+l_2=1$, the Lie algebra
$\mcr{S}(l_{1},l_{2},l_{3},0;\G)$ does not contain the whole vector
space $\mcr{D}$, which means we don't get the generalized weight
module. The case $l_1+l_2=2$ will be studied later due to a
technical difficulty. The main result of this paper is the following
theorem:
\begin{theorem}\label{t:1.3}
Suppose $l_1+l_2\geq 3$ and $l_2+l_3\geq 3$. Assume that $V =\bigoplus_{\be\in\mbb{F}^{l_2+l_3}}V_{\be}$ is a
generalized weight module of $\mcr{S}(l_1,l_2,l_3,0;\G)$ with $\dim V_{\be}^{(0)}\leq 1$ for all
$\be \in \mbb{F}^{l_2+l_3}$. If $V$ is irreducible or indecomposable, then $V$ is either isomorphic to the 1-dimensional trivial module $\mathbb{F}v_{0}$, or isomorphic to the module $A_{\mu}$ for appropriate $\mu\in \mbb{F}^{l_2+l_3}$, where $A_{\mu}$ is the module with $\mbb{F}$-basis $\{v_{\al,\textbf{i}}\mid \al\in\G,\; \textbf{i}\in\mbb{N}^{l_1+l_2}\}$ and the following action:
\begin{eqnarray}
& & D_{p,q}(x^{\al}t^{\textbf{i}}).v_{\be,\textbf{j}}\nonumber\\
& = & (x^{\al}t^{\textbf{i}}(\al_p\ptl_q-\al_q\ptl_p) + i_p x^{\al}t^{\textbf{i}-1_{[p]}}\ptl_q - i_q x^{\al}t^{\textbf{i}-1_{[q]}}\ptl_p).v_{\be,\textbf{j}}\nonumber\\
&= & (\al_p(\be_q+\mu_q)-\al_q(\be_p+\mu_p) )v_{\be+\al,\textbf{i}+\textbf{j}}
+(j_q\al_p  - i_q (\be_p+\mu_p))  v_{\be+\al,\textbf{i}+\textbf{j}-1_{[q]}}\nonumber\\
& &+(i_p (\be_q+\mu_q)- j_p\al_q  ) v_{\be+\al,\textbf{i}+\textbf{j}-1_{[p]}}
+(i_p j_q- i_q j_p) v_{\be+\al,\textbf{i}+\textbf{j}-1_{[p]}-1_{[q]}}
\end{eqnarray}
for $\al, \be\in\G$, $\textbf{i},\textbf{j}\in \mbb{N}^{l_1+l_2}$ and $p,q \in\ove{1,l}$ with $p\not=q$.
\end{theorem}

This paper is organized as follows. In Section 2, we give some
properties of the Lie algebra $\mcr{S}(l_{1},l_{2},l_{3},0;\G)$ and
the $\mcr{S}$-module $A_{\mu}$. Moreover, we review our earlier
classification theorem  on the multiplicity-one representations over
the graded divergence-free Lie algebras (cf. \cite{C}). In Sections
3 and 4, we prove the main theorem.

\section{Some properties}

In this section,  we would like to digress to give some properties of the Lie algebra $\mcr{S}(l_{1},l_{2},l_{3},0;\G)$ and the $\mcr{S}$-module $A_{\mu}$, which are useful for the proof of the main theorem. We always assume that $l_1+l_2\geq 3$ and $l_2+l_3\geq 3$ in the present and the following sections.

For $\mcr{S}(l_{1},l_{2},l_{3},0;\G)$-module $A_{\mu}$, it is easy to verify:

\begin{pro}
(1) The module $A_{\mu}$ is irreducible if and only if $\mu\not\in\G$.

(2) When $\mu\in\G$, $A_{\mu}\simeq A_{0}$ is indecomposable, and it has the trivial submodule $\mbb{F}v_{-\mu,\textbf{0}}$. The quotient module $A_{\mu}'=A_{\mu}/\mbb{F}v_{-\mu,\textbf{0}}$ is simple.

(3) The possible isomorphisms between modules of type $A_{\mu}$ are:
$A_{\mu}\simeq A_{\eta}$ iff $\mu-\eta\in\G$.
\end{pro}

We would like to remark that, the simple quotient module
$V=A_{0}'=A_0/\mbb{F}v_{0,\textbf{0}}$ is a generalized weight module
of $\mcr{S}(l_1,l_2,l_3,0;\G)$ with weight set $\G$; and we have $\dim V_{0}^{(0)}=l_1+l_2$,
$\dim V_{\be}^{(0)}= 1$ for all $\be \in \G\backslash\{0\}$. 
The
following conclusion will be frequently used later:

\begin{lemma}\label{le:1.6}
 Let $\mathbf{T}$ be a linear transformation on a vector space $U$,
and let $U_1$ be a subspace of $U$ such that $\mathbf{T}
(U_1)\subset U_1$. Suppose that $u_1,\ u_2,\ \cdots,\ u_n$ are
eigenvectors of $\mathbf{T}$ corresponding to different eigenvalues.
If $\sum_{j=1}^n u_j \in U_1$, then $u_1,\ u_2,\ \cdots,\ u_n\ \in
U_1$.
\end{lemma}

\begin{pro}\label{t:2.1}
The Lie algebra $\mcr{S}(l_{1},l_{2},l_{3},0;\G)$ is generated by
\begin{equation}\label{1.1}
\{D_{p,q}(x^{\al}),D_{p,q}(t^{\textbf{j}})\mid p,q\in\overline{1,l},\; \al\in\G\backslash\{0\}, \; \textbf{j}\in\mbb{N}^{l_1+l_2}\textrm{ with }|\textbf{j}|\leq 2 \}.
\end{equation}
\end{pro}

\noindent{\bf Proof.} Denote by $\mcr{K}$ the Lie subalgebra of $\mcr{S}(l_{1},l_{2},l_{3},0;\G)$ which is generated by the set (\ref{1.1}). It suffices to prove $\mcr{K}=\mcr{S}(l_{1},l_{2},l_{3},0;\G)$. By
(\ref{0.1}) and (\ref{0.5}), we have
\begin{equation}\label{1.2}
D_{p,q}(x^{\al}t^{\textbf{i}})=x^{\al}t^{\textbf{i}}(\al_p\ptl_q-\al_q\ptl_p)
+i_p x^{\al}t^{\textbf{i}-1_{[p]}}\ptl_q -i_q x^{\al}t^{\textbf{i}-1_{[q]}}\ptl_p
\end{equation}
for $p,q\in\overline{1,l}$, $\al\in\G$ and $\textbf{i}\in\mbb{N}^{l_1+l_2}$.
Then we proceed our proof in two steps.\vspace{0.2cm}

{\it Step 1}. $D_{p,q}(x^{\al}t^{\textbf{i}})\in \mcr{K}$ for any $p,q\in\overline{1,l}$, $\al\in\G\backslash\{0\}$ and $\textbf{i}\in\mbb{N}^{l_1+l_2}$.\vspace{0.2cm}

We prove this step by induction on $|\textbf{i}|$. By (\ref{1.1}), we know that this holds when $|\textbf{i}|=0$. Assume that
\begin{equation}\label{0.8}
D_{p,q}(x^{\al}t^{\textbf{i}})\in \mcr{K} \textrm{ for any } p,q\in\overline{1,l},\; \al\in\G\backslash\{0\},\; \textbf{i}\in\mbb{N}^{l_1+l_2}  \textrm{ with }|\textbf{i}|\leq k,
\end{equation}
where $k\geq 0$. It suffices to prove 
$D_{p,q}(x^{\al}t^{\textbf{i}+1_{[r]}})\in \mcr{K}$ for all
$p,q\in\overline{1,l}$, $\al\in\G\backslash\{0\}$,
$r\in\overline{1,l_1+l_2}$ and $\textbf{i}\in\mbb{N}^{l_1+l_2}$ with
$|\textbf{i}|= k$.

Fix any $\al\in\G\backslash\{0\}$ and $\textbf{i}\in\mbb{N}^{l_1+l_2}$ with $|\textbf{i}|= k$. Choose $s\in\ove{1,l}$ such that $\al_s\not=0$. Then for any $p,q\in\overline{1,l}$ and $r\in\overline{1,l_1+l_2}\backslash\{s\}$, we have
\begin{eqnarray}
&&[D_{p,q}(x^{\al}t^{\textbf{i}}),D_{r,s}(t^{2_{[r]}})]\nonumber\\
&=&2[x^{\al}t^{\textbf{i}}(\al_p\ptl_q-\al_q\ptl_p)+i_p x^{\al}t^{\textbf{i}-1_{[p]}}\ptl_q-i_q x^{\al}t^{\textbf{i}-1_{[q]}}\ptl_p,t^{1_{[r]}}\ptl_s]\nonumber\\
&=&2(\delta_{p,r}D_{s,q}(x^{\al}t^{\textbf{i}})-\delta_{q,r}D_{s,p}(x^{\al}t^{\textbf{i}})
-i_sD_{p,q}(x^{\al}t^{\textbf{i}+1_{[r]}-1_{[s]}})-\al_sD_{p,q}(x^{\al}t^{\textbf{i}+1_{[r]}})).\label{0.7}
\end{eqnarray}
Thus by (\ref{0.7}) and the induction hypothesis (\ref{0.8}), we
find
\begin{equation}\label{0.9}
D_{p,q}(x^{\al}t^{\textbf{i}+1_{[r]}})\in \mcr{K} \ \textrm{ for }p,q\in\overline{1,l},\; r\in\overline{1,l_1+l_2}\backslash\{s\}.
 \end{equation}
 If $s\not\in\overline{1,l_1+l_2}$, then it is done. Consider the case that $s\in\overline{1,l_1+l_2}$.
 If $\al_{s'}\not=0$ for some $s'\in\overline{1,l}\backslash\{s\}$, then with $s$ replaced by $s'$ in the above arguments,
 we can get $D_{p,q}(x^{\al}t^{\textbf{i}+1_{[s]}})\in \mcr{K}$ for $p,q\in\overline{1,l}$. Otherwise, we have $\al_{r}=0$
 for all $r\in\overline{1,l}\backslash\{s\}$. Pick $r\in\overline{1,l_1+l_2}\backslash\{s\}$. Expressions (\ref{1.1}) and (\ref{0.9}) give
\begin{eqnarray}
\mcr{K}&\ni&[D_{p,q}(x^{\al}t^{\textbf{i}+1_{[r]}}),D_{s,r}(t^{2_{[s]}})]\nonumber\\
&=&-2((i_r+1)D_{p,q}(x^{\al}t^{\textbf{i}+1_{[s]}})+\delta_{p,s}D_{q,r}(x^{\al}t^{\textbf{i}+1_{[r]}})
-\delta_{q,s}D_{p,r}(x^{\al}t^{\textbf{i}+1_{[r]}}))
\end{eqnarray}
for $p,q\in\overline{1,l}$. So this and (\ref{0.9}) indicate
\begin{equation}\label{0.15}
D_{p,q}(x^{\al}t^{\textbf{i}+1_{[s]}})\in \mcr{K} \  \textrm{ for }p,q\in\overline{1,l}.
\end{equation}
To summarize, (\ref{0.9}) and (\ref{0.15}) show
\begin{equation}\label{1.3}
D_{p,q}(x^{\al}t^{\textbf{i}+1_{[r]}})\in \mcr{K}
\end{equation}
 for $p,q\in\overline{1,l}$, $\al\in\G\backslash\{0\}$, $r\in\overline{1,l_1+l_2}$ and $\textbf{i}\in\mbb{N}^{l_1+l_2}$ with $|\textbf{i}|= k$, which completes the proof of the step.\vspace{0.2cm}

{\it Step 2. $D_{p,q}(t^{\textbf{i}})\in \mcr{K}$ for any $p,q\in\overline{1,l}$ and  $\textbf{i}\in\mbb{N}^{l_1+l_2}$.}\vspace{0.2cm}

We shall prove this step by induction on $|\textbf{i}|$. When $|\textbf{i}|=0$, $D_{p,q}(t^{\textbf{i}})=0$ for any $p,q\in\overline{1,l}$. Assume that
\begin{equation}\label{0.4}
D_{p,q}(t^{\textbf{j}})\in \mcr{K} \textrm{ for any }p,q\in\overline{1,l},\; \textbf{j}\in\mbb{N}^{l_1+l_2} \textrm{ with } |\textbf{j}|\leq k,
 \end{equation}
 where $k \geq 0$. It suffices to prove $D_{p,q}(t^{\textbf{i}})\in \mcr{K}$ for all $p,q\in\overline{1,l}$ and $\textbf{i}\in\mbb{N}^{l_1+l_2}$ with $|\textbf{i}|= k+1$. Suppose $p\not=q$ and $|\textbf{i}|= k+1$ in the following. We give the proof in several cases.\vspace{0.2cm}

Case 1. $p,q\in\overline{l_1+l_2+1,l}$.\vspace{0.2cm}

In this case, $D_{p,q}(t^{\textbf{i}})=0$.\vspace{0.2cm}

Case 2. $p,q\in\overline{1,l_1+l_2}$.\vspace{0.2cm}

Pick $r\in\ove{l_1+1,l}\backslash\{p,q\}$. Choose $\al\in\G\backslash\{0\}$ such that $\al_r\not=0$. Then (\ref{1.2}) and Step 1 indicate
\begin{eqnarray}
\mcr{K} &\ni & [D_{r,q}(x^{\al}t^{\textbf{i}}),D_{r,p}(x^{-\al})]\nonumber\\
&=&[x^{\al}t^{\textbf{i}}(\al_r\ptl_q-\al_q\ptl_r)+i_r x^{\al}t^{\textbf{i}-1_{[r]}}\ptl_q-i_q x^{\al}t^{\textbf{i}-1_{[q]}}\ptl_r,x^{-\al}(\al_p\ptl_r-\al_r\ptl_p)]\nonumber\\
&=&\al_r^2 D_{p,q}(t^{\textbf{i}})+\al_r \al_p D_{q,r}(t^{\textbf{i}})+\al_r \al_q D_{r,p}(t^{\textbf{i}})+\al_p i_r D_{q,r}(t^{\textbf{i}-1_{[r]}})\nonumber\\
& &
+\al_r i_p D_{r,q}(t^{\textbf{i}-1_{[p]}}).
\end{eqnarray}
By the induction hypothesis (\ref{0.4}), we get
\begin{eqnarray}
\al_r^2 D_{p,q}(t^{\textbf{i}})+\al_r \al_p D_{q,r}(t^{\textbf{i}})+\al_r \al_q D_{r,p}(t^{\textbf{i}})\in\mcr{K}.
\end{eqnarray}
Since
\begin{equation}
[D_{p,q}(t^{1_{[p]}+1_{[q]}}),D_{p,q}(t^{\textbf{i}})]=(i_q-i_p)D_{p,q}(t^{\textbf{i}}),
\end{equation}
\begin{equation}
[D_{p,q}(t^{1_{[p]}+1_{[q]}}),D_{q,r}(t^{\textbf{i}})]=(i_q-i_p-1)D_{q,r}(t^{\textbf{i}}),
\end{equation}
\begin{equation}
[D_{p,q}(t^{1_{[p]}+1_{[q]}}),D_{r,p}(t^{\textbf{i}})]=(i_q-i_p+1)D_{r,p}(t^{\textbf{i}}),
\end{equation}
we have $D_{p,q}(t^{\textbf{i}})\in\mcr{K}$ by (\ref{1.1}) and Lemma
\ref{le:1.6}.\vspace{0.2cm}

Case 3. $p\in\overline{1,l_1+l_2}$, $q\in\overline{l_1+l_2+1,l}$.\vspace{0.2cm}

Pick $s\in\overline{1,l_1+l_2}\backslash\{p\}$. Case 2 tells that
$D_{p,s}(t^{\textbf{i}})\in\mcr{K}$. So by (\ref{1.1}) and
(\ref{1.2}), we have
\begin{equation}
\mcr{K}\ni [D_{p,s}(t^{\textbf{i}}),\frac{1}{2}D_{s,q}(t^{2_{[s]}})]=i_pt^{\textbf{i}-1_{[p]}}\ptl_q=D_{p,q}(t^{\textbf{i}}).
\end{equation}

\vspace{0.2cm}

Case 4. $p\in\overline{l_1+l_2+1,l}$, $q\in\overline{1,l_1+l_2}$.\vspace{0.2cm}

We have $D_{p,q}(t^{\textbf{i}})=-D_{q,p}(t^{\textbf{i}})\in\mcr{K}$ by Case 3.\vspace{0.2cm}

So this step holds.\vspace{0.2cm}

Thus, from (\ref{0.6}) we see that $\mcr{K}=\mcr{S}(l_{1},l_{2},l_{3},0;\G)$. So this proposition holds. $\qquad\Box$\vspace{0.3cm}

We also have:

\begin{coral}\label{c:2.2}
The set
\begin{equation}\label{1.4}
\{D_{p,q}(x^{\al}t^{\textbf{i}})\mid p,q\in\overline{1,l},\; \al\in\G\backslash\{0\}, \textbf{i}\in\mbb{N}^{l_1+l_2}\}
\end{equation}
generates the Lie algebra $\mcr{S}(l_{1},l_{2},l_{3},0;\G)$.
\end{coral}

\noindent{\bf Proof.} Denote by $\mcr{K}'$ the Lie subalgebra of
$\mcr{S}(l_{1},l_{2},l_{3},0;\G)$ which is generated by the set
(\ref{1.4}). By Proposition \ref{t:2.1}, it suffices to prove
$D_{p,q}(t^{\textbf{i}})\in \mcr{K}'$ for all $p,q\in\overline{1,l}$
and  $\textbf{i}\in\mbb{N}^{l_1+l_2}$ with $|\textbf{i}|\leq 2$.
Notice that $\textrm{Span}_{\mbb{F}}\{D_{p,q}(t^{1_{[r]}})\mid
p,q\in\overline{1,l},\; r\in\overline{1,l_1+l_2}\}=\mcr{D}$. So
first,  we need to prove $\mcr{D} \subseteq \mcr{K}'$.

 Fix some $r\in\ove{l_1+1,l}$. Choose $\al\in\G$ such that $\al_r\not=0$. Pick $p\in\ove{1,l_1+l_2}\backslash\{r\}$. Then for any $q\in\ove{1,l}\backslash\{r,p\}$, we have
\begin{eqnarray}
\mcr{K}' &\ni & [D_{r,q}(x^{\al}t^{1_{[p]}}),D_{r,p}(x^{-\al})]\nonumber\\
&=&[x^{\al}t^{1_{[p]}}(\al_r\ptl_q-\al_q\ptl_r),x^{-\al}(\al_p\ptl_r-\al_r\ptl_p)]\nonumber\\
&=&\al_r(\al_r\ptl_q-\al_q\ptl_r)\label{1.5}
\end{eqnarray}
by (\ref{1.4}). Namely, we obtain
\begin{equation}
\al_r\ptl_q-\al_q\ptl_r\in\mcr{K}'\quad \textrm{ for }q\in\ove{1,l}\backslash\{r,p\}.
\end{equation}
 By picking another $p\in\ove{1,l_1+l_2}\backslash\{r\}$, we can get
\begin{equation}\label{1.6}
\al_r\ptl_q-\al_q\ptl_r\in\mcr{K}' \quad \textrm{ for }q\in\ove{1,l}\backslash\{r\}.
\end{equation}
Take $r'\in\ove{l_1+1,l}\backslash\{r\}$. Choose $\be\in\G$ such
that $\be_r\not=0$ and $\al_r\be_{r'}-\al_{r'}\be_r\not=0$. With
$\al$ replaced by $\be$ in (\ref{1.5})--(\ref{1.6}), we can get that
$\be_r\ptl_{r'}-\be_{r'}\ptl_r\in\mcr{K}'$. Since (\ref{1.6}) also
gives $\al_r\ptl_{r'}-\al_{r'}\ptl_r\in\mcr{K}'$, we have
\begin{equation}\label{1.7}
\ptl_r=\frac{1}{\al_r\be_{r'}-\al_{r'}\be_r}(\be_r(\al_r\ptl_{r'}-\al_{r'}\ptl_r)-\al_r(\be_r\ptl_{r'}-\be_{r'}\ptl_r))\in\mcr{K}'.
\end{equation}
Thus, from (\ref{1.6}) and (\ref{1.7}) we deduce 
\begin{equation}\label{0.13}
\mcr{D} \subseteq \mcr{K}'.
\end{equation}

 Second, we want to prove $t^{1_{[p]}}\ptl_q\in\mcr{K}'$ for any $p\in\ove{1,l_1+l_2}$ and $q\in\ove{1,l}\backslash\{p\}$.

Fix any $p\in\ove{1,l_1+l_2}$. Pick $r\in\ove{l_1+1,l}\backslash\{p\}$. Choose $\al\in\G$ such that $\al_r\not=0$. Then for any $q\in\ove{1,l}\backslash\{r,p\}$, we have
\begin{eqnarray}
\mcr{K}' &\ni & [D_{r,q}(x^{\al}t^{2_{[p]}}),D_{r,p}(x^{-\al})]\nonumber\\
&=&[x^{\al}t^{2_{[p]}}(\al_r\ptl_q-\al_q\ptl_r),x^{-\al}(\al_p\ptl_r-\al_r\ptl_p)]\nonumber\\
&=& 2 \al_r t^{1_{[p]}} (\al_r\ptl_q-\al_q\ptl_r)
\end{eqnarray}
by (\ref{1.4}). Namely,
 \begin{equation}\label{1.8}
 t^{1_{[p]}} (\al_r\ptl_q-\al_q\ptl_r)\in\mcr{K}'  \quad \textrm{ for }q\in\ove{1,l}\backslash\{r,p\}.
  \end{equation}
  Pick $r'\in\ove{l_1+1,l}\backslash\{r,p\}$. Choose $\be\in\G$ such that $\be_r\not=0$ and $\al_r\be_{r'}-\al_{r'}\be_r\not=0$. Likewise, we can get $t^{1_{[p]}} (\be_r\ptl_{r'}-\be_{r'}\ptl_r) \in\mcr{K}'$. Since (\ref{1.8}) gives $t^{1_{[p]}}(\al_r\ptl_{r'}-\al_{r'}\ptl_r)\in\mcr{K}'$, we have
\begin{equation}\label{1.9}
t^{1_{[p]}}\ptl_r=\frac{1}{\al_r\be_{r'}-\al_{r'}\be_r}(\be_r t^{1_{[p]}}(\al_r\ptl_{r'}-\al_{r'}\ptl_r)-\al_r t^{1_{[p]}}(\be_r\ptl_{r'}-\be_{r'}\ptl_r))\in\mcr{K}'.
\end{equation}
So (\ref{1.8}) and (\ref{1.9}) show 
\begin{equation}\label{0.11}
t^{1_{[p]}}\ptl_q\in\mcr{K}' \textrm{ for all } p\in\ove{1,l_1+l_2} \textrm{ and } q\in\ove{1,l}\backslash\{p\}.
\end{equation}

Third, we want to prove that $t^{1_{[q]}}\ptl_q-t^{1_{[p]}}\ptl_p \in\mcr{K}'$ for any $p,q\in\ove{1,l_1+l_2}$.

Fix any $p,q\in\ove{1,l_1+l_2}$. We further suppose that $p\not=q$. Pick $r\in\ove{1,l}\backslash\{p,q\}$. Choose $\al\in\G$ such that $\al_r\not=0$. Then (\ref{1.4}) indicates
\begin{eqnarray}
\mcr{K}' &\ni & [D_{r,q}(x^{\al}t^{1_{[p]}+1_{[q]}}),D_{r,p}(x^{-\al})]\nonumber\\
&=&[x^{\al}t^{1_{[p]}+1_{[q]}}(\al_r\ptl_q-\al_q\ptl_r)- x^{\al}t^{1_{[p]}}\ptl_r,x^{-\al}(\al_p\ptl_r-\al_r\ptl_p)]\nonumber\\
&=&\al_r^2 (t^{1_{[q]}}\ptl_q-t^{1_{[p]}}\ptl_p)+ \al_r \al_p t^{1_{[p]}}\ptl_r - \al_r \al_q t^{1_{[q]}}\ptl_r-\al_r\ptl_r.\label{0.12}
\end{eqnarray}
So by (\ref{0.13}), (\ref{0.11}) and (\ref{0.12}), we have
\begin{equation}\label{0.14}
t^{1_{[q]}}\ptl_q-t^{1_{[p]}}\ptl_p \in\mcr{K}' \quad \textrm{ for all } p,q\in\ove{1,l_1+l_2}.
\end{equation}

In summary, (\ref{0.13}), (\ref{0.11}) and (\ref{0.14}) give
\begin{equation}\label{0.10}
\{D_{p,q}(t^{\textbf{i}})\mid   p,q\in\overline{1,l},\; \textbf{i}\in\mbb{N}^{l_1+l_2} \textrm{ with } |\textbf{i}|\leq 2\} \subseteq \mcr{K}'.
\end{equation}
So the corollary follows from (\ref{1.4}), (\ref{0.10}) and Proposition \ref{t:2.1}.
$\qquad\Box$\vspace{0.3cm}

Last, for convenience of the reader, we shall review below the multiplicity-one representations over the graded divergence-free Lie algebras \cite{C}.

By (\ref{0.1}) and (\ref{4.12}), we have
\begin{equation}
\ptl(x^\al)=\ptl(\al)x^\al \ \textrm{ for }\al\in\G \textrm{ and }\ptl\in\mcr{D};
\end{equation}
\begin{equation}\label{0.19}
(\mcr{D}_2+\mcr{D}_3)(\al)=\{0\} \textrm{ with }\al\in\mbb{F}^{l_2+l_3}, \textrm{ iff }\al=0;
\end{equation}
\begin{equation}
\ptl (\G)=\{0\} \textrm{ with } \ptl \in \mcr{D}_2+\mcr{D}_3, \textrm{ iff }\ptl=0.
\end{equation}
Denote
\begin{equation}\label{2.3}
\textrm{ker} \be=\{\ptl\in\mcr{D} \mid \ptl(\be)=0\}\ \textrm{ for }\be\in\mbb{F}^{l_2+l_3}.
\end{equation}
Then the main classification theorem  in \cite{C} is:

\begin{pro}\label{t:1.2}
 Suppose $l\geq 3$. Assume that $V = \bigoplus_{\theta \in \G} V_\theta$ is a
$\G$-graded $\mcr{S}(0,0,l,0;\G)$-module with $\dim V_\theta\leq 1$ for
$\theta \in \G$. If $V$ is irreducible or indecomposable, then $V$ is isomorphic to one of the following modules
for appropriate $\mu\in \mbb{F}^l$ and $\eta \in \mbb{F}^l\backslash\{0\}$:
\begin{equation}
(\rmnum{1})  \textrm{the trivial module } \mathbb{F}v_{0},\ \
(\rmnum{2})    \mathcal{M}_{\mu}',\ \ (\rmnum{3})   \mathcal{A}_{\eta},\ \ (\rmnum{4}) \mathcal{B}_{\eta},
\end{equation}
where $\mathcal{M}_{\mu}$, $\mathcal{A}_{\eta }$ and $\mathcal{B}_{\eta }$ are $\mcr{S}(0,0,l,0;\G)$-modules with basis $\{v_\al\mid \al\in\G\}$ and the following actions:
\begin{eqnarray}
\mathcal{M}_{\mu}: \ \ & &x^\al\ptl . v_\be=\ptl(\be+\mu)v_{\al+\be}  \textrm{ for }
\al \in \G\backslash\{0\},  \be\in \G,  \ptl \in \kn\al;
\end{eqnarray}
\begin{eqnarray}
\mathcal{A}_{\eta }: \ \ & &x^\al\ptl. v_\be=  \ptl(\be)v_{\al+\be}  \textrm{ for
} \al, \be\in \G \backslash \{0\},  \ptl \in \kn\al, \nonumber \\
& & x^\al\ptl . v_{0}  = \ptl(\eta)v_{\al}  \textrm{
for } \al\in \G\backslash \{0\},   \ptl \in
\kn\al;
\end{eqnarray}
\begin{eqnarray}
\mathcal{B}_{\eta }:\ \ & &x^\al\ptl . v_\be = \ptl(\be)v_{\al+\be} \textrm{ for
} \al\in \G \backslash \{0\},  \be\in
\G\backslash\{ -\al\},   \ptl \in \kn\al, \nonumber \\
& & x^\al\ptl . v_{-\al}  = \ptl(\eta)v_{0}  \textrm{
for } \al\in \G \backslash \{0\},  \ptl \in
\kn\al,
\end{eqnarray}
and $\mathcal{M}_{\mu}'$ is denoted as the only nontrivial irreducible quotient module of $\mathcal{M}_{\mu}$.
\end{pro}

\section{Proof of the main theorem (\Rmnum{1})}

 Suppose that $V$ is a generalized weight module of $\mcr{S}(l_{1},l_{2},l_{3},0;\G)$ with $\dim V_{\be}^{(0)}\leq 1$ for all $\be\in\mbb{F}^{l_2+l_3}$. Let
 \begin{equation}
 V(\mu)=\bigoplus_{\al\in\G}V_{\al+\mu}\quad \textrm{ for }\mu\in\mbb{F}^{l_2+l_3}.
  \end{equation}
  Then $V(\mu)$ is a submodule of $V$, and $V$ is a direct sum of $V(\mu)$ for different $\mu\in\mbb{F}^{l_2+l_3}$. So, if $V$ is irreducible or indecomposable, we must have $V=V(\mu)$ for some $\mu\in\mbb{F}^{l_2+l_3}$.

  By Proposition \ref{t:2.1}, we only need to derive the action of the set (\ref{1.1}) on $V=V(\mu)$.
   Except the trivial case in Lemma \ref{le:3.1}, we should determine a basis of $V=V(\mu)$, and
   then derive the action of the set (\ref{1.1}) on the basis. For the nontrivial cases, we shall first
    give a general discussion for any $\mu\in\mbb{F}^{l_2+l_3}$, then prove the main theorem for
     $\mu\in\mbb{F}^{l_2+l_3}\backslash\G$ at the end of this section, and last prove the main theorem for $\mu\in\G$ in the next section.

\begin{lemma}\label{le:3.1}
If $V_{\sgm+\mu}^{(0)}\not=\{0\}$ for some $\sgm\in\G$, and $V_{\rho+\mu}^{(0)}=\{0\}$ for all $\rho\in\G\backslash\{\sgm\}$, then $\sgm+\mu=0$, and $V=V_0^{(0)}$ is a 1-dimensional trivial $\mcr{S}(l_{1},l_{2},l_{3},0;\G)$-module.
\end{lemma}

\noindent{\bf Proof.} By (3.1), we have $V=V_{\sgm+\mu}$. Therefore,
\begin{equation}
D_{p,q}(x^{\ga}t^{\textbf{i}}).V_{\sgm+\mu}=\{0\}\ \textrm{ for all } p,q\in\overline{1,l},\; \ga\in\G\backslash\{0\}, \textbf{i}\in\mbb{N}^{l_1+l_2}
\end{equation}
by (\ref{2.2}). Moreover, Corollary \ref{c:2.2} gives rise to
\begin{equation}
\mcr{S}(l_{1},l_{2},l_{3},0;\G).V_{\sgm+\mu}=\{0\}.
\end{equation}
So $V=V_{\sgm+\mu}$ is a direct sum of trivial submodules. Since $V$ is irreducible or indecomposable, we must have $V=V_{\sgm+\mu}^{(0)}$ as a trivial $\mcr{S}$-module.
As one result, we get $\ptl .V_{\sgm+\mu}^{(0)}=\{0\}$ for all $\ptl\in \mcr{D}_2+\mcr{D}_3$, which implies $\sgm+\mu=0$. So $V=V_0^{(0)}$, and it is a trivial $\mcr{S}(l_{1},l_{2},l_{3},0;\G)$-module. This completes the proof of the lemma. $\qquad\Box$

\vspace{0.3cm}

In the rest of the paper, we shall always assume that
\begin{equation}
V_{\sgm+\mu}^{(0)}\not=\{0\} \textrm{ and } V_{\rho+\mu}^{(0)}\not=\{0\} \textrm{ for some } \sgm,\rho\in\G \textrm{ with }\sgm\not=\rho.
\end{equation}

\begin{lemma}\label{le:3.2}
If $V_{\sgm+\mu}^{(0)}\not=\{0\}$ and $V_{\rho+\mu}^{(0)}\not=\{0\}$ for some $\sgm,\rho\in\G$ with $\sgm\not=\rho$, then $V_{\ga+\mu}^{(0)}\not=\{0\}$ for all $\ga\in\G\backslash\{-\mu\}$. Moreover, there exist nonzero $v_{\ga,\textbf{0}}\in V_{\ga+\mu}^{(0)}$ with $\ga\in\G\backslash\{-\mu\}$, such that
\begin{equation}
D_{p,q}(x^{\al}).v_{\ga,\textbf{0}}=x^{\al}(\al_p\ptl_q-\al_q\ptl_p).v_{\ga,\textbf{0}}
=(\al_p(\ga_q+\mu_q)-\al_q(\ga_p+\mu_p))v_{\al+\ga,\textbf{0}}
\end{equation}
for $p,q\in\overline{1,l}$, $\al\in\G\backslash\{0\}$ and $\ga\in\G\backslash\{-\mu\}$.
\end{lemma}

\noindent{\bf Proof.} Since
\begin{equation}\label{1.10}
\mcr{S}_0:=\textrm{Span}_{\mbb{F}}\{D_{p,q}(x^{\al})\mid p,q\in\overline{l_1+1,l},\; \al\in\G\backslash\{0\}\}\simeq \mcr{S}(0,0,l_{2}+l_{3},0;\G),
\end{equation}
we can see $V^{(0)}=\bigoplus_{\al\in\G}V_{\al+\mu}^{(0)}$ as a $\G$-graded $\mcr{S}_0$-module with $\textrm{dim}V^{(0)}\geq 2$ and $\textrm{dim}V_{\al+\mu}^{(0)}\leq 1$ for all $\al\in\G$.
So first of all, we need to prove:\vspace{0.2cm}

{\it Claim 1.} $\mcr{S}_0.V^{(0)}\not=\{0\}$.\vspace{0.2cm}

Suppose $\mcr{S}_0.V^{(0)}=\{0\}$, namely, $V^{(0)}$ is a direct sum of some trivial $\mcr{S}_0$-submodules, then we will get a contradiction.

Since $V_{\sgm+\mu}^{(0)}\not=\{0\}$ and $V_{\rho+\mu}^{(0)}\not=\{0\}$ for some $\sgm,\rho\in\G$ with $\sgm\not=\rho$, we have $V_{\theta+\mu}^{(0)}\not=\{0\}$ for some $\theta\in\G\backslash\{-\mu\}$. Fix  $\theta\in\G\backslash\{-\mu\}$ satisfying $V_{\theta+\mu}^{(0)}\not=\{0\}$.

Let $\bar{\al}=\al+\mu$ for all $\al\in\G$. Since $\bar{\theta}\not=0$, we can choose $r\in\ove{l_1+1,l}$ such that $\bar{\theta}_r\not=0$. Pick $p\in\ove{1,l_1+l_2}\backslash\{r\}$ and $s\in\ove{l_1+1,l}\backslash\{r,p\}$. Choose $\al\in\G$ such that $\al_r\not=0$ and $\bar{\theta}_r \al_s-\bar{\theta}_s\al_r\not=0$. Moreover, we let
\begin{equation}
\ptl=\bar{\theta}_r \ptl_s-\bar{\theta}_s\ptl_r,\ \ \ptl'=\al_r\ptl_s-\al_s\ptl_r,\ \ \ptl''=(\bar{\theta}_r-\al_r) \ptl_s - (\bar{\theta}_s-\al_s)\ptl_r.
\end{equation}
From (\ref{2.2}) it can be deduced that
\begin{equation}\label{0.16}
t^{1_{[p]}}\ptl.V_{\theta+\mu}^{(0)}\subseteq V_{\theta+\mu}^{(0)},\ \ t^{1_{[p]}}\ptl''.V_{\theta-\al+\mu}^{(0)}\subseteq V_{\theta-\al+\mu}^{(0)}.
\end{equation}
Take $0\not=v_\theta \in V_{\theta+\mu}^{(0)}$. So on one hand, we have
\begin{equation}\label{1.11}
\ptl'.v_{\theta}=\ptl'(\bar{\theta})v_{\theta}=-(\bar{\theta}_r \al_s-\bar{\theta}_s\al_r)v_{\theta}\not=0
\end{equation}
by (\ref{2.2}).
On the other hand, (\ref{0.16}) and $\mcr{S}_0.V^{(0)}=\{0\}$ give
\begin{eqnarray}
& &x^{\al}t^{1_{[p]}}\ptl'.v_{\theta}\nonumber\\
&=&-\frac{1}{\ptl(\al)}[x^\al\ptl', t^{1_{[p]}}\ptl].v_{\theta}\nonumber\\
&=& -\frac{1}{\bar{\theta}_r \al_s-\bar{\theta}_s\al_r}(x^\al\ptl'. (t^{1_{[p]}}\ptl.v_{\theta}) - t^{1_{[p]}}\ptl. x^\al\ptl'.v_{\theta})\nonumber\\
&=&0
\end{eqnarray}
and
\begin{eqnarray}
x^{\al}t^{1_{[p]}}\ptl'.V_{\theta-\al+\mu}^{(0)}&=&-\frac{1}{\bar{\theta}_r \al_s-\bar{\theta}_s\al_r}[x^\al\ptl', t^{1_{[p]}}\ptl''].V_{\theta-\al+\mu}^{(0)}\nonumber\\
&=&\{0\},
\end{eqnarray}
which further implies
\begin{eqnarray}
& &\ptl'.v_{\theta}\nonumber\\
&=&-\frac{1}{\al_r}[x^{\al}t^{1_{[p]}}\ptl', x^{-\al}(\al_r\ptl_p-\al_p\ptl_r)].v_\theta\nonumber\\
&=& -\frac{1}{\al_r}(x^{\al}t^{1_{[p]}}\ptl'. (x^{-\al}(\al_r\ptl_p-\al_p\ptl_r).v_{\theta}) -  x^{-\al}(\al_r\ptl_p-\al_p\ptl_r). (x^{\al}t^{1_{[p]}}\ptl'.v_{\theta}))\nonumber\\
&=&0,\label{0.17}
\end{eqnarray}
where in the third line of (\ref{0.17}), we have $x^{-\al}(\al_r\ptl_p-\al_p\ptl_r).v_{\theta}\in V_{\theta-\al+\mu}^{(0)}$. Since (\ref{0.17}) contradicts (\ref{1.11}), we must have $\mcr{S}_0.V^{(0)}\not=\{0\}$.
Thus this claim holds.\vspace{0.4cm}

Note that $V^{(0)}=\bigoplus_{\al\in\G}V_{\al+\mu}^{(0)}$ is a $\G$-graded $\mcr{S}_0$-submodule with $\textrm{dim}V^{(0)}\geq 2$ and $\textrm{dim}V_{\al+\mu}^{(0)}\leq 1$ for all $\al\in\G$. Moreover, Claim 1 shows that $V^{(0)}$ is not a direct sum of some trivial $\mcr{S}_0$-submodules. So Proposition \ref{t:1.2} implies that, there exists $\eta\in\mbb{F}^{l_2+l_3}$ such that $V_{\ga+\mu}^{(0)}\not=\{0\}$ for all $\ga\in\G\backslash\{-\eta\}$, and that we can choose nonzero $v_{\ga,\textbf{0}}\in V_{\ga+\mu}^{(0)}$ with $\ga\in\G\backslash\{-\eta\}$
 such that
\begin{equation}\label{1.12}
D_{p,q}(x^{\al}).v_{\ga,\textbf{0}}=x^{\al}(\al_p\ptl_q-\al_q\ptl_p).v_{\ga,\textbf{0}}
=(\al_p(\ga_q+\eta_q)-\al_q(\ga_p+\eta_p)) v_{\al+\ga,\textbf{0}}
\end{equation}
for $p,q\in\overline{l_1+1,l}$, $\al\in\G\backslash\{0\}$ and $\ga\in\G\backslash\{-\eta, -\eta-\al\}$.

So next, we need to prove:\vspace{0.3cm}

{\it Claim 2.} $\eta=\mu$.\vspace{0.2cm}

It suffices to prove $\ptl(\eta-\mu)=0$ for all $\ptl\in\mcr{D}_2+\mcr{D}_3=\sum_{i\in\ove{l_1+1,l}}\mbb{F}\ptl_i$ (cf. (\ref{0.19})). Choose $v_{\ga,\textbf{0}}\in V_{\ga+\mu}^{(0)}$ for $\ga\in\G\backslash\{-\eta\}$ as in (\ref{1.12}). We give the proof in two cases.\vspace{0.3cm}

Case 1. $l\geq 4$.\vspace{0.3cm}

Pick $r\in\ove{1,l_1}$ if $l_1>0$ or $r\in\ove{1,l_2}$ if $l_1=0$. Take $s\in\ove{l_1+1,l}\backslash\{r\}$. Choose $\al\in\G\backslash\{0\}$ such that $\al_s\not=0$. Moreover, we choose $\ga\in\G$ such that $\ga+\eta\not=0\not=\ga+\al+\eta$. Then we must have
\begin{equation}\label{1.14}
\{0\} \not=\ker(\ga+\mu)\bigcap \ker\al\bigcap\big(\sum_{i\in\ove{l_1+1,l}\backslash\{r\}}\mbb{F}\ptl_i\big)\subseteq \ker(\eta-\mu).
\end{equation}
 Suppose there exists
\begin{equation}\label{1.13}
0\not=\ptl\in \ker(\ga+\mu)\bigcap \ker\al\bigcap\big(\sum_{i\in\ove{l_1+1,l}\backslash\{r\}}\mbb{F}\ptl_i\big)\backslash \ker(\eta-\mu).
 \end{equation}
We will see that this leads to a contradiction. By (\ref{2.2}) and (\ref{1.13}), we have
\begin{equation}
t^{1_{[r]}}\ptl.v_{\ga,\textbf{0}}\in V_{\ga+\mu}^{(0)} \textrm{ and } t^{1_{[r]}}\ptl.v_{\ga+\al,\textbf{0}}\in V_{\ga+\al+\mu}^{(0)}.
\end{equation}
 So $t^{1_{[r]}}\ptl.v_{\ga,\textbf{0}}=cv_{\ga,\textbf{0}}$ for some $c\in\mbb{F}$. Since (\ref{1.12}) implies
\begin{eqnarray}
0&=&[x^\al \ptl, t^{1_{[r]}}\ptl].v_{\ga,\textbf{0}}\nonumber\\
&=&x^\al \ptl. (t^{1_{[r]}}\ptl.v_{\ga,\textbf{0}})- t^{1_{[r]}}\ptl. x^\al \ptl.v_{\ga,\textbf{0}}\nonumber\\
&=&c \ptl(\ga+\eta)v_{\ga+\al,\textbf{0}}- \ptl(\ga+\eta)t^{1_{[r]}}\ptl.v_{\ga+\al,\textbf{0}},
\end{eqnarray}
and (\ref{1.13}) indicates $\ptl(\ga+\eta)=\ptl(\eta-\mu)\not=0$, we find $t^{1_{[r]}}\ptl.v_{\ga+\al,\textbf{0}}=cv_{\ga+\al,\textbf{0}}$. Thus (\ref{1.12}) and (\ref{1.13}) give rise to
\begin{eqnarray}
0&\not=&\ptl(\ga+\eta)v_{\ga+\al,\textbf{0}}\nonumber\\
&=& x^\al\ptl.v_{\ga,\textbf{0}}\nonumber\\
&=&-\frac{1}{\al_s}[D_{r,s}(x^\al), t^{1_{[r]}}\ptl].v_{\ga,\textbf{0}}\nonumber\\
&=&-\frac{1}{\al_s}(c D_{r,s}(x^\al).v_{\ga,\textbf{0}}- t^{1_{[r]}}\ptl. (D_{r,s}(x^\al).v_{\ga,\textbf{0}}))\nonumber\\
&=&-\frac{1}{\al_s}(c D_{r,s}(x^\al).v_{\ga,\textbf{0}}- c D_{r,s}(x^\al).v_{\ga,\textbf{0}})\nonumber\\
&=&0,
\end{eqnarray}
where in the fourth line, $D_{r,s}(x^\al).v_{\ga,\textbf{0}}=(\al_r(\ga_s+\eta_s)-\al_s(\ga_r+\eta_r)) v_{\ga+\al,\textbf{0}}$ by (\ref{1.12}). This is a contradiction. Thus (\ref{1.14}) holds. By changing the choice of $\ga$ in (\ref{1.14}), we can get
\begin{equation}
\ker\al\bigcap\big(\sum_{i\in\ove{l_1+1,l}\backslash\{r\}}\mbb{F}\ptl_i\big)\subseteq \ker(\eta-\mu).
\end{equation}
Moreover, by changing the choice of $\al$, we further find
\begin{equation}\label{1.15}
\sum_{i\in\ove{l_1+1,l}\backslash\{r\}}\mbb{F}\ptl_i\subseteq \ker(\eta-\mu).
\end{equation}
Recall that $r\in\ove{1,l_1}$ if $l_1>0$ and that $r\in\ove{1,l_2}$ if $l_1=0$. In the case that $r\in\ove{1,l_2}$, by picking another $r\in\ove{1,l_2}$, we can always get
\begin{equation}\label{1.16}
\sum_{i\in\ove{l_1+1,l}}\mbb{F}\ptl_i\subseteq \ker(\eta-\mu).
\end{equation}
So (\ref{1.15}) and (\ref{1.16}) show
\begin{equation}
\ptl(\eta-\mu)=0 \textrm{ for all }\ptl\in\sum_{i\in\ove{l_1+1,l}}\mbb{F}\ptl_i,
\end{equation}
which implies $\eta=\mu$.\vspace{0.3cm}

Case 2. $l=3$.  \vspace{0.3cm}

Since $l_1+l_2\geq 3$ and $l_2+l_3\geq 3$, we have $l_1=l_3=0$ and $l_2=l=3$.

Pick $r\in\ove{1,3}$. Take $s\in\ove{1,3}\backslash\{r\}$. Choose $\al\in\G\backslash\{0\}$ such that $\al_s\not=0$. Then 
\begin{equation}
\textrm{dim} \big(\ker\al\bigcap(\sum_{i\in\ove{1,3}\backslash\{r\}}\mbb{F}\ptl_i)\big)=1.
 \end{equation}
 Picking $0\not=\tilde{\ptl}_1\in\ker\al\bigcap\big(\sum_{i\in\ove{1,3}\backslash\{r\}}\mbb{F}\ptl_i\big)$, we shall first prove that
\begin{equation} \label{1.103}
\tilde{\ptl}_1(\eta-\mu)=0.
\end{equation}
Choose $\ga\in\G$ such that $\ga+\eta\not=0\not=\ga+\al+\eta$ and
\begin{equation}\label{1.108}
\tilde{\ptl}_1(2\ga+\eta+\mu)\not=0.
\end{equation}
If
\begin{equation}
\tilde{\ptl}_1(\ga+\mu)=0,
\end{equation}
then similar arguments as those in Case 1 give (\ref{1.103}).
If $\tilde{\ptl}_1(\ga+\mu)\not=0$, namely,
\begin{equation}\label{1.104}
\ker(\ga+\mu)\bigcap \ker\al\bigcap\big(\sum_{i\in\ove{1,3}\backslash\{r\}}\mbb{F}\ptl_i\big)=\{0\},
\end{equation}
we pick
\begin{equation} \label{1.110} 0\not=\tilde{\ptl}_2\in\ker(\ga-\al+\mu)\bigcap\big(\sum_{i\in\ove{1,3}\backslash\{r\}}\mbb{F}\ptl_i\big).
\end{equation}
Then $\{\tilde{\ptl}_1,\tilde{\ptl}_2\}$ is a basis of
$\sum_{i\in\ove{1,3}\backslash\{r\}}\mbb{F}\ptl_i$. Take
\begin{equation}\label{1.107}
0\not=\ptl\in\ker(\ga+\mu)\bigcap \big(\sum_{i\in\ove{1,3}\backslash\{r\}}\mbb{F}\ptl_i\big),\
0\not=\ptl'\in \ker\al \bigcap \ker(\ga+\eta), \
0\not=\ptl''\in \ker\al.
 \end{equation}
It follows that $\ptl=a_1\tilde{\ptl}_1+a_2\tilde{\ptl}_2$ for some $a_1,a_2\in\mbb{F}\backslash\{0\}$. So on one hand, by (\ref{2.2}) we have
\begin{eqnarray}
& &[x^{-\al}\ptl'',[x^\al \ptl', t^{1_{[r]}}\ptl]].v_{\ga,\textbf{0}}\nonumber\\
&=&\ptl(\al)\cdot(\ptl'(t^{1_{[r]}})\ptl''-\ptl''(t^{1_{[r]}})\ptl').v_{\ga,\textbf{0}}\nonumber\\
&=& \ptl(\al-\ga-\mu)\cdot(\ptl'(t^{1_{[r]}})\ptl''-\ptl''(t^{1_{[r]}})\ptl')(\ga+\mu)v_{\ga,\textbf{0}}\nonumber\\
&=& -a_1\tilde{\ptl}_1(\ga+\mu)\cdot(\ptl'(t^{1_{[r]}})\ptl''-\ptl''(t^{1_{[r]}})\ptl')(\ga+\mu)v_{\ga,\textbf{0}},\label{1.105}
\end{eqnarray}
where in the third line we use the fact that $\ptl(\ga+\mu)=0$, and in the fourth line we use $\ptl=a_1\tilde{\ptl}_1+a_2\tilde{\ptl}_2$ and (\ref{1.110}). On the other hand, by (\ref{1.12}) and (\ref{1.107}), we have
\begin{eqnarray}
& &[x^{-\al}\ptl'',[x^\al \ptl', t^{1_{[r]}}\ptl]].v_{\ga,\textbf{0}}\nonumber\\
&=&x^{-\al}\ptl''.x^\al \ptl'.(t^{1_{[r]}}\ptl.v_{\ga,\textbf{0}}) -x^{-\al}\ptl''. t^{1_{[r]}}\ptl.x^\al \ptl'.v_{\ga,\textbf{0}}\nonumber\\
& & - x^\al \ptl'.t^{1_{[r]}}\ptl.x^{-\al}\ptl''.v_{\ga,\textbf{0}} +t^{1_{[r]}}\ptl.x^\al \ptl'. (x^{-\al}\ptl''.v_{\ga,\textbf{0}})\nonumber\\
&=& - x^\al \ptl'.t^{1_{[r]}}\ptl.x^{-\al}\ptl''.v_{\ga,\textbf{0}} \nonumber\\
&=& - \ptl''(\ga+\eta)x^\al \ptl'.t^{1_{[r]}}(a_1\tilde{\ptl}_1+a_2\tilde{\ptl}_2).v_{\ga-\al,\textbf{0}} \nonumber\\
&=& - a_1\ptl''(\ga+\eta) x^\al \ptl'.t^{1_{[r]}}\tilde{\ptl}_1.v_{\ga-\al,\textbf{0}} \nonumber\\
&=& - a_1\ptl''(\ga+\eta) [x^\al \ptl',t^{1_{[r]}}\tilde{\ptl}_1].v_{\ga-\al,\textbf{0}} \nonumber\\
&=& - a_1\ptl''(\ga+\eta) \cdot \ptl'(t^{1_{[r]}}) x^\al\tilde{\ptl}_1.v_{\ga-\al,\textbf{0}} \nonumber\\
&=& - a_1\ptl'(t^{1_{[r]}})\cdot \tilde{\ptl}_1(\ga+\eta) \cdot \ptl''(\ga+\eta) v_{\ga,\textbf{0}}, \label{1.106}
\end{eqnarray}
where in the second equation $t^{1_{[r]}}\ptl.v_{\ga,\textbf{0}}\in V_{\ga+\mu}^{(0)}$ and $x^{-\al}\ptl''.v_{\ga,\textbf{0}}\in V_{\ga-\al+\mu}^{(0)}$, in the fourth equation $t^{1_{[r]}}\tilde{\ptl}_2.v_{\ga-\al,\textbf{0}}\in V_{\ga-\al+\mu}^{(0)}$ and $x^\al \ptl'. (t^{1_{[r]}}\tilde{\ptl}_2.v_{\ga-\al,\textbf{0}})=0$. If $\ptl'(t^{1_{[r]}})=0$, by taking $0\not=\ptl''\in  \ker\al\bigcap \ker(\ga+\mu)$ and by comparing (\ref{1.105}) with (\ref{1.106}), we get $\ptl''(t^{1_{[r]}})\not=0$ (cf. (\ref{1.104})) and $\ptl'(\ga+\mu)=0$, which together with (\ref{1.104}) and (\ref{1.107}) further implies
\begin{equation}
0\not=\ptl'\in \ker\al \bigcap \ker(\ga+\mu)\bigcap\big(\sum_{i\in\ove{1,3}\backslash\{r\}}\mbb{F}\ptl_i\big)=\{0\}.
\end{equation}
This is absurd. So we must have $\ptl'(t^{1_{[r]}})\not=0$. Therefore, taking $\ptl''=\tilde{\ptl}_1$ and comparing (\ref{1.105}) with (\ref{1.106}), we get
\begin{equation}\label{1.109}
(\tilde{\ptl}_1(\ga+\eta))^2=(\tilde{\ptl}_1(\ga+\mu))^2.
 \end{equation}
Thus by (\ref{1.108}) and (\ref{1.109}), we find
\begin{equation}
\tilde{\ptl}_1(\eta-\mu)=0.
\end{equation}
So (\ref{1.103}) holds. In other words, we have
\begin{equation}
\ker\al\bigcap(\sum_{i\in\ove{1,3}\backslash\{r\}}\mbb{F}\ptl_i)\subseteq \ker(\eta-\mu).
\end{equation}

So by changing the choice of $\al$, we can get
\begin{equation}
\sum_{i\in\ove{1,3}\backslash\{r\}}\mbb{F}\ptl_i\subseteq \ker(\eta-\mu).
\end{equation}
Moreover, by picking another $r\in\ove{1,3}$, we further obtain
\begin{equation}
\mcr{D}=\sum_{i\in\ove{1,3}}\mbb{F}\ptl_i\subseteq \ker(\eta-\mu),
\end{equation}
which indicates
$\eta=\mu$. So this claim holds.\vspace{0.4cm}

Claim 2 shows that
$V_{\ga+\mu}^{(0)}\not=\{0\}$ for all $\ga\in\G\backslash\{-\mu\}$, and that we can choose nonzero $v_{\ga,\textbf{0}}\in V_{\ga+\mu}^{(0)}$ with $\ga\in\G\backslash\{-\mu\}$ such that
\begin{equation}
D_{p,q}(x^{\al}).v_{\ga,\textbf{0}}=x^{\al}(\al_p\ptl_q-\al_q\ptl_p).v_{\ga,\textbf{0}}
=(\al_p(\ga_q+\mu_q)-\al_q(\ga_p+\mu_p)) v_{\al+\ga,\textbf{0}}
\end{equation}
for $p,q\in\overline{l_1+1,l}$, $\al\in\G\backslash\{0\}$ and $\ga\in\G\backslash\{-\mu,-\al-\mu\}$. Then we want to prove: \vspace{0.4cm}

{\it Claim 3.} If $\mu\in\G$, then
\begin{equation}
D_{p,q}(x^{\al}). v_{-\mu-\al,\textbf{0}}=0 \textrm{ for all } \al\in\G\backslash\{0\} \textrm{ and } p,q\in\overline{l_1+1,l}.
\end{equation}

Suppose
\begin{equation}\label{1.17}
D_{p,q}(x^{\al}). v_{-\mu-\al,\textbf{0}}\not=0
\end{equation}
for some $\al\in\G\backslash\{0\}$ and $p,q\in\overline{l_1+1,l}$, then this will lead to a contradiction.
Choose $r\in\ove{l_1+1,l}$ such that $\al_r\not=0$. Then (\ref{1.17}) implies that, there exists some $s\in\ove{l_1+1,l}\backslash\{r\}$ such that
\begin{equation}\label{1.18}
x^{\al}(\al_r\ptl_s-\al_s\ptl_r). v_{-\mu-\al,\textbf{0}}\not=0.
\end{equation}
 Fix such $s$. By (\ref{1.17}) we have $V_0^{(0)}\not=\{0\}$. Take nonzero $v_{-\mu,\textbf{0}}\in V_0^{(0)}$. Moreover, we pick $p\in\ove{1,l_1+l_2}\backslash\{r,s\}$. It follows from (\ref{2.2}) that
\begin{equation}
t^{1_{[p]}}(\al_r\ptl_s-\al_s\ptl_r).v_{-\mu-\al,\textbf{0}}\in V_{-\al}^{(0)} \textrm{ and } t^{1_{[p]}}(\al_r\ptl_s-\al_s\ptl_r).v_{-\mu,\textbf{0}}\in V_{0}^{(0)}.
\end{equation}
So $t^{1_{[p]}}(\al_r\ptl_s-\al_s\ptl_r).v_{-\mu-\al,\textbf{0}}=cv_{-\mu-\al,\textbf{0}}$ for some $c\in\mbb{F}$. Since
\begin{eqnarray}
&&t^{1_{[p]}}(\al_r\ptl_s-\al_s\ptl_r).(x^{\al}(\al_r\ptl_s-\al_s\ptl_r). v_{-\mu-\al,\textbf{0}})\nonumber\\
&=&x^{\al}(\al_r\ptl_s-\al_s\ptl_r).(t^{1_{[p]}}(\al_r\ptl_s-\al_s\ptl_r). v_{-\mu-\al,\textbf{0}})\nonumber\\
&=&c (x^{\al}(\al_r\ptl_s-\al_s\ptl_r).v_{-\mu-\al,\textbf{0}}),
\end{eqnarray}
we find $t^{1_{[p]}}(\al_r\ptl_s-\al_s\ptl_r).v_{-\mu,\textbf{0}}=cv_{-\mu,\textbf{0}}$ by (\ref{1.18}). Thus
\begin{eqnarray}
&&x^{\al}(\al_r\ptl_s-\al_s\ptl_r). v_{-\mu-\al,\textbf{0}}\nonumber\\
&=&\frac{1}{\al_r}[x^{\al}(\al_r\ptl_p-\al_p\ptl_r),t^{1_{[p]}}(\al_r\ptl_s-\al_s\ptl_r)]. v_{-\mu-\al,\textbf{0}}\nonumber\\
&=&\frac{1}{\al_r}(x^{\al}(\al_r\ptl_p-\al_p\ptl_r). (c v_{-\mu-\al,\textbf{0}})- c(x^{\al}(\al_r\ptl_p-\al_p\ptl_r). v_{-\mu-\al,\textbf{0}}))\nonumber\\
&=&0,
\end{eqnarray}
which contradicts (\ref{1.18}). So this claim holds.\vspace{0.4cm}

Moreover, we have: \vspace{0.4cm}

{\it Claim 4.} If $l_1>0$, then
\begin{equation}
x^{\al}\ptl_p . v_{\ga,\textbf{0}}=0 \textrm{ for all } p\in\ove{1,l_1},\ \al\in\G\backslash\{0\} \textrm{ and } \ga \in\G\backslash\{-\mu\}.
\end{equation}\vspace{0.3cm}

Suppose that there exist $p\in\ove{1,l_1}$, $\al\in\G\backslash\{0\}$ and $\ga \in\G\backslash\{-\mu\}$ such that
\begin{equation}\label{1.19}
x^{\al}\ptl_p . v_{\ga,\textbf{0}}\not=0,
\end{equation}
then we will get a contradiction. Choose $r\in\ove{l_1+1,l}$ such that $\al_r \not=0$.
 Pick $q\in\ove{1,l_1+l_2}\backslash\{p,r\}$. By (\ref{2.2}), we get
\begin{equation}
t^{1_{[q]}}\ptl_p.v_{\ga,\textbf{0}}\in V_{\ga+\mu}^{(0)} \textrm{ and } t^{1_{[q]}}\ptl_p.v_{\ga+\al,\textbf{0}}\in V_{\ga+\al+\mu}^{(0)}.
\end{equation}
So $t^{1_{[q]}}\ptl_p.v_{\ga,\textbf{0}}=cv_{\ga,\textbf{0}}$ for some $c\in\mbb{F}$. Moreover,
\begin{equation}
t^{1_{[q]}}\ptl_p. (x^{\al}\ptl_p. v_{\ga,\textbf{0}})=x^{\al}\ptl_p.(t^{1_{[q]}}\ptl_p. v_{\ga,\textbf{0}})=c (x^{\al}\ptl_p.v_{\ga,\textbf{0}}),
\end{equation}
which together with (\ref{1.19}) implies $t^{1_{[q]}}\ptl_p.v_{\ga+\al,\textbf{0}}=cv_{\ga+\al,\textbf{0}}$. Thus
\begin{eqnarray}
x^{\al}\ptl_p. v_{\ga,\textbf{0}}
&=&\frac{1}{\al_r}[x^{\al}(\al_r\ptl_q-\al_q\ptl_r),t^{1_{[q]}}\ptl_p]. v_{\ga,\textbf{0}}\nonumber\\
&=&\frac{1}{\al_r}(x^{\al}(\al_r\ptl_q-\al_q\ptl_r). (c v_{\ga,\textbf{0}})- c(x^{\al}(\al_r\ptl_q-\al_q\ptl_r). v_{\ga,\textbf{0}}))\nonumber\\
&=&0,
\end{eqnarray}
which contradicts (\ref{1.19}). So this claim holds.\vspace{0.4cm}

In summary, Claims 2--4 tell that $V_{\ga+\mu}^{(0)}\not=\{0\}$ for all $\ga\in\G\backslash\{-\mu\}$, and that there exist nonzero $v_{\ga,\textbf{0}}\in V_{\ga+\mu}^{(0)}$ with $\ga\in\G\backslash\{-\mu\}$ such that
\begin{equation}
D_{p,q}(x^{\al}).v_{\ga,\textbf{0}}=x^{\al}(\al_p\ptl_q-\al_q\ptl_p).v_{\ga,\textbf{0}}
=(\al_p(\ga_q+\mu_q)-\al_q(\ga_p+\mu_p))v_{\al+\ga,\textbf{0}}
\end{equation}
for $p,q\in\overline{1,l}$, $\al\in\G\backslash\{0\}$ and $\ga\in\G\backslash\{-\mu\}$.
This completes the proof of the lemma.  $\qquad\Box$ \vspace{0.4cm}

Lemma \ref{le:3.2} enables us to choose nonzero $v_{\ga,\textbf{0}}\in V_{\ga+\mu}^{(0)}$ with $\ga\in\G\backslash\{-\mu\}$ such that
\begin{equation}\label{0.18}
D_{p,q}(x^{\al}).v_{\ga,\textbf{0}}=x^{\al}(\al_p\ptl_q-\al_q\ptl_p).v_{\ga,\textbf{0}}
=(\al_p(\ga_q+\mu_q)-\al_q(\ga_p+\mu_p))v_{\al+\ga,\textbf{0}}
\end{equation}
for $p,q\in\overline{1,l}$, $\al\in\G\backslash\{0\}$ and $\ga\in\G\backslash\{-\mu\}$. Next, we define
$v_{\ga,\textbf{i}}$ for $\ga\in\G\backslash\{-\mu\}$ and $\textbf{0}\not=\textbf{i}\in \mbb{N}^{l_1+l_2}$ inductively as follows:

Let $\bar{\ga}=\ga+\mu$ for all $\ga\in\G$. For $\ga\in\G\backslash\{-\mu\}$ and $\textbf{0}\not=\textbf{i}\in \mbb{N}^{l_1+l_2}$, we define
\begin{equation}\label{1.20}
p_{\textbf{i}}=\min\{r\in\ove{1,l_1+l_2}\mid i_r\not=0\},
\end{equation}
\begin{equation}\label{1.21}
P(\textbf{i})=\{r\in\ove{1,l_1+l_2}\backslash\{p_{\textbf{i}}\}\mid i_r\not=0\},
\end{equation}
\begin{equation}\label{1.22}
Q(\ga,\textbf{i})=\{r\in\ove{l_1+1,l}\backslash\{p_{\textbf{i}}\}\mid \bar{\ga}_r\not=0\}.
\end{equation}

\begin{defi}\label{d:3.3}
We choose $v_{\ga,\textbf{0}}$ for $\ga\in\G\backslash\{-\mu\}$ as those in (\ref{0.18}).
Suppose that we have defined
$v_{\ga,\textbf{i}}$ for $\ga\in\G\backslash\{-\mu\}$ and $\textbf{i}\in \mbb{N}^{l_1+l_2}$ with $|\textbf{i}|\leq k$, where $k\geq 0$. We then define $v_{\ga,\textbf{i}}$ for $\ga\in\G\backslash\{-\mu\}$ and $\textbf{i}\in \mbb{N}^{l_1+l_2}$ with $|\textbf{i}|=k+1$ in three steps:\vspace{0.2cm}

\noindent{Step 1.} Defining $v_{\ga,\textbf{i}}$ for the case that $Q(\ga,\textbf{i})\not= \emptyset $.\vspace{0.2cm}

Set $q(\ga,\textbf{i})=\min Q(\ga,\textbf{i})$. We define
\begin{equation}\label{1.23}
v_{\ga,\textbf{i}}=\frac{1}{\bar{\ga}_{q(\ga,\textbf{i})}}\big(t^{1_{[p_{\textbf{i}}]}}\ptl_{q(\ga,\textbf{i})}.
v_{\ga,\textbf{i}-1_{[p_{\textbf{i}}]}}-i_{q(\ga,\textbf{i})}v_{\ga,\textbf{i}-1_{[q(\ga,\textbf{i})]}}\big).
\end{equation}

\noindent{Step 2.} Defining $v_{\ga,\textbf{i}}$ for the case that $Q(\ga,\textbf{i})= \emptyset $ and $P(\textbf{i})\not= \emptyset $.\vspace{0.2cm}

Set $p(\textbf{i})=\min P(\textbf{i})$. Since $\ga+\mu\not=0$ and $Q(\ga,\textbf{i})= \emptyset $, (\ref{1.22}) indicates $\bar{\ga}_{p_{\textbf{i}}}\not=0$. We define
\begin{equation}\label{1.24}
v_{\ga,\textbf{i}}=\frac{1}{\bar{\ga}_{p_{\textbf{i}}}} \big(t^{1_{[p(\textbf{i})]}}\ptl_{p_{\textbf{i}}}.
v_{\ga,\textbf{i}-1_{[p(\textbf{i})]}}-i_{p_{\textbf{i}}}v_{\ga,\textbf{i}-1_{[p_{\textbf{i}}]}}\big).
\end{equation}

\noindent{Step 3.} Defining $v_{\ga,\textbf{i}}$ for the case that $Q(\ga,\textbf{i})= \emptyset $ and $P(\textbf{i})= \emptyset $.\vspace{0.2cm}

Set $s(\textbf{i})=\min \{\ove{1,l_1+l_2}\backslash\{p_{\textbf{i}}\}\}$.
Since $(\ga,\textbf{i}-1_{[p_{\textbf{i}}]}+1_{[s(\textbf{i})]})$ belongs to the Case in Step 1 or Step 2, $v_{\ga,\textbf{i}-1_{[p_{\textbf{i}}]}+1_{[s(\textbf{i})]}}$ was defined in (\ref{1.23}) or (\ref{1.24}). Thus we can define
\begin{equation}\label{1.25}
v_{\ga,\textbf{i}}=t^{1_{[p_{\textbf{i}}]}}\ptl_{s(\textbf{i})}.v_{\ga,\textbf{i}-1_{[p_{\textbf{i}}]}+1_{[s(\textbf{i})]}}.
\end{equation}
\end{defi}\vspace{0.3cm}

Then we want to derive the action of the set (\ref{1.1}) on $\{v_{\ga,\textbf{i}}\mid \ga\in\G\backslash\{-\mu\}, \; \textbf{i}\in \mbb{N}^{l_1+l_2}\}$ defined in Definition \ref{d:3.3} in several lemmas.

\begin{lemma}\label{le:3.4}
For any $p \in\ove{1,l_1+l_2}$, $q \in\ove{1,l}\backslash\{p\}$ and $\be\in\G\backslash\{-\mu\}$, we have
\begin{equation}
t^{1_{[p]}}\ptl_q .v_{\be,\textbf{0}}=(\be_q+\mu_q) v_{\be,1_{[p]}}.
\end{equation}
\end{lemma}

\noindent{\bf Proof.} Let $\bar{\al}=\al+\mu$ for all $\al\in\G$. Define $Q(\be,1_{[p]})$ for $\be\in\G\backslash\{-\mu\}$ and $p \in\ove{1,l_1+l_2}$ as in (\ref{1.20})--(\ref{1.22}).
We divide $\be\in\G\backslash\{-\mu\}$ and $p\in \ove{1,l_1+l_2}$ into two cases.\vspace{0.2cm}

Case 1. $Q(\be,1_{[p]})\not=\emptyset$.\vspace{0.2cm}

Let $r=\min Q(\be,1_{[p]})$. Then (\ref{1.23}) shows
\begin{equation}\label{1.43}
t^{1_{[p]}}\ptl_r .v_{\be,\textbf{0}}=\bar{\be}_r v_{\be,1_{[p]}}.
\end{equation}
Fix some $s\in\ove{1,l_1+l_2}\backslash\{p,r\}$. Let
$\tilde{\ptl}_1=\bar{\be}_r\ptl_p - \bar{\be}_p\ptl_r$ and
$\tilde{\ptl}_2=\bar{\be}_r\ptl_s - \bar{\be}_s\ptl_r$. Then
\begin{equation}
\tilde{\ptl}_1(t^{1_{[p]}+1_{[s]}})\tilde{\ptl}_2-\tilde{\ptl}_2(t^{1_{[p]}+1_{[s]}})\tilde{\ptl}_1
=\bar{\be}_r (t^{1_{[s]}}\tilde{\ptl}_2 - t^{1_{[p]}}\tilde{\ptl}_1) \in\mcr{S}(l_{1},l_{2},l_{3},0;\G).
\end{equation}
Since (\ref{2.2}) gives
\begin{equation}
(t^{1_{[s]}}\tilde{\ptl}_2 - t^{1_{[p]}}\tilde{\ptl}_1) .v_{\be,\textbf{0}}=a_1 v_{\be,\textbf{0}},\ \  t^{1_{[p]}}(\bar{\be}_r\ptl_s - \bar{\be}_s\ptl_r)  .v_{\be,\textbf{0}}=a_2 v_{\be,\textbf{0}}
\end{equation}
for some $a_1,a_2\in\mbb{F}$, we have
\begin{equation}\label{1.46}
t^{1_{[p]}}(\bar{\be}_r\ptl_s - \bar{\be}_s\ptl_r)  .v_{\be,\textbf{0}}=\frac{1}{2\bar{\be}_r}[t^{1_{[p]}}\tilde{\ptl}_2, t^{1_{[s]}}\tilde{\ptl}_2 - t^{1_{[p]}}\tilde{\ptl}_1]  .v_{\be,\textbf{0}}=0.
\end{equation}
Moreover, from (\ref{2.2}) it follows that
\begin{equation}\label{1.45}
t^{1_{[s]}}(\bar{\be}_r\ptl_q - \bar{\be}_q\ptl_r)  .v_{\be,\textbf{0}}\in V_{\be+\mu}^{(0)} \quad \textrm{ for } q\in\ove{1,l}\backslash\{p,r,s\}.
\end{equation}
 So
\begin{equation}\label{1.44}
t^{1_{[p]}}(\bar{\be}_r\ptl_q - \bar{\be}_q\ptl_r)  .v_{\be,\textbf{0}}=\frac{1}{\bar{\be}_r}[t^{1_{[p]}}(\bar{\be}_r\ptl_s - \bar{\be}_s\ptl_r), t^{1_{[s]}}(\bar{\be}_r\ptl_q - \bar{\be}_q\ptl_r)]  .v_{\be,\textbf{0}}=0
\end{equation}
for $q\in\ove{1,l}\backslash\{p,r,s\}$.
Since $\bar{\be}_r\not=0$, (\ref{1.43}), (\ref{1.46}) and (\ref{1.44}) imply
\begin{equation}
t^{1_{[p]}}\ptl_q .v_{\be,\textbf{0}}=(\be_q+\mu_q) v_{\be,1_{[p]}} \ \ \textrm{ for } q\in\ove{1,l}\backslash\{p\}.
\end{equation}

Case 2. $Q(\be,1_{[p]})=\emptyset$.\vspace{0.2cm}

Pick $s \in\ove{1,l_1+l_2}\backslash\{p\}$. For any $q \in\ove{1,l}\backslash\{p,s\}$, (\ref{2.2}) and $Q(\be,1_{[p]})=\emptyset$ give rise to
\begin{equation}
t^{1_{[p]}}\ptl_s .v_{\be,\textbf{0}}=b_1 v_{\be,\textbf{0}} \textrm{ and }t^{1_{[s]}}\ptl_q  .v_{\be,\textbf{0}}=b_q v_{\be,\textbf{0}}
\end{equation}
with some $b_1,b_q\in\mbb{F}$. So
\begin{equation}
t^{1_{[p]}}\ptl_q. v_{\be,\textbf{0}}=[t^{1_{[p]}}\ptl_s, t^{1_{[s]}}\ptl_q ]  .v_{\be,\textbf{0}}=0
\end{equation}
for $q\in\ove{1,l}\backslash\{p,s\}$. By picking another $s\in\ove{1,l_1+l_2}\backslash\{p\}$, we can get
\begin{equation}
t^{1_{[p]}}\ptl_q. v_{\be,\textbf{0}}=0=(\be_q+\mu_q) v_{\be,1_{[p]}} \textrm{ for } q\in\ove{1,l}\backslash\{p\}.
\end{equation}
This completes the proof of the lemma. $\qquad\Box$\vspace{0.3cm}

\begin{lemma}\label{le:3.5}
For any $\be\in\G\backslash\{-\mu\}$ and $p,q \in\ove{1,l_1+l_2}$ with $p\not=q$, we have
\begin{eqnarray}
D_{q,p}(t^{1_{[p]}+1_{[q]}}) .v_{\be,\textbf{0}}&=&(t^{1_{[p]}}\ptl_p - t^{1_{[q]}}\ptl_q) .v_{\be,\textbf{0}}\nonumber\\
&=&(\be_p+\mu_p) v_{\be,1_{[p]}}-(\be_q+\mu_q) v_{\be,1_{[q]}}.\label{1.47}
\end{eqnarray}
\end{lemma}

\noindent{\bf Proof.} Let $\bar{\al}=\al+\mu$ for all $\al\in\G$. Pick any $\be\in\G\backslash\{-\mu\}$ and $p,q \in\ove{1,l_1+l_2}$ with $p\not=q$. Then we shall give the proof in two cases.\vspace{0.2cm}

Case 1. There exists $s\in\ove{1,l}\backslash\{p,q\}$ such that $\bar{\be}_s\not=0$.\vspace{0.2cm}

Let
$\tilde{\ptl}_1=\bar{\be}_s\ptl_p - \bar{\be}_p\ptl_s$ and
$\tilde{\ptl}_2=\bar{\be}_s\ptl_q - \bar{\be}_q\ptl_s$. Then
\begin{eqnarray}
[t^{1_{[p]}}\tilde{\ptl}_2 , t^{1_{[q]}}\tilde{\ptl}_1]&=&\bar{\be}_s (t^{1_{[p]}}\tilde{\ptl}_1 - t^{1_{[q]}}\tilde{\ptl}_2)\nonumber\\
&=&\bar{\be}_s^2 (t^{1_{[p]}}\ptl_p - t^{1_{[q]}}\ptl_q)+ \bar{\be}_s(\bar{\be}_q t^{1_{[q]}}\ptl_s-\bar{\be}_p t^{1_{[p]}}\ptl_s) \nonumber\\
&=&\bar{\be}_s^2 D_{q,p}(t^{1_{[p]}+1_{[q]}})+ \bar{\be}_s(\bar{\be}_q t^{1_{[q]}}\ptl_s-\bar{\be}_p t^{1_{[p]}}\ptl_s)\label{4.13}
\end{eqnarray}
So by (\ref{4.13}) and Lemma \ref{le:3.4}, we have
\begin{eqnarray}
&& D_{q,p}(t^{1_{[p]}+1_{[q]}}) .v_{\be,\textbf{0}}\nonumber\\
&=&\frac{1}{\bar{\be}_s^2}([t^{1_{[p]}}\tilde{\ptl}_2 , t^{1_{[q]}}\tilde{\ptl}_1] - \bar{\be}_s(\bar{\be}_q t^{1_{[q]}}\ptl_s-\bar{\be}_p t^{1_{[p]}}\ptl_s)).v_{\be,\textbf{0}}\nonumber\\
&=& \bar{\be}_p v_{\be,1_{[p]}}-\bar{\be}_q v_{\be,1_{[q]}},
\end{eqnarray}
which coincides with (\ref{1.47}).\vspace{0.2cm}

Case 2. $\bar{\be}_s=0$ for all $s\in\ove{1,l}\backslash\{p,q\}$.\vspace{0.2cm}

Subcase 2.1. $\bar{\be}_p\not=0$ and $\bar{\be}_q=0$.\vspace{0.1cm}

Take $r=\min\{\ove{1,l_1+l_2}\backslash\{p\}\}$. By (\ref{1.25}) we know that
\begin{equation}\label{1.48}
v_{\be,1_{[p]}}=t^{1_{[p]}}\ptl_r. v_{\be,1_{[r]}}.
\end{equation}
Since $\bar{\be}_q=0$ and $\bar{\be}_r=0$, Lemma \ref{le:3.4} tells that
\begin{equation}\label{1.49}
D_{r,q}(t^{1_{[q]}+1_{[r]}}) .v_{\be,\textbf{0}}=(t^{1_{[q]}}\ptl_q-t^{1_{[r]}}\ptl_r) .v_{\be,\textbf{0}}=[t^{1_{[q]}}\ptl_r, t^{1_{[r]}}\ptl_q] .v_{\be,\textbf{0}}=0.
\end{equation}
So (\ref{1.48}), (\ref{1.49}) and Lemma \ref{le:3.4} imply
\begin{eqnarray}
D_{q,p}(t^{1_{[p]}+1_{[q]}}) .v_{\be,\textbf{0}}&=&(D_{q,p}(t^{1_{[p]}+1_{[q]}})+D_{r,q}(t^{1_{[q]}+1_{[r]}})).v_{\be,\textbf{0}}\nonumber\\
&=&(t^{1_{[p]}}\ptl_p-t^{1_{[r]}}\ptl_r) .v_{\be,\textbf{0}}\nonumber\\
&=& [t^{1_{[p]}}\ptl_r,t^{1_{[r]}}\ptl_p] .v_{\be,\textbf{0}}\nonumber\\
&=& t^{1_{[p]}}\ptl_r. (\bar{\be}_p v_{\be,1_{[r]}})\nonumber\\
&=& \bar{\be}_p v_{\be,1_{[p]}},
\end{eqnarray}
which coincides with (\ref{1.47}).\vspace{0.1cm}

Subcase 2.2. $\bar{\be}_p=0$ and $\bar{\be}_q\not=0$.\vspace{0.1cm}

We go back to case 2.1 by interchanging $p$ and $q$, namely
\begin{equation}
D_{q,p}(t^{1_{[p]}+1_{[q]}}) .v_{\be,\textbf{0}}=-D_{p,q}(t^{1_{[p]}+1_{[q]}}) .v_{\be,\textbf{0}}=-\bar{\be}_q v_{\be,1_{[q]}},
\end{equation}
which coincides with (\ref{1.47}).\vspace{0.1cm}

Subcase 2.3. $\bar{\be}_p\not=0$ and $\bar{\be}_q\not=0$.\vspace{0.1cm}

Pick $r\in\ove{1,l_1+l_2}\backslash\{p,q\}$, then $\bar{\be}_r=0$. By Subcase 2.1, we know that
\begin{equation}
D_{r,p}(t^{1_{[p]}+1_{[r]}}) .v_{\be,\textbf{0}}=\bar{\be}_p v_{\be,1_{[p]}},
\end{equation}
\begin{equation}
D_{r,q}(t^{1_{[q]}+1_{[r]}}) .v_{\be,\textbf{0}}=\bar{\be}_q v_{\be,1_{[p]}}.
\end{equation}
So the difference of the above two equations is
\begin{eqnarray}
D_{q,p}(t^{1_{[p]}+1_{[q]}}) .v_{\be,\textbf{0}} &=&(D_{r,p}(t^{1_{[p]}+1_{[r]}})-D_{r,q}(t^{1_{[r]}+1_{[q]}})) .v_{\be,\textbf{0}}\nonumber\\
&=& \bar{\be}_p v_{\be,1_{[p]}}-\bar{\be}_q v_{\be,1_{[q]}},
\end{eqnarray}
which coincides with (\ref{1.47}).

This completes the proof of the lemma. $\qquad\Box$\vspace{0.3cm}

\begin{lemma}\label{le:3.6}
For any $p \in\ove{1,l_1+l_2}$, $q \in\ove{1,l}\backslash\{p\}$, $\be\in\G\backslash\{-\mu\}$ and $\textbf{i}\in \mbb{N}^{l_1+l_2}$, we have
\begin{equation}\label{1.50}
t^{1_{[p]}}\ptl_q .v_{\be,\textbf{i}}=(\be_q+\mu_q) v_{\be,\textbf{i}+1_{[p]}} + i_q v_{\be,\textbf{i}+1_{[p]}-1_{[q]}}.
\end{equation}
\end{lemma}

\noindent{\bf Proof.} Let $\bar{\al}=\al+\mu$ for all $\al\in\G$. We give the proof by induction on $|\textbf{i}|$. Recall that Lemma \ref{le:3.4} has given the proof for the case $|\textbf{i}|=0$. Suppose this lemma holds for $\textbf{i}\in \mbb{N}^{l_1+l_2}$ with $|\textbf{i}|\leq k$, where $k \geq 0$. Then it suffices to prove that
(\ref{1.50}) holds for $\textbf{i}\in \mbb{N}^{l_1+l_2}$ with $|\textbf{i}|= k+1$.

Fix any $\be\in\G\backslash\{-\mu\}$ and $\textbf{i}\in \mbb{N}^{l_1+l_2}$ with $|\textbf{i}|= k+1$. It suffices to prove that (\ref{1.50}) holds for all $p \in\ove{1,l_1+l_2}$ and $q \in\ove{1,l}\backslash\{p\}$. Define $p_{\textbf{i}}$, $P(\textbf{i})$ and $Q(\be,\textbf{i})$ as in (\ref{1.20})--(\ref{1.22}). Then we give the proof in three cases.\vspace{0.2cm}

Case 1. $Q(\be,\textbf{i})\not= \emptyset $.\vspace{0.2cm}

Let $q(\be,\textbf{i})=\min Q(\be,\textbf{i})$. Recall that (cf. (\ref{1.23}))
\begin{equation}\label{1.54}
v_{\be,\textbf{i}} = \frac{1}{\bar{\be}_{q(\be,\textbf{i})}} ( t^{1_{[p_{\textbf{i}}]}}\ptl_{q(\be,\textbf{i})}.
v_{\be,\textbf{i}-1_{[p_{\textbf{i}}]}}-i_{q(\be,\textbf{i})} v_{\be,\textbf{i}-1_{[q(\be,\textbf{i})]}}).
\end{equation}
Fixing any $p \in\ove{1,l_1+l_2}$, we derive the action of $t^{1_{[p]}}\ptl_q$ for $q \in\ove{1,l}\backslash\{p\}$ in two subcases.\vspace{0.2cm}

Subcase 1.1. $\bar{\be}_q=0$ for all $q \in\ove{1,l}\backslash\{p\}$.\vspace{0.2cm}

Expression (\ref{1.54}) and the induction hypothesis give rise to
\begin{eqnarray}
t^{1_{[p]}}\ptl_q .v_{\be,\textbf{i}} &= &\frac{1}{\bar{\be}_{q(\be,\textbf{i})}}\big(t^{1_{[p]}}\ptl_q .t^{1_{[p_{\textbf{i}}]}}\ptl_{q(\be,\textbf{i})}.
v_{\be,\textbf{i}-1_{[p_{\textbf{i}}]}}-i_{q(\be,\textbf{i})}t^{1_{[p]}}\ptl_q .v_{\be,\textbf{i}-1_{[q(\be,\textbf{i})]}}\big)\nonumber\\
&= &\frac{1}{\bar{\be}_{q(\be,\textbf{i})}}\big(
t^{1_{[p_{\textbf{i}}]}}\ptl_{q(\be,\textbf{i})}.t^{1_{[p]}}\ptl_q.
v_{\be,\textbf{i}-1_{[p_{\textbf{i}}]}}
-i_{q(\be,\textbf{i})}t^{1_{[p]}}\ptl_q .v_{\be,\textbf{i}-1_{[q(\be,\textbf{i})]}}\nonumber\\
& & + (\delta_{q ,p_{\textbf{i}}}t^{1_{[p]}}\ptl_{q(\be,\textbf{i})}
- \delta_{q(\be,\textbf{i}),p}t^{1_{[p_{\textbf{i}}]}}\ptl_q ).
v_{\be,\textbf{i}-1_{[p_{\textbf{i}}]}}\big) \nonumber\\
&= &\frac{1}{\bar{\be}_{q(\be,\textbf{i})}}\big( (i_q-\delta_{q,p_{\textbf{i}}})t^{1_{[p_{\textbf{i}}]}}\ptl_{q(\be,\textbf{i})}.
v_{\be,\textbf{i}+1_{[p]}-1_{[p_{\textbf{i}}]}-1_{[q]}}
-i_{q(\be,\textbf{i})}t^{1_{[p]}}\ptl_q .v_{\be,\textbf{i}-1_{[q(\be,\textbf{i})]}}\nonumber\\
& & + (\delta_{q ,p_{\textbf{i}}}t^{1_{[p]}}\ptl_{q(\be,\textbf{i})}
- \delta_{q(\be,\textbf{i}),p}t^{1_{[p_{\textbf{i}}]}}\ptl_q ).
v_{\be,\textbf{i}-1_{[p_{\textbf{i}}]}} \big)\nonumber\\
&= & i_q v_{\be,\textbf{i}+1_{[p]}-1_{[q]}}
\end{eqnarray}
for $q \in\ove{1,l}\backslash\{p\}$, which coincides with (\ref{1.50}).\vspace{0.2cm}

Subcase 1.2. There exists $q \in\ove{1,l}\backslash\{p\}$ such that $\bar{\be}_q\not=0$.\vspace{0.2cm}

Take $r=\min \{q \in\ove{1,l}\backslash\{p\} \mid \bar{\be}_q\not=0\}$. First we want to derive the action of $t^{1_{[p]}}\ptl_r$.

If $p<p_{\textbf{i}}$, then the definition of (\ref{1.23}) gives
\begin{equation}\label{1.111}
t^{1_{[p]}}\ptl_r .v_{\be,\textbf{i}}=\bar{\be}_r v_{\be,\textbf{i}+1_{[p]}} + i_r v_{\be,\textbf{i}+1_{[p]}-1_{[r]}}.
\end{equation}

If $p\geq p_{\textbf{i}}$, then by (\ref{1.54}) and the induction hypothesis, we get
\begin{eqnarray}
t^{1_{[p]}}\ptl_r .v_{\be,\textbf{i}} &= &\frac{1}{\bar{\be}_{q(\be,\textbf{i})}}\big(t^{1_{[p]}}\ptl_r .t^{1_{[p_{\textbf{i}}]}}\ptl_{q(\be,\textbf{i})}.
v_{\be,\textbf{i}-1_{[p_{\textbf{i}}]}}-i_{q(\be,\textbf{i})}t^{1_{[p]}}\ptl_r .v_{\be,\textbf{i}-1_{[q(\be,\textbf{i})]}}\big)\nonumber\\
&= &\frac{1}{\bar{\be}_{q(\be,\textbf{i})}}\big(
t^{1_{[p_{\textbf{i}}]}}\ptl_{q(\be,\textbf{i})}.t^{1_{[p]}}\ptl_r.
v_{\be,\textbf{i}-1_{[p_{\textbf{i}}]}}
-i_{q(\be,\textbf{i})}t^{1_{[p]}}\ptl_r .v_{\be,\textbf{i}-1_{[q(\be,\textbf{i})]}}\nonumber\\
& & + (\delta_{r ,p_{\textbf{i}}}t^{1_{[p]}}\ptl_{q(\be,\textbf{i})}
- \delta_{q(\be,\textbf{i}),p}t^{1_{[p_{\textbf{i}}]}}\ptl_r ).
v_{\be,\textbf{i}-1_{[p_{\textbf{i}}]}}\big) \nonumber\\
&= &\frac{1}{\bar{\be}_{q(\be,\textbf{i})}}\big(\bar{\be}_r t^{1_{[p_{\textbf{i}}]}}\ptl_{q(\be,\textbf{i})}.v_{\be,\textbf{i}+1_{[p]}-1_{[p_{\textbf{i}}]}} +(i_r-\delta_{r,p_{\textbf{i}}})t^{1_{[p_{\textbf{i}}]}}\ptl_{q(\be,\textbf{i})}.
v_{\be,\textbf{i}+1_{[p]}-1_{[p_{\textbf{i}}]}-1_{[r]}}\nonumber\\
& &
-i_{q(\be,\textbf{i})}t^{1_{[p]}}\ptl_r .v_{\be,\textbf{i}-1_{[q(\be,\textbf{i})]}} + (\delta_{r ,p_{\textbf{i}}}t^{1_{[p]}}\ptl_{q(\be,\textbf{i})}
- \delta_{q(\be,\textbf{i}),p}t^{1_{[p_{\textbf{i}}]}}\ptl_r ).
v_{\be,\textbf{i}-1_{[p_{\textbf{i}}]}} \big)\nonumber\\
& =& \bar{\be}_r v_{\be,\textbf{i}+1_{[p]}} + i_r v_{\be,\textbf{i}+1_{[p]}-1_{[r]}},\label{1.51}
\end{eqnarray}
where
\begin{eqnarray}
& & t^{1_{[p_{\textbf{i}}]}}\ptl_{q(\be,\textbf{i})}.v_{\be,\textbf{i}+1_{[p]}-1_{[p_{\textbf{i}}]}}\nonumber\\
&=& \bar{\be}_{q(\be,\textbf{i})}v_{\be,\textbf{i}+1_{[p]}}+ (i_{q(\be,\textbf{i})}+\delta_{q(\be,\textbf{i}) ,p})v_{\be,\textbf{i}+1_{[p]}-1_{[q(\be,\textbf{i})]}}
\end{eqnarray}
because of $p\geq p_{\textbf{i}}$ and (\ref{1.23}).

Next we want to derive the action of $t^{1_{[p]}}\ptl_q$ for $q \in\ove{1,l}\backslash\{p,r\}$. By (\ref{1.54}) and the induction hypothesis, we have
\begin{eqnarray}
& & t^{1_{[p]}}(\bar{\be}_r \ptl_q-\bar{\be}_q \ptl_r) .v_{\be,\textbf{i}} \nonumber\\
&= &\frac{1}{\bar{\be}_{q(\be,\textbf{i})}}\big(t^{1_{[p]}}(\bar{\be}_r \ptl_q-\bar{\be}_q \ptl_r) .t^{1_{[p_{\textbf{i}}]}}\ptl_{q(\be,\textbf{i})}.
v_{\be,\textbf{i}-1_{[p_{\textbf{i}}]}}
-i_{q(\be,\textbf{i})}t^{1_{[p]}}(\bar{\be}_r \ptl_q-\bar{\be}_q \ptl_r) .v_{\be,\textbf{i}-1_{[q(\be,\textbf{i})]}}\big)\nonumber\\
&= &\frac{1}{\bar{\be}_{q(\be,\textbf{i})}}\big(
t^{1_{[p_{\textbf{i}}]}}\ptl_{q(\be,\textbf{i})}.t^{1_{[p]}}(\bar{\be}_r \ptl_q-\bar{\be}_q \ptl_r).
v_{\be,\textbf{i}-1_{[p_{\textbf{i}}]}}
-i_{q(\be,\textbf{i})}t^{1_{[p]}}(\bar{\be}_r \ptl_q-\bar{\be}_q \ptl_r) .v_{\be,\textbf{i}-1_{[q(\be,\textbf{i})]}}\nonumber\\
& & + ((\bar{\be}_r \delta_{q ,p_{\textbf{i}}}-\bar{\be}_q \delta_{r ,p_{\textbf{i}}})t^{1_{[p]}}\ptl_{q(\be,\textbf{i})}
- \delta_{q(\be,\textbf{i}),p}t^{1_{[p_{\textbf{i}}]}}(\bar{\be}_r \ptl_q-\bar{\be}_q \ptl_r) ).
v_{\be,\textbf{i}-1_{[p_{\textbf{i}}]}} \big)\nonumber\\
&= &\frac{1}{\bar{\be}_{q(\be,\textbf{i})}}\big(
t^{1_{[p_{\textbf{i}}]}}\ptl_{q(\be,\textbf{i})}.(\bar{\be}_r (i_q-\delta_{q,p_{\textbf{i}}})v_{\be,\textbf{i}+1_{[p]}-1_{[p_{\textbf{i}}]}-1_{[q]}}-\bar{\be}_q (i_r-\delta_{r,p_{\textbf{i}}})v_{\be,\textbf{i}+1_{[p]}-1_{[p_{\textbf{i}}]}-1_{[r]}})
\nonumber\\
& & + ((\bar{\be}_r \delta_{q ,p_{\textbf{i}}}-\bar{\be}_q \delta_{r ,p_{\textbf{i}}})t^{1_{[p]}}\ptl_{q(\be,\textbf{i})}
- \delta_{q(\be,\textbf{i}),p}t^{1_{[p_{\textbf{i}}]}}(\bar{\be}_r \ptl_q-\bar{\be}_q \ptl_r) ).
v_{\be,\textbf{i}-1_{[p_{\textbf{i}}]}} \nonumber\\
& &
-i_{q(\be,\textbf{i})}t^{1_{[p]}}(\bar{\be}_r \ptl_q-\bar{\be}_q \ptl_r) .v_{\be,\textbf{i}-1_{[q(\be,\textbf{i})]}}\big)\nonumber\\
&= & \bar{\be}_r i_q v_{\be,\textbf{i}+1_{[p]}-1_{[q]}}-\bar{\be}_q i_r v_{\be,\textbf{i}+1_{[p]}-1_{[r]}}\label{1.52}
\end{eqnarray}
for any $q \in\ove{1,l}\backslash\{p,r\}$. Since $\bar{\be}_r\not=0$, (\ref{1.111}), (\ref{1.51}) and (\ref{1.52}) indicate
\begin{equation}\label{1.53}
t^{1_{[p]}}\ptl_q .v_{\be,\textbf{i}}=\bar{\be}_q v_{\be,\textbf{i}+1_{[p]}} + i_q v_{\be,\textbf{i}+1_{[p]}-1_{[q]}}
\end{equation}
for any $q \in\ove{1,l}\backslash\{p\}$, which coincides with (\ref{1.50}).\vspace{0.2cm}

Case 2. $Q(\be,\textbf{i})= \emptyset $ and $P(\textbf{i})\not= \emptyset $.\vspace{0.2cm}

Let $p(\textbf{i})=\min P(\textbf{i})$. Recall that (cf. (\ref{1.24}))
\begin{equation}
v_{\be,\textbf{i}}=\frac{1}{\bar{\be}_{p_{\textbf{i}}}} \big(t^{1_{[p(\textbf{i})]}}\ptl_{p_{\textbf{i}}}.
v_{\be,\textbf{i}-1_{[p(\textbf{i})]}}-i_{p_{\textbf{i}}}v_{\be,\textbf{i}-1_{[p_{\textbf{i}}]}}\big).
\end{equation}
The rest proof of this case is analogous to Case 1. We omit the details.\vspace{0.2cm}

Case 3. $Q(\be,\textbf{i})= \emptyset $ and $P(\textbf{i})= \emptyset $.\vspace{0.2cm}

Observe that $Q(\be,\textbf{i})= \emptyset $ and $P(\textbf{i})= \emptyset $ mean
\begin{equation}\label{1.55}
i_r=0 \textrm{ and } \bar{\be}_r=0 \textrm{ for all } r\in\ove{1,l}\backslash\{p_{\textbf{i}}\}.
\end{equation}
Pick any $r\in \ove{1,l_1+l_2}\backslash\{p_{\textbf{i}}\}$.
Then $(\be,\textbf{i}-1_{[p_{\textbf{i}}]}+1_{[r]})$ belongs to Case 1 or Case 2, which implies
\begin{equation}\label{1.56}
t^{1_{[p]}}\ptl_{q}.v_{\be,\textbf{i}-1_{[p_{\textbf{i}}]}+1_{[r]}}=\bar{\be}_q v_{\be,\textbf{i}-1_{[p_{\textbf{i}}]}+1_{[r]}+1_{[p]}} + (i_q -\delta_{p_{\textbf{i}},q}+\delta_{r,q}) v_{\be,\textbf{i}-1_{[p_{\textbf{i}}]}+1_{[r]}+1_{[p]}-1_{[q]}}
\end{equation}
for $p \in\ove{1,l_1+l_2}$ and $q \in\ove{1,l}\backslash\{p\}$.
In particular, by (\ref{1.55}) we have
\begin{equation}\label{1.112}
t^{1_{[p_{\textbf{i}}]}}\ptl_{r}.v_{\be,\textbf{i}-1_{[p_{\textbf{i}}]}+1_{[r]}}=v_{\be,\textbf{i}}.
\end{equation}
So for any $q\in\ove{1,l}\backslash\{p_{\textbf{i}},r\}$, we obtain
\begin{eqnarray}
t^{1_{[p_{\textbf{i}}]}}\ptl_q .v_{\be,\textbf{i}}& = &t^{1_{[p_{\textbf{i}}]}}\ptl_q .
t^{1_{[p_{\textbf{i}}]}}\ptl_r.
v_{\be,\textbf{i}-1_{[p_{\textbf{i}}]}+1_{[r]}}\nonumber\\
&=& t^{1_{[p_{\textbf{i}}]}}\ptl_r.t^{1_{[p_{\textbf{i}}]}}\ptl_q .v_{\be,\textbf{i}-1_{[p_{\textbf{i}}]}+1_{[r]}}
\nonumber\\
&=& 0\label{1.113}
\end{eqnarray}
by (\ref{1.55}), (\ref{1.56}) and (\ref{1.112}). By changing the choice of $r\in \ove{1,l_1+l_2}\backslash\{p_{\textbf{i}}\}$, we get
\begin{equation}\label{1.114}
t^{1_{[p_{\textbf{i}}]}}\ptl_q .v_{\be,\textbf{i}}=0 \textrm{ for }q\in\ove{1,l}\backslash\{p_{\textbf{i}}\},
\end{equation}
which coincides with (\ref{1.50}).

Fix any $p\in \ove{1,l_1+l_2}\backslash\{p_{\textbf{i}}\}$. Then (\ref{1.23}) or (\ref{1.24}) gives
\begin{equation}\label{1.115}
t^{1_{[p]}}\ptl_{p_{\textbf{i}}} .v_{\be,\textbf{i}}=\bar{\be}_{p_{\textbf{i}}} v_{\be,\textbf{i}+1_{[p]}} + i_{p_{\textbf{i}}} v_{\be,\textbf{i}+1_{[p]}-1_{[p_{\textbf{i}}]}}.
\end{equation}
Pick $s\in \ove{1,l_1+l_2}\backslash\{p,p_{\textbf{i}}\}$. By (\ref{1.55}), (\ref{1.56}) and (\ref{1.112}), we have
\begin{equation}\label{1.116}
t^{1_{[p_{\textbf{i}}]}}\ptl_{s}.v_{\be,\textbf{i}-1_{[p_{\textbf{i}}]}+1_{[s]}}=v_{\be,\textbf{i}}
\end{equation}
and
\begin{equation}\label{1.117}
t^{1_{[p_{\textbf{i}}]}}\ptl_s .v_{\be,\textbf{i}-1_{[p_{\textbf{i}}]}+1_{[p]}}=0.
\end{equation}
So for any $q\in\ove{1,l}\backslash\{p,p_{\textbf{i}}\}$, we get
\begin{eqnarray}
t^{1_{[p]}}\ptl_q .v_{\be,\textbf{i}}
& = &t^{1_{[p]}}\ptl_q .
t^{1_{[p_{\textbf{i}}]}}\ptl_s.
v_{\be,\textbf{i}-1_{[p_{\textbf{i}}]}+1_{[s]}}\nonumber\\
&=& t^{1_{[p_{\textbf{i}}]}}\ptl_s.t^{1_{[p]}}\ptl_q .v_{\be,\textbf{i}-1_{[p_{\textbf{i}}]}+1_{[s]}}
\nonumber\\
&=& \delta_{q,s} t^{1_{[p_{\textbf{i}}]}}\ptl_s .v_{\be,\textbf{i}-1_{[p_{\textbf{i}}]}+1_{[p]}}
\nonumber\\
&=& 0\label{1.118}
\end{eqnarray}
by (\ref{1.116}) and (\ref{1.117}). So in this case, (\ref{1.50}) follows from
(\ref{1.114}), (\ref{1.115}) and (\ref{1.118}).

This completes the proof of the lemma. $\qquad\Box$

\begin{coral}\label{c:3.7}
For any $\be\in\G\backslash\{-\mu\}$, $\textbf{i}\in \mbb{N}^{l_1+l_2}$ and $p,q \in\ove{1,l_1+l_2}$ with $p\not=q$, we have
\begin{eqnarray}
& & D_{q,p}(t^{1_{[p]}+1_{[q]}}) .v_{\be,\textbf{i}}\nonumber\\
&=&(t^{1_{[p]}}\ptl_p - t^{1_{[q]}}\ptl_q) .v_{\be,\textbf{i}}\nonumber\\
&=&(\be_p+\mu_p) v_{\be,\textbf{i}+1_{[p]}}-(\be_q+\mu_q) v_{\be,\textbf{i}+1_{[q]}}+(i_p-i_q)v_{\be,\textbf{i}}.\label{1.57}
\end{eqnarray}
\end{coral}

\noindent{\bf Proof.} For any $\be\in\G\backslash\{-\mu\}$, $\textbf{i}\in \mbb{N}^{l_1+l_2}$ and $p,q \in\ove{1,l_1+l_2}$ with $p\not=q$, we have
\begin{eqnarray}
& & D_{q,p}(t^{1_{[p]}+1_{[q]}}) .v_{\be,\textbf{i}}\nonumber\\
&=&[t^{1_{[p]}}\ptl_q, t^{1_{[q]}}\ptl_p] .v_{\be,\textbf{i}}\nonumber\\
&=&t^{1_{[p]}}\ptl_q.t^{1_{[q]}}\ptl_p.v_{\be,\textbf{i}}-  t^{1_{[q]}}\ptl_p.t^{1_{[p]}}\ptl_q .v_{\be,\textbf{i}}\nonumber\\
&=&(\be_p+\mu_p) v_{\be,\textbf{i}+1_{[p]}}-(\be_q+\mu_q) v_{\be,\textbf{i}+1_{[q]}}+(i_p-i_q)v_{\be,\textbf{i}}
\end{eqnarray}
by Lemma \ref{le:3.6}.    $\qquad\Box$

\begin{lemma}\label{le:3.8}
For any $r \in\ove{1,l}$, $\be\in\G\backslash\{-\mu\}$ and $\textbf{i}\in \mbb{N}^{l_1+l_2}$, we have
\begin{equation}\label{1.26}
\ptl_r.v_{\be,\textbf{i}}=(\be_r+\mu_r) v_{\be,\textbf{i}}+ i_r v_{\be,\textbf{i}-1_{[r]}}.
\end{equation}
\end{lemma}

\noindent{\bf Proof.} Let $\bar{\al}=\al+\mu$ for all $\al\in\G$. We shall prove this lemma by induction on $|\textbf{i}|$. We have already known that (\ref{1.26}) holds when $|\textbf{i}|=0$ (cf. (\ref{2.2})). Suppose that it holds for $|\textbf{i}|\leq k$, where $k\geq 0$. It suffices to prove (\ref{1.26}) for $r \in\ove{1,l}$, $\be\in\G\backslash\{-\mu\}$ and $\textbf{i}\in \mbb{N}^{l_1+l_2}$ with $|\textbf{i}|= k+1$.

For any $\be\in\G\backslash\{-\mu\}$ and $\textbf{i}\in \mbb{N}^{l_1+l_2}$ with $|\textbf{i}|= k+1$, we define $p_{\textbf{i}}$, $P(\textbf{i})$ and $Q(\be,\textbf{i})$ as in (\ref{1.20})--(\ref{1.22}). Then we give the proof in three cases.\vspace{0.2cm}

Case 1. $Q(\be,\textbf{i})\not= \emptyset $.\vspace{0.2cm}

Let $q(\be,\textbf{i})=\min Q(\be,\textbf{i})$. Then (\ref{1.23}), Lemma \ref{le:3.6} and the induction hypothesis give
\begin{eqnarray}
\ptl_r.v_{\be,\textbf{i}}
&=&\frac{1}{\bar{\be}_{q(\be,\textbf{i})}}\big( t^{1_{[p_{\textbf{i}}]}}\ptl_{q(\be,\textbf{i})}. \ptl_r.
v_{\be,\textbf{i}-1_{[p_{\textbf{i}}]}}
+\delta_{r,p_{\textbf{i}}}\ptl_{q(\be,\textbf{i})}. v_{\be,\textbf{i}-1_{[p_{\textbf{i}}]}}
-i_{q(\be,\textbf{i})}\ptl_r.v_{\be,\textbf{i}-1_{[q(\be,\textbf{i})]}}
\big)\nonumber\\
&=& \frac{1}{\bar{\be}_{q(\be,\textbf{i})}}\big( t^{1_{[p_{\textbf{i}}]}}\ptl_{q(\be,\textbf{i})}. (\bar{\be}_r v_{\be,\textbf{i}-1_{[p_{\textbf{i}}]}}
+ (i_r-\delta_{r,p_{\textbf{i}}}) v_{\be,\textbf{i}-1_{[p_{\textbf{i}}]}-1_{[r]}})\nonumber\\
& & +\delta_{r,p_{\textbf{i}}}(\bar{\be}_{q(\be,\textbf{i})} v_{\be,\textbf{i}-1_{[p_{\textbf{i}}]}}+ i_{q(\be,\textbf{i})} v_{\be,\textbf{i}-1_{[p_{\textbf{i}}]}-1_{[q(\be,\textbf{i})]}})
 \nonumber\\
& &-i_{q(\be,\textbf{i})}(\bar{\be}_r v_{\be,\textbf{i}-1_{[q(\be,\textbf{i})]}}+ (i_r-\delta_{r,q(\be,\textbf{i})}) v_{\be,\textbf{i}-1_{[q(\be,\textbf{i})]}-1_{[r]}})\big)\nonumber\\
&=& \bar{\be}_r v_{\be,\textbf{i}}+ i_r v_{\be,\textbf{i}-1_{[r]}}  \label{1.28}
\end{eqnarray}
for $ r \in\ove{1,l}$,
which coincides with (\ref{1.26}).\vspace{0.2cm}

Case 2. $Q(\be,\textbf{i})= \emptyset $ and $P(\textbf{i})\not= \emptyset $.\vspace{0.2cm}

Let $p(\textbf{i})=\min P(\textbf{i})$.
Then (\ref{1.24}), Lemma \ref{le:3.6} and the induction hypothesis imply
\begin{eqnarray}
\ptl_r.v_{\be,\textbf{i}}
&=& \frac{1}{\bar{\be}_{p_{\textbf{i}}}} \big(t^{1_{[p(\textbf{i})]}}\ptl_{p_{\textbf{i}}}. \ptl_r.
v_{\be,\textbf{i}-1_{[p(\textbf{i})]}}
+\delta_{r,p(\textbf{i})}\ptl_{p_{\textbf{i}}}.v_{\be,\textbf{i}-1_{[p(\textbf{i})]}}
-i_{p_{\textbf{i}}} \ptl_r. v_{\be,\textbf{i}-1_{[p_{\textbf{i}}]}}\big)\nonumber\\
&=& \frac{1}{\bar{\be}_{p_{\textbf{i}}}}\big( t^{1_{[p(\textbf{i})]}}\ptl_{p_{\textbf{i}}}. (\bar{\be}_r v_{\be,\textbf{i}-1_{[p(\textbf{i})]}}
+ (i_r-\delta_{r,p(\textbf{i})}) v_{\be,\textbf{i}-1_{[p(\textbf{i})]}-1_{[r]}})\nonumber\\
& &+\delta_{r,p(\textbf{i})}(\bar{\be}_{p_{\textbf{i}}} v_{\be,\textbf{i}-1_{[p(\textbf{i})]}}+ i_{p_{\textbf{i}}} v_{\be,\textbf{i}-1_{[p(\textbf{i})]}-1_{[p_{\textbf{i}}]}})
 \nonumber\\
& &-i_{p_{\textbf{i}}}(\bar{\be}_r v_{\be,\textbf{i}-1_{[p_{\textbf{i}}]}}+ (i_r-\delta_{r,p_{\textbf{i}}}) v_{\be,\textbf{i}-1_{[p_{\textbf{i}}]}-1_{[r]}})\big)\nonumber\\
&=& \bar{\be}_r v_{\be,\textbf{i}}+ i_r v_{\be,\textbf{i}-1_{[r]}}\label{1.32}
\end{eqnarray}
for $ r \in\ove{1,l}$,
which coincides with (\ref{1.26}).\vspace{0.2cm}

Case 3. $Q(\be,\textbf{i})= \emptyset $ and $P(\textbf{i})= \emptyset $.\vspace{0.2cm}

Let $s(\textbf{i})=\min \{\ove{1,l_1+l_2}\backslash\{p_{\textbf{i}}\}\}$.
Note that $(\be,\textbf{i}-1_{[p_{\textbf{i}}]}+1_{[s(\textbf{i})]})$ belongs to Case 1 or Case 2. So we have
\begin{eqnarray}
& & \ptl_r.v_{\be,\textbf{i}-1_{[p_{\textbf{i}}]}+1_{[s(\textbf{i})]}}\nonumber\\
 &=& \bar{\be}_r v_{\be,\textbf{i}-1_{[p_{\textbf{i}}]}+1_{[s(\textbf{i})]}}+ (i_r -\delta_{r,p_{\textbf{i}}}+\delta_{r,s(\textbf{i})} ) v_{\be,\textbf{i}-1_{[p_{\textbf{i}}]}+1_{[s(\textbf{i})]}-1_{[r]}}\label{1.39}
\end{eqnarray}
for $r \in\ove{1,l}$. Since $\bar{\be}_s=0$ and $i_s=0$ for $s \in\ove{1,l}\backslash\{p_{\textbf{i}}\}$, from 
(\ref{1.25}), (\ref{1.39}) and Lemma \ref{le:3.6} we derive
\begin{eqnarray}
\ptl_r.v_{\be,\textbf{i}}
&=&t^{1_{[p_{\textbf{i}}]}}\ptl_{s(\textbf{i})}.\ptl_r.
v_{\be,\textbf{i}-1_{[p_{\textbf{i}}]}+1_{[s(\textbf{i})]}}
+  \delta_{r,p_{\textbf{i}}}\ptl_{s(\textbf{i})}. v_{\be,\textbf{i}-1_{[p_{\textbf{i}}]}+1_{[s(\textbf{i})]}}\nonumber\\
&=&t^{1_{[p_{\textbf{i}}]}}\ptl_{s(\textbf{i})}.(\bar{\be}_r v_{\be,\textbf{i}-1_{[p_{\textbf{i}}]}+1_{[s(\textbf{i})]}}+ (i_r -\delta_{r,p_{\textbf{i}}}+\delta_{r,s(\textbf{i})}) v_{\be,\textbf{i}-1_{[p_{\textbf{i}}]}+1_{[s(\textbf{i})]}-1_{[r]}})\nonumber\\
& &+\delta_{r,p_{\textbf{i}}}v_{\be,\textbf{i}-1_{[p_{\textbf{i}}]}}\nonumber\\
&=& \bar{\be}_r v_{\be,\textbf{i}}+ i_r v_{\be,\textbf{i}-1_{[r]}}\label{1.38}
\end{eqnarray}
for $ r \in\ove{1,l}$, which coincides with (\ref{1.26}). 

So this lemma holds. $\qquad\Box$

\begin{lemma}\label{le:3.9}
For any $\al\in\G\backslash\{0\}$, $\be\in\G\backslash\{-\mu, -\mu-\al\}$, $\textbf{i}\in \mbb{N}^{l_1+l_2}$ and $p,q \in\ove{1,l}$ with $p\not=q$, we have
\begin{eqnarray}
 D_{p,q}(x^{\al}).v_{\be,\textbf{i}}
& = & x^{\al}(\al_p\ptl_q-\al_q\ptl_p).v_{\be,\textbf{i}}\nonumber\\
&= & (\al_p(\be_q+\mu_q)-\al_q(\be_p+\mu_p))v_{\be+\al,\textbf{i}}+ \al_p i_q v_{\be+\al,\textbf{i}-1_{[q]}} \nonumber\\
& & -\al_q i_p v_{\be+\al,\textbf{i}-1_{[p]}}.\label{1.58}
\end{eqnarray}

\end{lemma}

\noindent{\bf Proof.} Let $\bar{\ga}=\ga+\mu$ for all $\ga\in\G$. We give the proof by induction on $|\textbf{i}|$. Since Lemma \ref{le:3.2} has given the proof for $|\textbf{i}|=0$, supposing that this lemma holds for $\textbf{i}\in \mbb{N}^{l_1+l_2}$ with $|\textbf{i}|\leq k$, where $k \geq 0$, we only need to prove (\ref{1.58}) for $\textbf{i}\in \mbb{N}^{l_1+l_2}$ with $|\textbf{i}| = k+1$. Pick any $\al\in\G\backslash\{0\}$, $\be\in\G\backslash\{-\mu, -\mu-\al\}$, $\textbf{i}\in \mbb{N}^{l_1+l_2}$ with $|\textbf{i}|= k+1$, and $p,q \in\ove{1,l}$ with $p\not=q$. First of all, we define $p_{\textbf{i}}$ as in (\ref{1.20}), namely,
\begin{equation}
p_{\textbf{i}}=\min\{r\in\ove{1,l_1+l_2}\mid i_r\not=0\}.
\end{equation}
 Then we proceed the proof in several cases.\vspace{0.2cm}

Case 1. $\ker \al \bigcap \big(\sum_{j\in\ove{1,l}\backslash\{p_{\textbf{i}}\}}\mbb{F}\ptl_j\big)\backslash\ker\bar{\be}\not=\emptyset$.\vspace{0.2cm}

If there exists $s\in\ove{1,l}\backslash\{p_{\textbf{i}}\}$ such that $\al_s\not=0$, then $\al_s\bar{\be}_r-\al_r\bar{\be}_s\not=0$ for some $r\in\ove{1,l}\backslash\{s,p_{\textbf{i}}\}$. Fix such $s$ and $r$. Since Lemma \ref{le:3.6} gives
\begin{eqnarray}
& & t^{1_{[p_{\textbf{i}}]}}(\al_s\ptl_r-\al_r\ptl_s).v_{\be,\textbf{i}-1_{[p_{\textbf{i}}]}}\nonumber\\
& = & (\al_s\bar{\be}_r-\al_r\bar{\be}_s)v_{\be,\textbf{i}}+\al_s i_r v_{\be,\textbf{i}-1_{[r]}}
-\al_r i_s v_{\be,\textbf{i}-1_{[s]}},\label{1.120}
\end{eqnarray}
we have
\begin{eqnarray}
 D_{p,q}(x^{\al}).v_{\be,\textbf{i}}
& = & \frac{1}{\al_s\bar{\be}_r-\al_r\bar{\be}_s}x^{\al}(\al_p\ptl_q-\al_q\ptl_p).\big(t^{1_{[p_{\textbf{i}}]}}(\al_s\ptl_r-\al_r\ptl_s).v_{\be,\textbf{i}-1_{[p_{\textbf{i}}]}}
\nonumber\\
& &-\al_s i_r v_{\be,\textbf{i}-1_{[r]}}
+\al_r i_s v_{\be,\textbf{i}-1_{[s]}}\big)\nonumber\\
& = & \frac{1}{\al_s\bar{\be}_r-\al_r\bar{\be}_s} \Big( t^{1_{[p_{\textbf{i}}]}}(\al_s\ptl_r-\al_r\ptl_s).
x^{\al}(\al_p\ptl_q-\al_q\ptl_p).v_{\be,\textbf{i}-1_{[p_{\textbf{i}}]}}
\nonumber\\
& & +(\al_p\delta_{q,p_{\textbf{i}}}-\al_q\delta_{p,p_{\textbf{i}}})x^{\al}(\al_s\ptl_r-\al_r\ptl_s).v_{\be,\textbf{i}-1_{[p_{\textbf{i}}]}}
\nonumber\\
& & -\al_s i_r x^{\al}(\al_p\ptl_q-\al_q\ptl_p).v_{\be,\textbf{i}-1_{[r]}}
+\al_r i_s x^{\al}(\al_p\ptl_q-\al_q\ptl_p).v_{\be,\textbf{i}-1_{[s]}}\Big)\nonumber\\
&= & (\al_p \bar{\be}_q-\al_q \bar{\be}_p)v_{\be+\al,\textbf{i}}+ \al_p i_q v_{\be+\al,\textbf{i}-1_{[q]}} -\al_q i_p v_{\be+\al,\textbf{i}-1_{[p]}}
\end{eqnarray}
by (\ref{1.120}), Lemma \ref{le:3.6} and the induction hypothesis. So (\ref{1.58}) holds.

If $\al_s=0$ for all $s\in\ove{1,l}\backslash\{p_{\textbf{i}}\}$, then $\bar{\be}_r\not=0$ for some $r\in\ove{1,l}\backslash\{p_{\textbf{i}}\}$ because that $\ker \al \bigcap \big(\sum_{j\in\ove{1,l}\backslash\{p_{\textbf{i}}\}}\mbb{F}\ptl_j\big)\backslash\ker\bar{\be}\not=\emptyset$. Fix such $r$. Since Lemma \ref{le:3.6} indicates
\begin{equation}\label{1.121}
t^{1_{[p_{\textbf{i}}]}}\ptl_r.v_{\be,\textbf{i}-1_{[p_{\textbf{i}}]}}= \bar{\be}_r v_{\be,\textbf{i}}+ i_r v_{\be,\textbf{i}-1_{[r]}},
\end{equation}
we obtain
\begin{eqnarray}
& & D_{p,q}(x^{\al}).v_{\be,\textbf{i}}\nonumber\\
& = & \frac{1}{\bar{\be}_r}x^{\al}(\al_p\ptl_q-\al_q\ptl_p).(t^{1_{[p_{\textbf{i}}]}}\ptl_r.v_{\be,\textbf{i}-1_{[p_{\textbf{i}}]}}
- i_r v_{\be,\textbf{i}-1_{[r]}})\nonumber\\
& = & \frac{1}{\bar{\be}_r} \Big( t^{1_{[p_{\textbf{i}}]}}\ptl_r.
x^{\al}(\al_p\ptl_q-\al_q\ptl_p).v_{\be,\textbf{i}-1_{[p_{\textbf{i}}]}} +(\al_p\delta_{q,p_{\textbf{i}}}-\al_q\delta_{p,p_{\textbf{i}}})x^{\al}\ptl_r.v_{\be,\textbf{i}-1_{[p_{\textbf{i}}]}}
\nonumber\\
& & -  i_r x^{\al}(\al_p\ptl_q-\al_q\ptl_p).v_{\be,\textbf{i}-1_{[r]}}\Big)\nonumber\\
&= & (\al_p \bar{\be}_q-\al_q \bar{\be}_p)v_{\be+\al,\textbf{i}}+ \al_p i_q v_{\be+\al,\textbf{i}-1_{[q]}} -\al_q i_p v_{\be+\al,\textbf{i}-1_{[p]}}
\end{eqnarray}
by (\ref{1.121}), Lemma \ref{le:3.6} and the induction hypothesis. So (\ref{1.58}) holds.
\vspace{0.2cm}

Case 2.  $\ker \al \bigcap \big(\sum_{j\in\ove{1,l}\backslash\{p_{\textbf{i}}\}}\mbb{F}\ptl_j\big)\subseteq \ker\bar{\be}$ and $\ker \al\not=\ker\bar{\be}$.\vspace{0.2cm}

Since $\al\not=0$ and $\bar{\be}\not=0\not=\bar{\be}+\al$, we have $\al_s\not=0$ for some $s\in\ove{1,l}\backslash\{p_{\textbf{i}}\}$, and $\al_s\bar{\be}_{p_{\textbf{i}}}-\al_{p_{\textbf{i}}}\bar{\be}_s\not=0$ in this case. Fix such $s$. Pick $r\in\ove{1,l_1+l_2}\backslash\{p_{\textbf{i}},s\}$. Set $\tilde{\ptl}_1=\al_s\ptl_{p_{\textbf{i}}}-\al_{p_{\textbf{i}}}\ptl_s$ and $\tilde{\ptl}_2=\al_s\ptl_r-\al_r\ptl_s$. Observe that
\begin{equation}
\tilde{\ptl}_1(t^{1_{[r]}+1_{[p_{\textbf{i}}]}}) \tilde{\ptl}_2-\tilde{\ptl}_2(t^{1_{[r]}+1_{[p_{\textbf{i}}]}}) \tilde{\ptl}_1
=\al_s (t^{1_{[r]}}\tilde{\ptl}_2-  t^{1_{[p_{\textbf{i}}]}} \tilde{\ptl}_1).
\end{equation}
So Lemma \ref{le:3.6} and Corollary \ref{c:3.7} imply
\begin{eqnarray}
& &(t^{1_{[p_{\textbf{i}}]}} \tilde{\ptl}_1 -t^{1_{[r]}}\tilde{\ptl}_2 ). v_{\be,\textbf{i}-1_{[p_{\textbf{i}}]}}\nonumber\\
&=&(\al_s\bar{\be}_{p_{\textbf{i}}}-\al_{p_{\textbf{i}}}\bar{\be}_s)
v_{\be,\textbf{i}}+ \al_s (i_{p_{\textbf{i}}}-1) v_{\be,\textbf{i}-1_{[p_{\textbf{i}}]}} -\al_{p_{\textbf{i}}}i_s v_{\be,\textbf{i}-1_{[s]}}\nonumber\\
& & -\al_s i_r v_{\be,\textbf{i}-1_{[p_{\textbf{i}}]}} + \al_r i_s v_{\be,\textbf{i}+1_{[r]}-1_{[p_{\textbf{i}}]}-1_{[s]}}.\label{1.122}
\end{eqnarray}
Thus, we have
\begin{eqnarray}
& & D_{p,q}(x^{\al}).v_{\be,\textbf{i}}\nonumber\\
& = & \frac{1}{\al_s\bar{\be}_{p_{\textbf{i}}}-\al_{p_{\textbf{i}}}\bar{\be}_s}x^{\al}(\al_p\ptl_q-\al_q\ptl_p).
\Big((t^{1_{[p_{\textbf{i}}]}} \tilde{\ptl}_1 -t^{1_{[r]}}\tilde{\ptl}_2 ). v_{\be,\textbf{i}-1_{[p_{\textbf{i}}]}}\nonumber\\
& &- \al_s (i_{p_{\textbf{i}}}-1) v_{\be,\textbf{i}-1_{[p_{\textbf{i}}]}}+\al_{p_{\textbf{i}}}i_s v_{\be,\textbf{i}-1_{[s]}} +\al_s i_r v_{\be,\textbf{i}-1_{[p_{\textbf{i}}]}}  -\al_r i_s v_{\be,\textbf{i}+1_{[r]}-1_{[p_{\textbf{i}}]}-1_{[s]}}\Big)\nonumber\\
& = & \frac{1}{\al_s\bar{\be}_{p_{\textbf{i}}}-\al_{p_{\textbf{i}}}\bar{\be}_s}
\Big((t^{1_{[p_{\textbf{i}}]}} \tilde{\ptl}_1 -t^{1_{[r]}}\tilde{\ptl}_2 ). x^{\al}(\al_p\ptl_q-\al_q\ptl_p). v_{\be,\textbf{i}-1_{[p_{\textbf{i}}]}}\nonumber\\
& & + (\al_p\delta_{q,p_{\textbf{i}}}-\al_q\delta_{p,p_{\textbf{i}}})
x^{\al}\tilde{\ptl}_1.v_{\be,\textbf{i}-1_{[p_{\textbf{i}}]}} - (\al_p\delta_{q,r}-\al_q\delta_{p,r})x^{\al}\tilde{\ptl}_2 . v_{\be,\textbf{i}-1_{[p_{\textbf{i}}]}}\nonumber\\
& & -x^{\al}(\al_p\ptl_q-\al_q\ptl_p).
( \al_s (i_{p_{\textbf{i}}}-1) v_{\be,\textbf{i}-1_{[p_{\textbf{i}}]}}- \al_{p_{\textbf{i}}}i_s v_{\be,\textbf{i}-1_{[s]}} -\al_s i_r v_{\be,\textbf{i}-1_{[p_{\textbf{i}}]}} \nonumber\\
& & +\al_r i_s v_{\be,\textbf{i}+1_{[r]}-1_{[p_{\textbf{i}}]}-1_{[s]}})\Big)\nonumber\\
&= & (\al_p \bar{\be}_q-\al_q \bar{\be}_p)v_{\be+\al,\textbf{i}}+ \al_p i_q v_{\be+\al,\textbf{i}-1_{[q]}} -\al_q i_p v_{\be+\al,\textbf{i}-1_{[p]}}
\end{eqnarray}
by (\ref{1.122}), Lemma \ref{le:3.6} and the induction hypothesis. So (\ref{1.58}) holds.\vspace{0.2cm}

Case 3.  $\ker \al =\ker\bar{\be}$.\vspace{0.2cm}

Since $\al\not=0$, $\bar{\be}\not=0\not=\bar{\be}+\al$ and $\ker \al =\ker\bar{\be}$, we have $\al_s\not=0$, $\bar{\be}_s\not=0$ and $\bar{\be}_s+\al_s\not=0$ for some $s\in\ove{1,l}$. Fix such $s$.

If there exists some $r\in\ove{1,l}\backslash\{s\}$ such that $\bar{\be}_r=0$, then $\al_r=0$. Fix such $r$.
Choose $\sigma\in\G\backslash\{0\}$ such that $\sigma_r\not=0$. Then $\ker \sigma \not=\ker(\bar{\be}-\sigma)$. So Case 1 or Case 2 gives
\begin{equation}\label{1.59}
x^{\sigma}(\sigma_r\ptl_s-\sigma_s\ptl_r). v_{\be-\sigma,\textbf{i}}
= \sigma_r \bar{\be}_s v_{\be,\textbf{i}} + \sigma_r i_s v_{\be,\textbf{i}-1_{[s]}}- \sigma_s i_r v_{\be,\textbf{i}-1_{[r]}}.
\end{equation}
Thus we have
\begin{eqnarray}
& & D_{p,q}(x^{\al}).v_{\be,\textbf{i}}\nonumber\\
& = & \frac{1}{\sigma_r\bar{\be}_s}x^{\al}(\al_p\ptl_q-\al_q\ptl_p).
(x^{\sigma}(\sigma_r\ptl_s-\sigma_s\ptl_r). v_{\be-\sigma,\textbf{i}}
-  \sigma_r i_s v_{\be,\textbf{i}-1_{[s]}}+\sigma_s i_r v_{\be,\textbf{i}-1_{[r]}})\nonumber\\
& = & \frac{1}{\sigma_r\bar{\be}_s}
(x^{\sigma}(\sigma_r\ptl_s-\sigma_s\ptl_r).x^{\al}(\al_p\ptl_q-\al_q\ptl_p). v_{\be-\sigma,\textbf{i}}\nonumber\\
& &+x^{\al+\sigma}((\al_p\sigma_q-\al_q\sigma_p)(\sigma_r\ptl_s-\sigma_s\ptl_r)
-(\sigma_r\al_s-\sigma_s\al_r)(\al_p\ptl_q-\al_q\ptl_p)).v_{\be-\sigma,\textbf{i}}\nonumber\\
& & - \sigma_r i_s x^{\al}(\al_p\ptl_q-\al_q\ptl_p).v_{\be,\textbf{i}-1_{[s]}}+\sigma_s i_r x^{\al}(\al_p\ptl_q-\al_q\ptl_p).v_{\be,\textbf{i}-1_{[r]}}) \label{1.63}
\end{eqnarray}
by (\ref{1.59}). Moreover, since $\ker \al \not=\ker(\bar{\be}-\sigma)$, $\ker(\al+ \sigma) \not=\ker(\bar{\be}-\sigma)$ and $\ker \sigma \not=\ker(\bar{\be}+\al-\sigma)$, Case 1 and Case 2 give
\begin{eqnarray}
& &x^{\al}(\al_p\ptl_q-\al_q\ptl_p). v_{\be-\sigma,\textbf{i}}\nonumber\\
&=& (\al_q\sigma_p-\al_p\sigma_q) v_{\be+\al-\sigma,\textbf{i}} + \al_p i_q v_{\be+\al-\sigma,\textbf{i}-1_{[q]}}
- \al_q i_p v_{\be+\al-\sigma,\textbf{i}-1_{[p]}},\label{1.60}
\end{eqnarray}
\begin{eqnarray}
& &x^{\sigma}(\sigma_r\ptl_s-\sigma_s\ptl_r). v_{\be+\al-\sigma,\textbf{i}}\nonumber\\
&=& \sigma_r(\bar{\be}_s+\al_s) v_{\be+\al,\textbf{i}} + \sigma_r i_s v_{\be+\al,\textbf{i}-1_{[s]}} - \sigma_s i_r v_{\be+\al,\textbf{i}-1_{[r]}}\label{1.61}
\end{eqnarray}
and
\begin{eqnarray}
& &x^{\al+\sigma}((\al_p\sigma_q-\al_q\sigma_p)(\sigma_r\ptl_s-\sigma_s\ptl_r)
-(\sigma_r\al_s-\sigma_s\al_r)(\al_p\ptl_q-\al_q\ptl_p)).v_{\be-\sigma,\textbf{i}}\nonumber\\
&=& \sigma_r(\bar{\be}_s+\al_s)(\al_p\sigma_q-\al_q\sigma_p) v_{\be+\al,\textbf{i}}
+\sigma_r i_s(\al_p\sigma_q-\al_q\sigma_p)v_{\be+\al,\textbf{i}-1_{[s]}}\nonumber\\
& &
-\sigma_s i_r(\al_p\sigma_q-\al_q\sigma_p)v_{\be+\al,\textbf{i}-1_{[r]}}
-\al_p i_q(\sigma_r\al_s-\sigma_s\al_r)v_{\be+\al,\textbf{i}-1_{[q]}}\nonumber\\
& &
+\al_q i_p(\sigma_r\al_s-\sigma_s\al_r)v_{\be+\al,\textbf{i}-1_{[p]}}.\label{1.62}
\end{eqnarray}
So (\ref{1.63}) becomes
\begin{eqnarray}
& & D_{p,q}(x^{\al}).v_{\be,\textbf{i}}\nonumber\\
& = & \al_p i_q v_{\be+\al,\textbf{i}-1_{[q]}} -\al_q i_p v_{\be+\al,\textbf{i}-1_{[p]}}\nonumber\\
& = & (\al_p \bar{\be}_q-\al_q \bar{\be}_p)v_{\be+\al,\textbf{i}}+ \al_p i_q v_{\be+\al,\textbf{i}-1_{[q]}} -\al_q i_p v_{\be+\al,\textbf{i}-1_{[p]}}\label{1.64}
\end{eqnarray}
by (\ref{1.60})--(\ref{1.62}) and the induction hypothesis, which coincides with (\ref{1.58}).

If $\bar{\be}_r\not=0$ for all $r\in\ove{1,l}\backslash\{s\}$, which also means $\al_r\not=0$ for all $r\in\ove{1,l}\backslash\{s\}$, we pick $r_1,r_2\in\ove{1,l}\backslash\{s\}$. Then replacing $\ptl_r$ by $\ptl=\al_{r_1}\ptl_{r_2}-\al_{r_2}\ptl_{r_1}$ and $\sigma_r$ by $\ptl(\sigma)$ respectively in the above discussion from (\ref{1.59}) to (\ref{1.64}), we can also get
\begin{eqnarray}
& & D_{p,q}(x^{\al}).v_{\be,\textbf{i}}\nonumber\\
& = & (\al_p \bar{\be}_q-\al_q \bar{\be}_p)v_{\be+\al,\textbf{i}}+ \al_p i_q v_{\be+\al,\textbf{i}-1_{[q]}} -\al_q i_p v_{\be+\al,\textbf{i}-1_{[p]}}.
\end{eqnarray}

Thus we complete the proof of this lemma. $\qquad\Box$ \vspace{0.3cm}

Next, we shall prove the main theorem for the case $\mu\in\mbb{F}^{l_2+l_3}\backslash\G$.

For convenience, we shall first give a total order on $\mbb{N}^{l_1+l_2}$:\vspace{0.3cm}
\begin{defi}\label{d:1.5}
We define a total order on $\mbb{N}^{l_1+l_2}$ by:
\begin{equation}\label{2.1}
\textbf{i}>\textbf{j} \Longleftrightarrow |\textbf{i}|>|\textbf{j}|, \textrm{ or, } |\textbf{i}|=|\textbf{j}| \textrm{ and } i_s>j_s  \textrm{ with } i_p=j_p \textrm{ for }p\in\ove{s+1,l_1+l_2}.
\end{equation}
\end{defi}

Then we have:\vspace{0.3cm}

\begin{lemma}\label{le:3.10}
If $\mu\not\in\G$, then $\{v_{\be,\textbf{i}}\mid \be\in\G, \textbf{i}\in \mbb{N}^{l_1+l_2}\}$ is an $\mbb{F}$-basis of $V$.
\end{lemma}

\noindent{\bf Proof.} It is straightforward to prove that $\{v_{\be,\textbf{i}}\mid \be\in\G, \textbf{i}\in \mbb{N}^{l_1+l_2} \}$ is a linearly independent set by Lemma \ref{le:3.8}. We omit the details.
We only prove that $V$ is spanned by $\{v_{\be,\textbf{i}}\mid \be\in\G, \textbf{i}\in \mbb{N}^{l_1+l_2}\}$ here. For any $\be\in\G$ and $\textbf{i}\in \mbb{N}^{l_1+l_2}$, we set
\begin{equation}\label{1.123}
V_{\be+\mu}^{[\textbf{i}]}=\{v\in V\mid \prod_{p=1}^l (\ptl_p-(\be_p+\mu_p))^{j_p} (v)=0 \textrm{ for }\textbf{j}\in \mbb{N}^{l_1+l_2} \textrm{ with }\textbf{j}>\textbf{i}\},
\end{equation}
\begin{equation}
V_{\be+\mu}^{(\textbf{i})}=\bigcup_{\textbf{j}\in \mbb{N}^{l_1+l_2}; \textbf{j}<\textbf{i}} V_{\be+\mu}^{[\textbf{j}]}
\end{equation}
and
\begin{equation}
V^{[\textbf{i}]}=\bigoplus_{\be\in\G}V_{\be+\mu}^{[\textbf{i}]},\quad V^{(\textbf{i})}=\bigoplus_{\be\in\G}V_{\be+\mu}^{(\textbf{i})}.
\end{equation}
Then Lemma \ref{le:3.8} implies
\begin{equation}\label{1.126}
V_{\be+\mu}=\bigcup_{\textbf{j}\in \mbb{N}^{l_1+l_2}} V_{\be+\mu}^{[\textbf{j}]} \textrm{ and } v_{\be,\textbf{i}}\in V_{\be+\mu}^{[\textbf{i}]}\backslash V_{\be+\mu}^{(\textbf{i})},\quad \forall \be\in\G,\; \textbf{i}\in \mbb{N}^{l_1+l_2}.
\end{equation}
Suppose that $V$ cannot be spanned by $\{v_{\be,\textbf{i}}\mid \be\in\G, \textbf{i}\in \mbb{N}^{l_1+l_2}\}$. Then this will lead to a contradiction. Let
 $\textbf{j}\in\mbb{N}^{l_1+l_2}$ be the minimal element such that
\begin{eqnarray}
& \textrm{ there exist } \al\in\G \textrm{ and } v\in V_{\al+\mu}^{[\textbf{j}]}\backslash V_{\al+\mu}^{(\textbf{j})} \textrm{  such that}&\nonumber\\
&  v\not\in \textrm{Span}_{\mbb{F}}\{v_{\be,\textbf{i}}\mid \be\in\G, \textbf{i}\in \mbb{N}^{l_1+l_2}\}.&\label{1.124}
\end{eqnarray}
 Since $V_{\be+\mu}^{[\textbf{0}]}=V_{\be+\mu}^{(0)}=\mbb{F}v_{\be,\textbf{0}}$ for $\be\in\G$ (cf. (\ref{2.2}), Definition \ref{d:3.3}), we have $\textbf{j}\not=\textbf{0}$. Let $q=\min\{p\in\ove{1,l_1+l_2}\mid j_p\not=0\}$. Then (\ref{1.123}) and Lemma \ref{le:3.8} show
\begin{equation}\label{1.127}
v'=(\ptl_q-(\al_q+\mu_q))v \in V_{\al+\mu}^{[\textbf{j}-1_{[q]}]}.
\end{equation}
By the minimality of $\textbf{j}$ in (\ref{1.124}), we have
\begin{equation}\label{1.125}
v' \in \textrm{Span}_{\mbb{F}}\{v_{\be,\textbf{i}}\mid \be\in\G, \textbf{i}\in \mbb{N}^{l_1+l_2}\}.
\end{equation}
Thus (\ref{1.123}), (\ref{1.126}), (\ref{1.127}) and (\ref{1.125}) give
\begin{equation}\label{1.128}
v''=v'-c v_{\al,\textbf{j}-1_{[q]}} \in V_{\al+\mu}^{(\textbf{j}-1_{[q]})} \textrm{ for some } c\in\mbb{F}.
\end{equation}
Let
\begin{equation}\label{1.129}
w=v -\frac{c}{j_q} v_{\al,\textbf{j}}.
\end{equation}
Then $w\in V_{\al+\mu}^{[\textbf{j}]}$ by (\ref{1.126}) and (\ref{1.124}). So for $\textbf{i}>\textbf{j}$, we have
\begin{equation}\label{1.130}
\prod_{p=1}^l (\ptl_p-(\al_p+\mu_p))^{i_p} (w)=0
\end{equation}
by (\ref{1.123}). Moreover, since (\ref{1.127}), (\ref{1.128}), (\ref{1.129}) and Lemma \ref{le:3.8} give
\begin{equation}
(\ptl_q-(\al_q+\mu_q))(w)=v'-c v_{\al,\textbf{j}-1_{[q]}} =v'',
\end{equation}
 we have
\begin{equation}\label{1.131}
\prod_{p=1}^l (\ptl_p-(\al_p+\mu_p))^{j_p} (w)=\prod_{p=1}^l (\ptl_p-(\al_p+\mu_p))^{j_p-\delta_{p,q}} (v'')=0
\end{equation}
by (\ref{1.128}). So (\ref{1.130}) and (\ref{1.131}) show that
\begin{equation}
w=v -\frac{c}{j_q} v_{\al,\textbf{j}}\in  V_{\al+\mu}^{(\textbf{j})},
\end{equation}
which indicates $v\in \textrm{Span}_{\mbb{F}}\{v_{\be,\textbf{i}}\mid \be\in\G, \textbf{i}\in \mbb{N}^{l_1+l_2}\}$. This contradicts (\ref{1.124}). So $V$ is spanned by $\{v_{\be,\textbf{i}}\mid \be\in\G, \textbf{i}\in \mbb{N}^{l_1+l_2}\}$. This lemma holds. $\qquad\Box$ \vspace{0.3cm}

For the case $\mu\not\in\G$, we have determined a basis of $V=V(\mu)$, and derived the action of the set (\ref{1.1}) on the basis. So Proposition \ref{t:2.1}, Lemma \ref{le:3.6}, Lemma \ref{le:3.8}, Lemma \ref{le:3.9}, Corollary \ref{c:3.7} and Lemma \ref{le:3.10} give

\begin{lemma}\label{le:3.11}
If $\mu\not\in\G$, then $V=V(\mu)\simeq A_{\mu}$.
\end{lemma}

\section{Proof of the main theorem (\Rmnum{2})}

In this section, we consider the case that $V=V(\mu)$ for some $\mu\in\G$. With a shift of the indices, we can always assume that $\mu=0$. Notice that the contents from Lemma \ref{le:3.2} to Lemma \ref{le:3.9} were discussed for general case $\mu\in\mbb{F}^{l_2+l_3}$. So they still hold for the case $\mu=0$. In this section, we need to complement the basis of $V=V(\mu)$, and to derive the action of the set (\ref{1.1}) on the basis which were missed from Lemma \ref{le:3.2} to Lemma \ref{le:3.9} under the condition $\mu=0$.

\begin{lemma}\label{le:4.1}
 $V_0\not=\{0\}$.
\end{lemma}

\noindent{\bf Proof.} Suppose that $V_0=\{0\}$. Then this will lead to a contradiction. Pick $p\in\ove{1,l_1+l_2}$ and $q\in\ove{l_1+1,l}\backslash\{p\}$. Choose $\rho\in\G\backslash\{0\}$ such that $\rho_q\not=0$. Then by Lemmas \ref{le:3.2} and \ref{le:3.8}, we have
\begin{equation}
\ptl_s.(x^{\rho}(\rho_p\ptl_q-\rho_q\ptl_p).v_{-\rho,2_{[p]}})=0 \quad \textrm{ for }s\in\ove{1,l}\backslash\{p\},
\end{equation}
\begin{equation}
\ptl_p^2.( x^{\rho}(\rho_p\ptl_q-\rho_q\ptl_p).v_{-\rho,2_{[p]}})=0.
\end{equation}
So $x^{\rho}(\rho_p\ptl_q-\rho_q\ptl_p).v_{-\rho,2_{[p]}}\in V_0$ by (\ref{2.2}), which indicates
\begin{equation}\label{4.2}
x^{\rho}(\rho_p\ptl_q-\rho_q\ptl_p).v_{-\rho,2_{[p]}}=0.
\end{equation}
Pick $r\in\ove{l_1+1,l}\backslash\{p,q\}$ and choose $\al\in\G\backslash\{0,\rho\}$ such that $\al_r\not=0$.
Then (\ref{4.2}) and Lemma \ref{le:3.9} give
\begin{eqnarray}
0&=&x^{\al}(\al_r\ptl_p-\al_p\ptl_r).x^{\rho}(\rho_p\ptl_q-\rho_q\ptl_p).v_{-\rho,2_{[p]}}\nonumber\\
&=& x^{\al+\rho}((\al_r\rho_p-\al_p\rho_r)(\rho_p\ptl_q-\rho_q\ptl_p)-(\rho_p\al_q-\rho_q\al_p)(\al_r\ptl_p-\al_p\ptl_r)).v_{-\rho,2_{[p]}}\nonumber\\
& & +x^{\rho}(\rho_p\ptl_q-\rho_q\ptl_p).x^{\al}(\al_r\ptl_p-\al_p\ptl_r).v_{-\rho,2_{[p]}}\nonumber\\
&=& -2\rho_q \al_r  v_{\al,\textbf{0}} \not=0,
\end{eqnarray}
which is absurd. So we must have $V_0\not=\{0\}$. This lemma holds. $\qquad\Box$\vspace{0.3cm}

Observe that
\begin{equation}
V_0\not=\{0\}\Leftrightarrow V_0^{(0)}\not=\{0\}.
\end{equation}
Since $\dim V_0^{(0)}\leq 1$, we have $\dim V_0^{(0)}= 1$ by the above lemma. Moreover, we can obtain:

\begin{lemma}\label{le:4.2}
$V_0^{(0)}$ is a trivial $\mcr{S}(l_1,l_2,l_3;0,\G)$-submodule of $V$.
\end{lemma}

\noindent{\bf Proof.} Pick $0\not=v \in V_0^{(0)}$. Then (\ref{2.2}) indicates
\begin{equation}\label{1.65}
\ptl_p.v =0 \quad \textrm{ for } p\in \ove{1,l}.
\end{equation}
Moreover, we obtain
\begin{equation}
t^{1_{[p]}}\ptl_{q'}.v ,\ (t^{1_{[p]}}\ptl_p-t^{1_{[q]}}\ptl_q).v \ \in V_0^{(0)}
\end{equation}
 for any $p,q\in\ove{1,l_1+l_2}$ and $q'\in \ove{1,l}$ with $p\not=q'$ and $p\not=q$ respectively. Namely, they all act on $v$ as scalars.
Thus we can derive, for $p\in\ove{1,l_1+l_2}$ and $q'\in \ove{1,l}$ with $p\not=q'$,
\begin{equation}\label{1.66}
t^{1_{[p]}}\ptl_{q'}.v =[t^{1_{[p]}}\ptl_r,t^{1_{[r]}}\ptl_{q'}].v =0,
\end{equation}
where $r\in\ove{1,l_1+l_2}\backslash\{p,q'\}$, and for $p,q\in\ove{1,l_1+l_2}$ with $p\not=q$,
\begin{equation}\label{1.67}
(t^{1_{[p]}}\ptl_p-t^{1_{[q]}}\ptl_q).v =[t^{1_{[p]}}\ptl_q,t^{1_{[q]}}\ptl_p].v =0.
\end{equation}
Take any $\al\in\G\backslash\{0\}$. Then $\al_r\not=0$ for some $r\in\ove{1,l}$. Fix such $r$. Then for any $s\in\ove{1,l}\backslash\{r\}$, picking $p\in\ove{1,l_1+l_2}\backslash\{r,s\}$, we obtain
\begin{eqnarray}
& & x^{\al}(\al_r\ptl_s-\al_s\ptl_r).v \nonumber\\
&=&\frac{1}{\al_r}[x^{\al}(\al_r\ptl_p-\al_p\ptl_r),t^{1_{[p]}}(\al_r\ptl_s-\al_s\ptl_r)].v \nonumber\\
&=& \frac{1}{\al_r}\Big(x^{\al}(\al_r\ptl_p-\al_p\ptl_r).t^{1_{[p]}}(\al_r\ptl_s-\al_s\ptl_r).v \nonumber\\
& &-t^{1_{[p]}}(\al_r\ptl_s-\al_s\ptl_r).(x^{\al}(\al_r\ptl_p-\al_p\ptl_r).v )\Big)\nonumber\\
&=&0
\end{eqnarray}
by (\ref{1.66}) and Lemma \ref{le:3.6}, where in the fourth line
$x^{\al}(\al_r\ptl_p-\al_p\ptl_r).v \ \in V_\al^{(0)}$.
 So it follows that
\begin{equation}\label{1.69}
D_{p,q}(x^{\al}).v =0 \quad \textrm{ for any }\al\in\G\backslash\{0\},\; p,q\in\ove{1,l}.
\end{equation}
By (\ref{1.65}), (\ref{1.66}), (\ref{1.67}), (\ref{1.69}) and Proposition \ref{t:2.1}, we can derive that $V_0^{(0)}$ is a trivial $\mcr{S}(l_1,l_2,l_3;0,\G)$-submodule of $V$. $\qquad\Box$ \vspace{0.3cm}

For any $p\in\ove{1,l_1+l_2}$, $\al\in\G\backslash\{0\}$ and $\ptl\in\ker\al$, we have
\begin{equation}\label{1.68}
x^{\al}\ptl.v_{-\al,1_{[p]}}\in V_0^{(0)}
\end{equation}
because
\begin{equation}
\ptl_q.(x^{\al}\ptl.v_{-\al,1_{[p]}})=[\ptl_q,x^{\al}\ptl].v_{-\al,1_{[p]}}+x^{\al}\ptl.\ptl_q.v_{-\al,1_{[p]}}
=0 \ \textrm{ for }   q\in\ove{1,l}
\end{equation}
by Lemma \ref{le:3.2} and Lemma \ref{le:3.8}.

Throughout the rest of the section, we fix some
\begin{equation}\label{4.8}
r_1,r_2\in\ove{l_1+1,l}\backslash\{1\} \textrm{ with } r_1\not=r_2,
\end{equation}
 and fix some
\begin{equation}\label{4.9}
\rho\in\G\backslash\{0\} \textrm{ satisfying }\rho_{r_1}\not=0  \textrm{ and }\rho_{{r_2}}\not=0.
 \end{equation}
 We are going to use $\rho,r_1,r_2$ to determine a basis of the vector space $V_0$. And the fixed elements $\rho,r_1,r_2$ will be frequently used throughout the definitions, lemmas and proofs of the rest of the section.
First we have:

\begin{lemma}\label{le:4.3}
 $x^{\rho}(\rho_{r_1}\ptl_1-\rho_1\ptl_{r_1}).v_{-\rho,1_{[1]}}\not=0$, where $\rho$ and $r_1$ are the fixed elements in (\ref{4.8}) and (\ref{4.9}).
\end{lemma}

\noindent{\bf Proof.} Suppose that
\begin{equation}\label{1.72}
x^{\rho}(\rho_{r_1}\ptl_1-\rho_1\ptl_{r_1}).v_{-\rho,1_{[1]}}=0,
\end{equation}
 then we have
\begin{equation}
\ptl_p.(x^{\rho}(\rho_{r_1}\ptl_1-\rho_1\ptl_{r_1}).v_{-\rho,2_{[1]}})=0 \quad \textrm{ for }p\in\ove{1,l}
\end{equation}
by Lemma \ref{le:3.8} and (\ref{1.72}). In other words,
\begin{equation}\label{1.73}
x^{\rho}(\rho_{r_1}\ptl_1-\rho_1\ptl_{r_1}).v_{-\rho,2_{[1]}}\in V_0^{(0)}.
\end{equation}
Choose $\al\in\G\backslash\{0,\rho\}$ such that $\al_{r_1}\not=0$.
So (\ref{1.73}), Lemma \ref{le:4.2} and Lemma \ref{le:3.9} give
\begin{eqnarray}
0&=&x^{\al}(\al_{r_1}\ptl_1-\al_1\ptl_{r_1}).(x^{\rho}(\rho_{r_1}\ptl_1-\rho_1\ptl_{r_1}).v_{-\rho,2_{[1]}})\nonumber\\
&=& x^{\al+\rho}\big((\al_{r_1}\rho_1-\al_1\rho_{r_1})(\rho_{r_1}\ptl_1-\rho_1\ptl_{r_1})-(\rho_{r_1}\al_1-\rho_1\al_{r_1})(\al_{r_1}\ptl_1-\al_1\ptl_{r_1})\big).v_{-\rho,2_{[1]}}\nonumber\\
& & +x^{\rho}(\rho_{r_1}\ptl_1-\rho_1\ptl_{r_1}).x^{\al}(\al_{r_1}\ptl_1-\al_1\ptl_{r_1}).v_{-\rho,2_{[1]}}\nonumber\\
&=& 2\rho_{r_1} \al_{r_1}  v_{\al,\textbf{0}} \not=0,
\end{eqnarray}
which is absurd. So we must have
\begin{equation}
x^{\rho}(\rho_{r_1}\ptl_1-\rho_1\ptl_{r_1}).v_{-\rho,1_{[1]}}\not=0.
\end{equation}
Thus the lemma holds. $\qquad\Box$ \vspace{0.3cm}

Let $\rho$ and $r_1$ be the fixed elements in (\ref{4.8}) and (\ref{4.9}). Since (\ref{1.68}) shows
\begin{equation}
x^{\rho}(\rho_{r_1}\ptl_1-\rho_1\ptl_{r_1}).v_{-\rho,1_{[1]}}\in V_0^{(0)},
\end{equation}
the above lemma enables us to choose $0\not=v_{0,\textbf{0}}\in V_0^{(0)}$ such that
\begin{equation}\label{1.74}
D_{r_1,1}(x^{\rho}).v_{-\rho,1_{[1]}}=x^{\rho}(\rho_{r_1}\ptl_1-\rho_1\ptl_{r_1}).v_{-\rho,1_{[1]}}=\rho_{r_1} v_{0,\textbf{0}}.
\end{equation}
Then we define $v_{0,\textbf{i}}$ for $\textbf{i}\in\mbb{N}^{l_1+l_2}$ as:

\begin{defi}\label{d:4.4}
Choose $0\not=v_{0,\textbf{0}}\in V_0^{(0)}$ as in (\ref{1.74}). For $k\geq 1$, we define
\begin{eqnarray}
v_{0,k_{[1]}}&=&\frac{1}{\rho_{r_1}(k+1)}x^{\rho}(\rho_{r_1}\ptl_1-\rho_1\ptl_{r_1}).v_{-\rho,(k+1)_{[1]}}\nonumber\\
&=&\frac{1}{\rho_{r_1}(k+1)}D_{r_1,1}(x^{\rho}).v_{-\rho,(k+1)_{[1]}},\label{1.76}
\end{eqnarray}
where $\rho$ and $r_1$ are the fixed elements in (\ref{4.8}) and (\ref{4.9}).
Moreover, we define
\begin{equation}\label{1.77}
v_{0,\textbf{i}}=\frac{i_1!}{k!}(t^{1_{[l_1+l_2]}}\ptl_1)^{i_{l_1+l_2}}.(t^{1_{[l_1+l_2-1]}}\ptl_1)^{i_{l_1+l_2-1}}.\cdots .(t^{1_{[2]}}\ptl_1)^{i_2}.v_{0, k_{[1]}}
\end{equation}
for $\textbf{i}=(i_1,i_2,\cdots,i_{l_1+l_2},0,\cdots,0)\in\mbb{N}^{l_1+l_2}$ with $|\textbf{i}|=k\geq 1$ and $\textbf{i}\not=k_{[1]}$.
\end{defi}\vspace{0.3cm}

It is not straightforward to verify the rationality of the definition (\ref{1.76}). It can be verified only until Lemma \ref{le:4.8}. \vspace{0.3cm}

\begin{lemma}\label{le:4.5}
For any $p\in\ove{1,l}$ and $\textbf{i}\in\mbb{N}^{l_1+l_2}$, we have
\begin{equation}
\ptl_p.v_{0,\textbf{i}}=i_p v_{0,\textbf{i}-1_{[p]}}.
\end{equation}
\end{lemma}

\noindent{\bf Proof.} When $\textbf{i}=\textbf{0}$, we have $v_{0,\textbf{0}}\in V_0^{(0)}$, which indicates
\begin{equation}\label{4.10}
\ptl_p.v_{0,\textbf{0}}=0 \quad \textrm{ for all } p\in\ove{1,l}
\end{equation}
by Lemma \ref{le:4.2}.
When $\textbf{i}=k_{[1]}$ with $k\geq 1$, we have
\begin{eqnarray}
\ptl_p.v_{0,k_{[1]}}
&=&\frac{1}{\rho_{r_1}(k+1)} \ptl_p. x^{\rho}(\rho_{r_1}\ptl_1-\rho_1\ptl_{r_1}).v_{-\rho,(k+1)_{[1]}}\nonumber\\
&=&\frac{1}{\rho_{r_1}(k+1)} (k+1)\delta_{p,1}x^{\rho}(\rho_{r_1}\ptl_1-\rho_1\ptl_{r_1}).v_{-\rho, k_{[1]}}\nonumber\\
&=& \delta_{p,1}k v_{0,(k-1)_{[1]}}, \qquad\forall \; p\in\ove{1,l} \label{1.80}
\end{eqnarray}
by (\ref{1.74}), (\ref{1.76})  and Lemma \ref{le:3.8}. When $\textbf{i}=(i_1,i_2,\cdots,i_{l_1+l_2},0,\cdots,0)$ with $|\textbf{i}|=k\geq 1$ and $\textbf{i}\not=k_{[1]}$, (\ref{1.77}), (\ref{1.80}) and Lemma \ref{le:4.2} show
\begin{eqnarray}
\ptl_p.v_{0,\textbf{i}}&=&\ptl_p.(\frac{i_1!}{k!}(t^{1_{[l_1+l_2]}}\ptl_1)^{i_{l_1+l_2}}.\cdots .(t^{1_{[2]}}\ptl_1)^{i_2}.v_{0, k_{[1]}})\nonumber\\
&=&\frac{i_1!}{k!}(t^{1_{[l_1+l_2]}}\ptl_1)^{i_{l_1+l_2}}.\cdots .(t^{1_{[2]}}\ptl_1)^{i_2}.(\ptl_p.v_{0, k_{[1]}})\nonumber\\
&=&\delta_{p,1}k \cdot\frac{i_1!}{k!} (t^{1_{[l_1+l_2]}}\ptl_1)^{i_{l_1+l_2}}.\cdots .(t^{1_{[2]}}\ptl_1)^{i_2}.v_{0, (k-1)_{[1]}}\nonumber\\
&=&\delta_{p,1}i_1 v_{0, \textbf{i}-1_{[1]}}\nonumber\\
&=&i_p v_{0,\textbf{i}-1_{[p]}} \qquad\textrm{ for } p\in \{1\}\cup\ove{l_1+l_2+1,l},\label{1.81}
\end{eqnarray}
and
\begin{eqnarray}
\ptl_p.v_{0,\textbf{i}}&=&\ptl_p.(\frac{i_1!}{k!}(t^{1_{[l_1+l_2]}}\ptl_1)^{i_{l_1+l_2}}.\cdots .(t^{1_{[2]}}\ptl_1)^{i_2}.v_{0, k_{[1]}})\nonumber\\
&=&\frac{i_1!}{k!}\Big((t^{1_{[l_1+l_2]}}\ptl_1)^{i_{l_1+l_2}}.\cdots .(t^{1_{[2]}}\ptl_1)^{i_2}.(\ptl_p.v_{0, k_{[1]}})\nonumber\\
& &+ i_p (t^{1_{[l_1+l_2]}}\ptl_1)^{i_{l_1+l_2}}.\cdots .(t^{1_{[p]}}\ptl_1)^{i_p-1}. \cdots.(t^{1_{[2]}}\ptl_1)^{i_2}.(\ptl_1.v_{0, k_{[1]}}) \Big)\nonumber\\
&=& i_p k \cdot\frac{i_1!}{k!} (t^{1_{[l_1+l_2]}}\ptl_1)^{i_{l_1+l_2}}.\cdots .(t^{1_{[p]}}\ptl_1)^{i_p-1}. \cdots.(t^{1_{[2]}}\ptl_1)^{i_2}.v_{0, (k-1)_{[1]}}\nonumber\\
&=&i_p v_{0,\textbf{i}-1_{[p]}} \qquad\textrm{ for } p\in \ove{2,l_1+l_2}.\label{1.82}
\end{eqnarray}
Thus this lemma follows from (\ref{4.10})--(\ref{1.82}). $\qquad\Box$ \vspace{0.3cm}

\begin{lemma}\label{le:4.6}
For any $\textbf{i}\in\mbb{N}^{l_1+l_2}$, $k\in\mbb{N}$, $p,q\in\ove{1,l_1+l_2}$ and $p',q'\in \ove{1,l}$, we have
\begin{equation}\label{4.6}
t^{1_{[p]}}\ptl_{q'}.v_{0,\textbf{i}}=i_{q'} v_{0,\textbf{i}+1_{[p]}-1_{[q']}},\ \ p\not=q',
\end{equation}
\begin{equation}\label{4.7}
(t^{1_{[p]}}\ptl_p - t^{1_{[q]}}\ptl_q).v_{0,\textbf{i}}=(i_p-i_q)v_{0,\textbf{i}},\ \ p\not=q,
\end{equation}
and
\begin{equation}\label{1.78}
D_{p',q'}(x^{\rho}).v_{-\rho,(k+1)_{[1]}}=(k+1)(\rho_{p'}\delta_{q',1}-\rho_{q'}\delta_{p',1})v_{0,k_{[1]}},\ \ p'\not=q',
\end{equation}
where $\rho$ is the fixed element in (\ref{4.9}).
\end{lemma}

\noindent{\bf Proof.} When $\textbf{i}=\textbf{0}$, (\ref{4.6}) and (\ref{4.7}) follow from Lemma \ref{le:4.2}. So we only need to consider the case $\textbf{i}\not=\textbf{0}$ for (\ref{4.6}) and (\ref{4.7}).
We divide the proof into several steps.\vspace{0.2cm}

{\it Step 1.} $t^{1_{[p]}}\ptl_1.v_{0,\textbf{i}}=i_1v_{0,\textbf{i}+1_{[p]}-1_{[1]}}$ for $\textbf{0}\not=\textbf{i}\in\mbb{N}^{l_1+l_2}$ and $p \in\ove{2,l_1+l_2}$.\vspace{0.2cm}

For any $p \in\ove{2,l_1+l_2}$, $\textbf{i}=(i_1,i_2,\cdots,i_{l_1+l_2},0,\cdots,0)\in\mbb{N}^{l_1+l_2}$ with $|\textbf{i}|=k>0$, we have
\begin{eqnarray}
t^{1_{[p]}}\ptl_1.v_{0,\textbf{i}} & = & \frac{i_1!}{k!}(t^{1_{[l_1+l_2]}}\ptl_1)^{i_{l_1+l_2}}.\cdots .(t^{1_{[p]}}\ptl_1)^{i_p+1}. \cdots.(t^{1_{[2]}}\ptl_1)^{i_2}.v_{0, k_{[1]}}\nonumber\\
& =& i_1v_{0,\textbf{i}+1_{[p]}-1_{[1]}}
\end{eqnarray}
by (\ref{1.77}).\vspace{0.2cm}

{\it Step 2.}  For any $\textbf{0}\not=\textbf{i}\in\mbb{N}^{l_1+l_2}$, $p,q \in\ove{2,l_1+l_2}$ and $q'\in \ove{2,l}$ with $p\not=q'$ and $p\not=q$ respectively, we have
\begin{equation}\label{1.84}
t^{1_{[p]}}\ptl_{q'}.v_{0,\textbf{i}}=i_{q'}v_{0,\textbf{i}+1_{[p]}-1_{[q']}} \textrm{ and }
(t^{1_{[p]}}\ptl_p - t^{1_{[q]}}\ptl_q).v_{0,\textbf{i}}=(i_p-i_q)v_{0,\textbf{i}}.
\end{equation}

{\it Firstly}, we want to prove
\begin{equation}\label{1.83}
t^{1_{[p]}}\ptl_{q'}.v_{0,k_{[1]}}=0 \textrm{ and }
(t^{1_{[p]}}\ptl_p - t^{1_{[q]}}\ptl_q).v_{0,k_{[1]}}=0
\end{equation}
 for $k>0$, $p,q \in\ove{2,l_1+l_2}$ and $q'\in \ove{2,l}$ with $p\not=q'$ and $p\not=q$ respectively. We shall manage that by induction on $k$.

When $k=1$, for any $p,q \in\ove{2,l_1+l_2}$ and $q'\in \ove{2,l}$ with $p\not=q'$ and $p\not=q$, we have
\begin{equation}\label{4.3}
t^{1_{[p]}}\ptl_{q'}.v_{0,1_{[1]}}, \ (t^{1_{[p]}}\ptl_p - t^{1_{[q]}}\ptl_q).v_{0,1_{[1]}}\; \in V_0^{(0)}
\end{equation}
by (\ref{2.2}), Lemma \ref{le:4.2} and Lemma \ref{le:4.5}.
So (\ref{4.3}) and Lemma \ref{le:4.2} imply
\begin{equation}
t^{1_{[p]}}\ptl_{q'}.v_{0,1_{[1]}}= -\frac{1}{2}[t^{1_{[p]}}\ptl_{q'},(t^{1_{[p]}}\ptl_p - t^{1_{[q']}}\ptl_{q'})].v_{0,1_{[1]}}=0
\end{equation}
for $p \in\ove{2,l_1+l_2}$ and $q' \in\ove{2,l_1+l_2}\backslash\{p\}$, and
\begin{equation}
t^{1_{[p]}}\ptl_{q'}.v_{0,1_{[1]}}=[t^{1_{[p]}}\ptl_s,t^{1_{[s]}}\ptl_{q'}].v_{0,1_{[1]}}=0
\end{equation}
for $p \in\ove{2,l_1+l_2}$ and $q' \in\ove{l_1+l_2+1,l}$, where $s \in\ove{2,l_1+l_2}\backslash\{p\}$.
Moreover,  by (\ref{4.3}) and Lemma \ref{le:4.2} again, we have
\begin{equation}
(t^{1_{[p]}}\ptl_p - t^{1_{[q]}}\ptl_q).v_{0,1_{[1]}}= [t^{1_{[p]}}\ptl_q,t^{1_{[q]}}\ptl_p] .v_{0,1_{[1]}}=0
\end{equation}
for $p,q \in\ove{2,l_1+l_2}$ with $p\not=q$. Thus, (\ref{1.83}) holds when $k=1$.

Assume that, with some $k\geq 1$,
\begin{equation}\label{4.5}
t^{1_{[p]}}\ptl_{q'}.v_{0,k_{[1]}}=0 \textrm{ and }
(t^{1_{[p]}}\ptl_p - t^{1_{[q]}}\ptl_q).v_{0,k_{[1]}}=0
 \end{equation}
 for any $p,q \in\ove{2,l_1+l_2}$ and $q'\in \ove{2,l}$ with $p\not=q'$ and $p\not=q$ respectively. Then the induction hypothesis (\ref{4.5}) and Lemma \ref{le:4.5} give
\begin{equation}\label{4.4}
t^{1_{[p]}}\ptl_{q'}.v_{0,(k+1)_{[1]}}, \ (t^{1_{[p]}}\ptl_p - t^{1_{[q]}}\ptl_q).v_{0,(k+1)_{[1]}}\; \in V_0^{(0)}
\end{equation}
 for any $p,q \in\ove{2,l_1+l_2}$ and $q'\in \ove{2,l}$ with $p\not=q'$ and $p\not=q$. Thus (\ref{4.4}) and Lemma \ref{le:4.2} imply
\begin{equation}
t^{1_{[p]}}\ptl_{q'}.v_{0,(k+1)_{[1]}}= -\frac{1}{2}[t^{1_{[p]}}\ptl_{q'},(t^{1_{[p]}}\ptl_p - t^{1_{[q']}}\ptl_{q'})].v_{0,(k+1)_{[1]}}=0
\end{equation}
for $p \in\ove{2,l_1+l_2}$ and $q' \in\ove{2,l_1+l_2}\backslash\{p\}$, and
\begin{equation}
t^{1_{[p]}}\ptl_{q'}.v_{0,(k+1)_{[1]}}= [t^{1_{[p]}}\ptl_s,t^{1_{[s]}}\ptl_{q'}].v_{0,(k+1)_{[1]}}=0
\end{equation}
for $p \in\ove{2,l_1+l_2}$ and $q' \in\ove{l_1+l_2+1,l}$, where $s \in\ove{2,l_1+l_2}\backslash\{p\}$.
Moreover, by (\ref{4.4}) and Lemma \ref{le:4.2} again, we have
\begin{equation}
(t^{1_{[p]}}\ptl_p - t^{1_{[q]}}\ptl_q).v_{0,(k+1)_{[1]}}= [t^{1_{[p]}}\ptl_q,t^{1_{[q]}}\ptl_p] .v_{0,(k+1)_{[1]}}=0
\end{equation}
for $p,q \in\ove{2,l_1+l_2}$ with $p\not=q$.
So (\ref{1.83}) holds with $k+1$. Thus (\ref{1.83}) follows from induction on $k$.

{\it Next}, we shall prove (\ref{1.84}) for $\textbf{i}=(i_1,i_2,\cdots,i_{l_1+l_2},0,\cdots,0)\in\mbb{N}^{l_1+l_2}$ with $i_s\not=0$ for some $s\in\ove{2,l_1+l_2}$.

For any $p\in\ove{2,l_1+l_2}$ and $q'\in \ove{2,l}\backslash\{p\}$, we have
\begin{eqnarray}
t^{1_{[p]}}\ptl_{q'}.v_{0,\textbf{i}}
&=&\frac{i_1!}{k!}t^{1_{[p]}}\ptl_{q'}.(t^{1_{[l_1+l_2]}}\ptl_1)^{i_{l_1+l_2}}.\cdots .(t^{1_{[2]}}\ptl_1)^{i_2}.v_{0, k_{[1]}}\nonumber\\
&=&\frac{i_1!}{k!}(t^{1_{[l_1+l_2]}}\ptl_1)^{i_{l_1+l_2}}.\cdots .(t^{1_{[2]}}\ptl_1)^{i_2}.(t^{1_{[p]}}\ptl_{q'}.v_{0, k_{[1]}})\nonumber\\
& &+\sum_{q=2}^{l_1+l_2} \delta_{q',q} i_q\frac{i_1!}{k!}(t^{1_{[l_1+l_2]}}\ptl_1)^{i_{l_1+l_2}}.\cdots .(t^{1_{[q]}}\ptl_1)^{i_q-1}.\nonumber\\
& &\cdots.(t^{1_{[p]}}\ptl_1)^{i_p+1}.\cdots .(t^{1_{[2]}}\ptl_1)^{i_2}.v_{0, k_{[1]}}\nonumber\\
&=& i_{q'}v_{0,\textbf{i}+1_{[p]}-1_{[q']}}\label{1.85}
\end{eqnarray}
by (\ref{1.77}) and (\ref{1.83}), where $k=|\textbf{i}|$ and $i_{q'}=0$ for $q' \in\ove{l_1+l_2+1,l}$. Thus (\ref{1.85}) imply
\begin{equation}\label{1.86}
(t^{1_{[p]}}\ptl_p - t^{1_{[q]}}\ptl_q).v_{0,\textbf{i}}= [t^{1_{[p]}}\ptl_q,t^{1_{[q]}}\ptl_p] .v_{0,\textbf{i}}=(i_p-i_q)v_{0,\textbf{i}}
\end{equation}
for any $p,q \in\ove{2,l_1+l_2}$ with $p\not=q$. So this step follows from (\ref{1.83}), (\ref{1.85}) and (\ref{1.86}).\vspace{0.2cm}

{\it Step 3.} $t^{1_{[1]}}\ptl_{q'}.v_{0,k_{[1]}}=0$ for all $q' \in\ove{2,l}$ and $k>0$.\vspace{0.2cm}

We shall give the proof by induction on $k$. When $k=1$, Lemmas \ref{le:4.2} and \ref{le:4.5} indicate
\begin{equation}\label{1.87}
t^{1_{[1]}}\ptl_{q'}.v_{0,1_{[1]}} \; \in V_0^{(0)},\quad \forall q' \in\ove{2,l}.
\end{equation}
So for any $q' \in\ove{2,l}$, it can be derived from (\ref{1.87}), Lemma \ref{le:4.2} and Step 2 of this lemma that
\begin{equation}
t^{1_{[1]}}\ptl_{q'}.v_{0,1_{[1]}}=[t^{1_{[1]}}\ptl_s, t^{1_{[s]}}\ptl_{q'}] .v_{0,1_{[1]}}= -t^{1_{[s]}}\ptl_{q'} . ( t^{1_{[1]}}\ptl_s .v_{0,1_{[1]}})=0,
\end{equation}
where $s\in\ove{2,l_1+l_2}\backslash\{q'\}$. Namely, this step holds for $k=1$.

Assume that, with some $k>0$,
\begin{equation}
t^{1_{[1]}}\ptl_{q'}.v_{0,k_{[1]}}=0 \quad\textrm{ for all } q' \in\ove{2,l}.
\end{equation}
Then it follows that
\begin{equation}\label{1.88}
t^{1_{[1]}}\ptl_{q'}.v_{0,(k+1)_{[1]}} \; \in V_0^{(0)},\quad \forall q' \in\ove{2,l}.
\end{equation}
Thus by Step 2 of this lemma, (\ref{1.88}) and Lemma \ref{le:4.2}, we obtain
\begin{eqnarray}
t^{1_{[1]}}\ptl_{q'}.v_{0,(k+1)_{[1]}}&=&[t^{1_{[1]}}\ptl_s, t^{1_{[s]}}\ptl_{q'}] .v_{0,(k+1)_{[1]}}\nonumber\\
&=& -t^{1_{[s]}}\ptl_{q'} . ( t^{1_{[1]}}\ptl_s .v_{0,(k+1)_{[1]}})=0
\end{eqnarray}
for $q' \in\ove{2,l}$, where $s\in\ove{2,l_1+l_2}\backslash\{q'\}$. So this step follows from induction on $k$.\vspace{0.2cm}

{\it Step 4.} For $k\geq 0$ and $p',q'\in\ove{1,l}$ with $p'\not=q'$, we have
\begin{equation}\label{1.92}
D_{p',q'}(x^{\rho}).v_{-\rho,(k+1)_{[1]}}=(k+1)(\rho_{p'}\delta_{q',1}-\rho_{q'}\delta_{p',1})v_{0,k_{[1]}},
\end{equation}
where $\rho$ is the fixed element in (\ref{4.9}).\vspace{0.2cm}

Notice that (\ref{1.74}) coincides with (\ref{1.76}) for $k=0$.
For any $k\geq 0$ and $q'\in\ove{2,l}\backslash\{r_1\}$, where $r_1$ is the fixed element in (\ref{4.8}), we obtain
\begin{eqnarray}
& &D_{r_1,q'}(x^{\rho}).v_{-\rho,(k+1)_{[1]}}\nonumber\\
&=&x^{\rho}(\rho_{r_1}\ptl_{q'}-\rho_{q'}\ptl_{r_1}).v_{-\rho,(k+1)_{[1]}}\nonumber\\
&=&\frac{1}{\rho_{r_1}}[x^{\rho}(\rho_{r_1}\ptl_1-\rho_1\ptl_{r_1}),t^{1_{[1]}}(\rho_{r_1}\ptl_{q'}-\rho_{q'}\ptl_{r_1})].v_{-\rho,(k+1)_{[1]}}\nonumber\\
&=&-\frac{1}{\rho_{r_1}} t^{1_{[1]}}(\rho_{r_1}\ptl_{q'}-\rho_{q'}\ptl_{r_1}).
(x^{\rho}(\rho_{r_1}\ptl_1-\rho_1\ptl_{r_1}).v_{-\rho,(k+1)_{[1]}})\nonumber\\
&=&-(k+1) t^{1_{[1]}}(\rho_{r_1}\ptl_{q'}-\rho_{q'}\ptl_{r_1}).
v_{0,k_{[1]}}\nonumber\\
&=& 0 \label{1.79}
\end{eqnarray}
by (\ref{1.74}), (\ref{1.76}), Lemma \ref{le:3.6} and Step 3. So (\ref{1.74}), (\ref{1.76}) and (\ref{1.79}) give
\begin{eqnarray}
 D_{p',q'}(x^{\rho}).v_{-\rho,(k+1)_{[1]}} &=& \frac{1}{\rho_{r_1}}(\rho_{p'} D_{r_1,q'}(x^{\rho}) - \rho_{q'} D_{r_1,p'}(x^{\rho})).v_{-\rho,(k+1)_{[1]}}\nonumber\\
 &=&(k+1)(\rho_{p'}\delta_{q',1}-\rho_{q'}\delta_{p',1})v_{0,k_{[1]}}
\end{eqnarray}
for $p',q'\in\ove{1,l}$ with $p'\not=q'$. This completes the proof of the step.\vspace{0.2cm}

{\it Step 5.} $(t^{1_{[1]}}\ptl_1-t^{1_{[2]}}\ptl_2).v_{0,k_{[1]}}=kv_{0,k_{[1]}}$ for $k>0$.\vspace{0.2cm}

Let $\rho$ and $r_1,r_2$ be the fixed elements in (\ref{4.8}) and (\ref{4.9}). Pick $r\in\{r_1,r_2\}\backslash\{2\}$. Then $r\not\in\{1,2\}$. Set
\begin{equation}
\tilde{\ptl}_1=\rho_r\ptl_2-\rho_2\ptl_r \ \textrm{ and } \tilde{\ptl}_2=\rho_r\ptl_1-\rho_1\ptl_r.
\end{equation}
 Then
\begin{equation}
t^{1_{[1]}}\ptl_1-t^{1_{[2]}}\ptl_2
=\frac{1}{\rho_r}( t^{1_{[1]} }\tilde{\ptl}_2 -t^{ 1_{[2]}}\tilde{\ptl}_1+ \rho_1 t^{1_{[1]}}\ptl_r-\rho_2 t^{1_{[2]}}\ptl_r),
\end{equation}
where $t^{1_{[1]} }\tilde{\ptl}_2 -t^{ 1_{[2]}}\tilde{\ptl}_1=\frac{1}{\rho_r}(\tilde{\ptl}_1(t^{1_{[1]}+1_{[2]}})\tilde{\ptl}_2
-\tilde{\ptl}_2(t^{1_{[1]}+1_{[2]}})\tilde{\ptl}_1)$. Thus (\ref{1.76}), Step 2--4, Lemma \ref{le:3.6} and Corollary \ref{c:3.7} give
\begin{eqnarray}
& &(t^{1_{[1]}}\ptl_1-t^{1_{[2]}}\ptl_2).v_{0,k_{[1]}}\nonumber\\
& = &\frac{1}{\rho_r}( t^{1_{[1]} }\tilde{\ptl}_2 -t^{ 1_{[2]}}\tilde{\ptl}_1+ \rho_1 t^{1_{[1]}}\ptl_r-\rho_2 t^{1_{[2]}}\ptl_r).v_{0,k_{[1]}}\nonumber\\
& = &\frac{1}{\rho_r\rho_{r_1}(k+1)}( t^{1_{[1]} }\tilde{\ptl}_2 -t^{ 1_{[2]}}\tilde{\ptl}_1).x^{\rho}(\rho_{r_1}\ptl_1-\rho_1\ptl_{r_1}).v_{-\rho,(k+1)_{[1]}}\nonumber\\
& = &\frac{1}{\rho_r\rho_{r_1}(k+1)}\Big( x^{\rho}(\rho_{r_1}\ptl_1-\rho_1\ptl_{r_1}). (t^{1_{[1]} }\tilde{\ptl}_2 -t^{1_{[2]}}\tilde{\ptl}_1).v_{-\rho,(k+1)_{[1]}} \nonumber\\
& & -x^{\rho}(\rho_{r_1}\tilde{\ptl}_2 +\delta_{r_1,2}\rho_1\tilde{\ptl}_1)
.v_{-\rho,(k+1)_{[1]}}\Big)\nonumber\\
& = & \frac{1}{\rho_{r_1}}x^{\rho}(\rho_{r_1}\ptl_1-\rho_1\ptl_{r_1}).v_{-\rho,(k+1)_{[1]}} -v_{0,k_{[1]}}\nonumber\\
& = &kv_{0,k_{[1]}}
\end{eqnarray}
for $k>0$.\vspace{0.2cm}

{\it Step 6.} $t^{1_{[1]}}\ptl_{q}.v_{0,\textbf{i}}=i_{q} v_{0,\textbf{i}+1_{[1]}-1_{[q]}}$ for $\textbf{0}\not=\textbf{i}=(i_1,i_2,\cdots,i_{l_1+l_2},0,\cdots,0)\in\mbb{N}^{l_1+l_2}$ and $q \in\ove{2,l}$.\vspace{0.2cm}

First of all, by (\ref{1.77}), Steps 2, 3 and 5, we obtain
\begin{eqnarray}
t^{1_{[1]}}\ptl_2.v_{0,\textbf{i}}
&=&\frac{i_1!}{k!}t^{1_{[1]}}\ptl_2.(t^{1_{[l_1+l_2]}}\ptl_1)^{i_{l_1+l_2}}.\cdots .(t^{1_{[2]}}\ptl_1)^{i_2}.v_{0, k_{[1]}}\nonumber\\
&=&\frac{i_1!}{k!}\big(-\sum_{s=3}^{l_1+l_2}i_s i_2 (t^{1_{[l_1+l_2]}}\ptl_1)^{i_{l_1+l_2}}.\cdots .(t^{1_{[2]}}\ptl_1)^{i_2-1}.v_{0, k_{[1]}}\nonumber\\
& & -i_2(i_2-1)(t^{1_{[l_1+l_2]}}\ptl_1)^{i_{l_1+l_2}}.\cdots .(t^{1_{[2]}}\ptl_1)^{i_2-1}.v_{0, k_{[1]}}\nonumber\\
& &+i_2(t^{1_{[l_1+l_2]}}\ptl_1)^{i_{l_1+l_2}}.\cdots .(t^{1_{[2]}}\ptl_1)^{i_2-1}.(t^{1_{[1]}}\ptl_1-t^{1_{[2]}}\ptl_2).v_{0, k_{[1]}}\big)\nonumber\\
&=&i_2\frac{(i_1+1)!}{k!}(t^{1_{[l_1+l_2]}}\ptl_1)^{i_{l_1+l_2}}.\cdots .(t^{1_{[2]}}\ptl_1)^{i_2-1}.v_{0, k_{[1]}}\nonumber\\
&=& i_2 v_{0,\textbf{i}+1_{[1]}-1_{[2]}},\label{1.89}
\end{eqnarray}
where $k=|\textbf{i}|$. Then (\ref{1.89}) and Step 2 tell
\begin{eqnarray}
t^{1_{[1]}}\ptl_q.v_{0,\textbf{i}}
&=&[t^{1_{[1]}}\ptl_2,t^{1_{[2]}}\ptl_q].v_{0,\textbf{i}}\nonumber\\
&=&t^{1_{[1]}}\ptl_2.(i_q v_{0,\textbf{i}+1_{[2]}-1_{[q]}})
-t^{1_{[2]}}\ptl_q.(i_2 v_{0,\textbf{i}+1_{[1]}-1_{[2]}})\nonumber\\
&=&i_q(i_2+1) v_{0,\textbf{i}+1_{[1]}-1_{[q]}}
- i_q i_2 v_{0,\textbf{i}+1_{[1]}-1_{[q]}}\nonumber\\
&=& i_q v_{0,\textbf{i}+1_{[1]}-1_{[q]}}
\end{eqnarray}
for $q\in\ove{3,l}$. So this step holds.\vspace{0.2cm}

{\it Step 7.} $(t^{1_{[1]}}\ptl_1 - t^{1_{[p]}}\ptl_p).v_{0,\textbf{i}}=(i_1-i_p)v_{0,\textbf{i}}$ for $\textbf{0}\not=\textbf{i} \in\mbb{N}^{l_1+l_2}$ and $p \in\ove{2,l_1+l_2}$.\vspace{0.2cm}

By Step 1 and Step 6, we have
\begin{eqnarray}
(t^{1_{[1]}}\ptl_1 - t^{1_{[p]}}\ptl_p).v_{0,\textbf{i}}
&=&[t^{1_{[1]}}\ptl_p,t^{1_{[p]}}\ptl_1].v_{0,\textbf{i}}\nonumber\\
&=& (i_1-i_p)v_{0,\textbf{i}}
\end{eqnarray}
 for $\textbf{0}\not=\textbf{i} \in\mbb{N}^{l_1+l_2}$ and $p \in\ove{2,l_1+l_2}$.

In summary, Steps 1, 2, 3 and 6 show (\ref{4.6}); Steps 2 and 7 show (\ref{4.6}); and Step 4 shows (\ref{1.78}). So we complete the proof of this lemma.
$\qquad\Box$

\begin{lemma}\label{le:4.7}
For $\textbf{i}\in\mbb{N}^{l_1+l_2}$ and $p,q\in\ove{1,l}$ with $p\not=q$, we have
\begin{equation}\label{1.90}
D_{p,q}(x^{\rho}).v_{-\rho,\textbf{i}}=\rho_p i_q v_{0,\textbf{i}-1_{[q]}}-\rho_q i_pv_{0,\textbf{i}-1_{[p]}},
\end{equation}
where $\rho$ is the fixed element in (\ref{4.9}).
\end{lemma}

\noindent{\bf Proof.} When $\textbf{i}=\textbf{0}$, (\ref{1.90}) follows from Lemma \ref{le:3.2}. So we only need to consider the case $\textbf{i}\not=\textbf{0}$.
Recall that Lemma \ref{le:4.6} has given the proof for $\textbf{i}=k_{[1]}$ with $k\geq 1$. Therefore, we can prove the lemma this way: for each $k\geq 1$, we verify (\ref{1.90}) for all $\textbf{0}\not=\textbf{i}\in\mbb{N}^{l_1+l_2}$ with $|\textbf{i}|=k$, by induction on $\textbf{i}$ with the total order defined in Definition \ref{d:1.5}.

Fix any $k\geq 1$. Assume that, (\ref{1.90}) holds for all $\textbf{j}<\textbf{i}=(i_1,i_2,\cdots,i_{l_1+l_2},0,\cdots,0)$ with $|\textbf{j}|=|\textbf{i}|=k$, where $i_p\not=0$ for some $p\in\ove{2,l_1+l_2}$. Next, we verify (\ref{1.90}) for $\textbf{i}$. Set
\begin{equation}
s=\max\{p\in\ove{2,l_1+l_2}\mid i_p\not=0\}.
 \end{equation}
 Then $s>1$.
 Pick $r\in\{r_1,r_2\}\backslash\{s\}$, where $r_1,r_2$ are the fixed element in (\ref{4.8}).
 Then $r\not=1$ and $\rho_r\not=0$. So by the total order defined in (\ref{2.1}), we have
\begin{equation}
 \textbf{i}+1_{[1]}-1_{[s]}<\textbf{i}, \quad \textbf{i}+1_{[1]}-1_{[r]}<\textbf{i}  \textrm{ or } i_r=0.
\end{equation}
Thus Lemma \ref{le:3.6}, Lemma \ref{le:4.6} and the induction hypothesis give
\begin{eqnarray}
& &D_{p,q}(x^{\rho}).v_{-\rho,\textbf{i}}\nonumber\\
&=&\frac{1}{\rho_r(i_1+1)}D_{p,q}(x^{\rho}).\big(t^{1_{[s]}}(\rho_r\ptl_1-\rho_1\ptl_r). v_{-\rho,\textbf{i}+1_{[1]}-1_{[s]}} + \rho_1 i_r v_{-\rho,\textbf{i}+1_{[1]}-1_{[r]}}\big)\nonumber\\
&=&\frac{1}{\rho_r(i_1+1)}\Big((\rho_p\delta_{q,s}-\rho_q\delta_{p,s})x^{\rho}(\rho_r\ptl_1-\rho_1\ptl_r). v_{-\rho,\textbf{i}+1_{[1]}-1_{[s]}}\nonumber\\
& & + t^{1_{[s]}}(\rho_r\ptl_1-\rho_1\ptl_r). D_{p,q}(x^{\rho}).v_{-\rho,\textbf{i}+1_{[1]}-1_{[s]}} + \rho_1 i_r D_{p,q}(x^{\rho}).v_{-\rho,\textbf{i}+1_{[1]}-1_{[r]}}\Big)\nonumber\\
&=&\rho_p i_q v_{0,\textbf{i}-1_{[q]}}-\rho_q i_pv_{0,\textbf{i}-1_{[p]}}
\end{eqnarray}
for $p,q\in\ove{1,l}$ with $p\not=q$, which coincides with (\ref{1.90}). Thus this lemma holds. $\qquad\Box$ \vspace{0.3cm}

\begin{lemma}\label{le:4.8}
For any $\al\in\G\backslash\{0\}$, $\textbf{i}\in\mbb{N}^{l_1+l_2}$ and $p,q\in\ove{1,l}$ with $p\not=q$, we have
\begin{equation}\label{1.91}
D_{p,q}(x^{\al}).v_{-\al,\textbf{i}}=\al_p i_q v_{0,\textbf{i}-1_{[q]}}-\al_q i_pv_{0,\textbf{i}-1_{[p]}}.
\end{equation}
\end{lemma}

\noindent{\bf Proof.} We shall give the proof by induction on $|\textbf{i}|$. When $|\textbf{i}|=0$, Lemma \ref{le:3.2} shows
 \begin{equation}
D_{p,q}(x^{\al}).v_{-\al,\textbf{0}}=0
\end{equation}
for all $\al\in\G\backslash\{0\}$ and $p,q\in\ove{1,l}$ with $p\not=q$, which coincides with (\ref{1.91}).
 Suppose that (\ref{1.91}) holds for $\textbf{i}\in\mbb{N}^{l_1+l_2}$ with $|\textbf{i}|\leq k$, where $k\geq 0$. Fix some $\sgm\in\G\backslash\{0\}$ such that
 \begin{equation}
 \ker\sgm\not=\ker\rho,
  \end{equation}
 where $\rho$ is the fixed element in (\ref{4.9}).
Then we prove (\ref{1.91}) for $\textbf{i}\in\mbb{N}^{l_1+l_2}$ with $|\textbf{i}|=k+1$ in two steps.\vspace{0.2cm}

{\it Step 1.}
For $\textbf{i}\in\mbb{N}^{l_1+l_2}$ with $|\textbf{i}|=k+1$, and $p,q\in\ove{1,l}$ with $p\not=q$, we have
\begin{equation}\label{1.93}
D_{p,q}(x^{\sgm}).v_{-\sgm,\textbf{i}}=\sgm_p i_q v_{0,\textbf{i}-1_{[q]}}-\sgm_q i_pv_{0,\textbf{i}-1_{[p]}}.
\end{equation}

Since $\ker\sgm\not=\ker\rho$, we have $\sgm\not=\pm \rho$, $\sgm_{s_1}\not=0$ for some $s_1\in\ove{l_1+1,l}$, and $\sgm_{s_1}\rho_{s_2}-\sgm_{s_2}\rho_{s_1}\not=0$ for some $s_2\in\ove{l_1+1,l}\backslash\{s_1\}$. Fix such $s_1$ and $s_2$.
Set
\begin{equation}\label{1.97}
\tilde{\ptl}_q=(\sgm_{s_1}\rho_{s_2}-\sgm_{s_2}\rho_{s_1})(\sgm_{s_1}\ptl_q-\sgm_q\ptl_{s_1})
-(\sgm_{s_1}\rho_q-\sgm_q\rho_{s_1})(\sgm_{s_1}\ptl_{s_2}-\sgm_{s_2}\ptl_{s_1})
 \end{equation}
 for $q\in\ove{1,l}\backslash\{s_1,s_2\}$, where $\rho$ is the fixed element in (\ref{4.9}).
 Then for any $q\in\ove{1,l}\backslash\{s_1,s_2\}$, we have \begin{equation}
 \tilde{\ptl}_q(\sgm)=\tilde{\ptl}_q(\rho)=0,
 \end{equation}
 \begin{equation}\label{1.94}
x^{\sgm}\tilde{\ptl}_q=((\sgm_{s_1}\rho_{s_2}-\sgm_{s_2}\rho_{s_1})D_{s_1,q}(x^{\sgm})
-(\sgm_{s_1}\rho_q-\sgm_q\rho_{s_1})D_{s_1,s_2}(x^{\sgm})),
\end{equation}
\begin{equation}\label{1.95}
x^{\rho}\tilde{\ptl}_q=\sgm_{s_1}\big(\sgm_{s_1}D_{s_2,q}(x^{\rho})
+\sgm_{s_2}D_{q,s_1}(x^{\rho})+\sgm_q D_{s_1,s_2}(x^{\rho})\big).
\end{equation}
  Let
 \begin{equation}
 \bar{\ptl}=\sgm_{s_1}\ptl_{s_2}-\sgm_{s_2}\ptl_{s_1} \textrm{ and } \bar{\ptl}'=\rho_{s_1}\ptl_{s_2}-\rho_{s_2}\ptl_{s_1}.
 \end{equation}
Then $\bar{\ptl}\in\ker\sgm\backslash\ker\rho$ and $\bar{\ptl}'\in\ker\rho\backslash\ker\sgm$. So Lemma \ref{le:3.9}, Lemma \ref{le:4.7} and the induction hypothesis give
\begin{eqnarray}
& &x^{\sgm}\tilde{\ptl}_q.v_{-\sgm,\textbf{i}}\nonumber\\
&=&\frac{1}{\bar{\ptl}(\rho)\bar{\ptl}'(\sgm)}
[x^{-\rho}\bar{\ptl}',[x^{\sgm}\bar{\ptl},x^{\rho}\tilde{\ptl}_q]].v_{-\sgm,\textbf{i}}\nonumber\\
&=&\frac{1}{\bar{\ptl}(\rho)\bar{\ptl}'(\sgm)}
\Big(\bar{\ptl}(\rho)x^{-\rho}\bar{\ptl}'.x^{\sgm+\rho}\tilde{\ptl}_q.v_{-\sgm,\textbf{i}}
-x^{\sgm}\bar{\ptl}.x^{\rho}\tilde{\ptl}_q.x^{-\rho}\bar{\ptl}'.v_{-\sgm,\textbf{i}}\nonumber\\
& & +x^{\rho}\tilde{\ptl}_q.x^{\sgm}\bar{\ptl}.x^{-\rho}\bar{\ptl}'.v_{-\sgm,\textbf{i}}\Big)\nonumber\\
&=&\sgm_{s_1}\Big((\sgm_{s_1}\rho_{s_2}-\sgm_{s_2}\rho_{s_1})i_q v_{0,\textbf{i}-1_{[q]}}
-(\sgm_{s_1}\rho_q-\sgm_q\rho_{s_1})i_{s_2} v_{0,\textbf{i}-1_{[s_2]}}\nonumber\\
& &+(\sgm_{s_2}\rho_q-\sgm_q\rho_{s_2})i_{s_1} v_{0,\textbf{i}-1_{[s_1]}}\Big)\label{1.96}
\end{eqnarray}
for $q\in\ove{1,l}\backslash\{s_1,s_2\}$, $\textbf{i}\in\mbb{N}^{l_1+l_2}$ with $|\textbf{i}|=k+1$. Thus by Lemma \ref{le:3.6}, Lemma \ref{le:4.6}, (\ref{1.97}) and (\ref{1.96}), we obtain
\begin{eqnarray}
& &D_{s_1,s_2}(x^{\sgm}).v_{-\sgm,\textbf{i}}\nonumber\\
&=&\frac{1}{\sgm_{s_1}(\sgm_{s_1}\rho_{s_2}-\sgm_{s_2}\rho_{s_1})}[x^{\sgm}\tilde{\ptl}_p,
t^{1_{[p]}}(\sgm_{s_1}\ptl_{s_2}-\sgm_{s_2}\ptl_{s_1})].v_{-\sgm,\textbf{i}}\nonumber\\
&=&\frac{1}{\sgm_{s_1}(\sgm_{s_1}\rho_{s_2}-\sgm_{s_2}\rho_{s_1})}
\Big(\sgm_{s_1}i_{s_2}x^{\sgm}\tilde{\ptl}_p.v_{-\sgm,\textbf{i}+1_{[p]}-1_{[s_2]}}
-\sgm_{s_2}i_{s_1}x^{\sgm}\tilde{\ptl}_p.v_{-\sgm,\textbf{i}+1_{[p]}-1_{[s_1]}}\nonumber\\
& &-t^{1_{[p]}}(\sgm_{s_1}\ptl_{s_2}-\sgm_{s_2}\ptl_{s_1}).
(x^{\sgm}\tilde{\ptl}_p.v_{-\sgm,\textbf{i}})\Big)\nonumber\\
&=& \sgm_{s_1}i_{s_2}v_{0,\textbf{i}-1_{[s_2]}}-\sgm_{s_2}i_{s_1}v_{0,\textbf{i}-1_{[s_1]}}\label{1.98}
\end{eqnarray}
for $\textbf{i}\in\mbb{N}^{l_1+l_2}$ with $|\textbf{i}|=k+1$, where $p\in\ove{1,l_1+l_2}\backslash\{s_1,s_2\}$. Moreover, (\ref{1.94}), (\ref{1.96}) and (\ref{1.98}) imply
\begin{eqnarray}
D_{s_1,q}(x^{\sgm}).v_{-\sgm,\textbf{i}}
&=&\frac{1}{\sgm_{s_1}\rho_{s_2}-\sgm_{s_2}\rho_{s_1}}
(x^{\sgm}\tilde{\ptl}_q+(\sgm_{s_1}\rho_q-\sgm_q\rho_{s_1})D_{s_1,s_2}(x^{\sgm})).v_{-\sgm,\textbf{i}}\nonumber\\
&=&\sgm_{s_1}i_q v_{0,\textbf{i}-1_{[q]}}-\sgm_q i_{s_1}v_{0,\textbf{i}-1_{[s_1]}}\label{1.99}
\end{eqnarray}
for $q\in\ove{1,l}\backslash\{s_1,s_2\}$, $\textbf{i}\in\mbb{N}^{l_1+l_2}$ with $|\textbf{i}|=k+1$. Therefore, (\ref{1.98}) and (\ref{1.99}) indicate
\begin{eqnarray}
D_{p,q}(x^{\sgm}).v_{-\sgm,\textbf{i}}
&=&\frac{1}{\sgm_{s_1}}(\sgm_p D_{s_1,q}(x^{\sgm})-\sgm_q D_{s_1,p}(x^{\sgm})).v_{-\sgm,\textbf{i}}\nonumber\\
&=&\sgm_pi_q v_{0,\textbf{i}-1_{[q]}}-\sgm_q i_pv_{0,\textbf{i}-1_{[p]}}.
\end{eqnarray}
for $\textbf{i}\in\mbb{N}^{l_1+l_2}$ with $|\textbf{i}|=k+1$, $p,q\in\ove{1,l}$ with $p\not=q$. So this step holds.\vspace{0.2cm}

{\it Step 2.} For any $\al\in\G\backslash\{0\}$, $\textbf{i}\in\mbb{N}^{l_1+l_2}$ with $|\textbf{i}|=k+1$, and $p,q\in\ove{1,l}$ with $p\not=q$, we have
\begin{equation}
D_{p,q}(x^{\al}).v_{-\al,\textbf{i}}=\al_p i_q v_{0,\textbf{i}-1_{[q]}}-\al_q i_pv_{0,\textbf{i}-1_{[p]}}.
\end{equation}

Fix an arbitrary $\al\in\G\backslash\{0\}$. Then we have $\ker\al\not=\ker\rho$ or $\ker\al\not=\ker\sgm$, where $\rho$ is the fixed element in (\ref{4.9}).
 Choose $\tau\in\{\rho,\sgm\}$ such that $\ker\al\not=\ker\tau$. Then, with $\sgm$ replaced by $\al$ and $\rho$ by $\tau$ in Step 1, we can get
 \begin{equation}
D_{p,q}(x^{\al}).v_{-\al,\textbf{i}}=\al_p i_q v_{0,\textbf{i}-1_{[q]}}-\al_q i_pv_{0,\textbf{i}-1_{[p]}}
\end{equation}
for $\textbf{i}\in\mbb{N}^{l_1+l_2}$ with $|\textbf{i}|=k+1$, $p,q\in\ove{1,l}$ with $p\not=q$.
As $\al\in\G\backslash\{0\}$ is arbitrary, this completes the proof of the step.\vspace{0.2cm}

Thus (\ref{1.91}) holds for $\textbf{i}\in\mbb{N}^{l_1+l_2}$ with $|\textbf{i}|=k+1$. So this lemma follows from induction on $|\textbf{i}|$.
  $\qquad\Box$ \vspace{0.3cm}

\noindent{\bf Remark.} The above lemma proves the rationality of Definition \ref{d:4.4}. \vspace{0.3cm}

\begin{lemma}\label{le:4.9}
For any $\al\in\G\backslash\{0\}$, $\textbf{i}\in\mbb{N}^{l_1+l_2}$ and $p,q\in\ove{1,l}$ with $p\not=q$, we have
\begin{equation}\label{4.11}
D_{p,q}(x^{\al}).v_{0,\textbf{i}}=\al_p i_q v_{\al,\textbf{i}-1_{[q]}}-\al_q i_pv_{\al,\textbf{i}-1_{[p]}}.
\end{equation}
\end{lemma}

\noindent{\bf Proof.} When $\textbf{i}=\textbf{0}$, (\ref{4.11}) follows from Lemma \ref{le:4.2}. So we only need to consider the case $\textbf{i}\not=\textbf{0}$. Fix any $k>0$. It suffices to prove
\begin{equation}\label{1.100}
D_{p,q}(x^{\al}).v_{0,\textbf{i}}=\al_p i_q v_{\al,\textbf{i}-1_{[q]}}-\al_q i_pv_{\al,\textbf{i}-1_{[p]}}
\end{equation}
for all $\al\in\G\backslash\{0\}$, $\textbf{i}\in\mbb{N}^{l_1+l_2}$ with $|\textbf{i}|=k$, and $p,q\in\ove{1,l}$ with $p\not=q$. We shall prove (\ref{1.100}) in two steps.\vspace{0.2cm}

{\it Step 1.} For any $\al\in\G\backslash\{0,\rho\}$, $\textbf{i}\in\mbb{N}^{l_1+l_2}$ with $|\textbf{i}|=k$, $p,q\in\ove{1,l}$ with $p\not=q$, we have
\begin{equation}\label{1.101}
D_{p,q}(x^{\al}).v_{0,\textbf{i}}=\al_p i_q v_{\al,\textbf{i}-1_{[q]}}-\al_q i_pv_{\al,\textbf{i}-1_{[p]}},
\end{equation}
where $\rho$ is the fixed element in (\ref{4.9}).\vspace{0.2cm}

We prove this step by induction on $\textbf{i}$ with the total order defined in Definition \ref{d:1.5}. When $\textbf{i}=k_{[1]}$, by (\ref{1.76}) and Lemma \ref{le:3.9} we have
\begin{eqnarray}
& &D_{p,q}(x^{\al}).v_{0,k_{[1]}}\nonumber\\
&=&\frac{1}{\rho_{r_1}(k+1)}D_{p,q}(x^{\al}).(x^{\rho}(\rho_{r_1}\ptl_1-\rho_1\ptl_{r_1}).
v_{-\rho,(k+1)_{[1]}})\nonumber\\
&=&\frac{1}{\rho_{r_1}(k+1)}\Big(x^{\rho}(\rho_{r_1}\ptl_1-\rho_1\ptl_{r_1}).D_{p,q}(x^{\al}).
v_{-\rho,(k+1)_{[1]}}\nonumber\\
& &+x^{\al+\rho}\big((\al_p \rho_q -\al_q \rho_p)(\rho_{r_1}\ptl_1-\rho_1\ptl_{r_1})\nonumber\\
& &-(\rho_{r_1}\al_1-\rho_1\al_{r_1})(\al_p \ptl_q -\al_q \ptl_p)\big).v_{-\rho,(k+1)_{[1]}}\Big)\nonumber\\
&=&k(\al_p \delta_{q,1}-\al_q \delta_{p,1}) v_{\al,(k-1)_{[1]}}
\end{eqnarray}
for all $\al\in\G\backslash\{0,\rho\}$, $p,q\in\ove{1,l}$ with $p\not=q$, where $\rho$ and $r_1$ are the fixed elements in (\ref{4.8}) and (\ref{4.9}).
Namely, (\ref{1.101}) holds when $\textbf{i}=k_{[1]}$. Suppose that, (\ref{1.101}) holds for all $\textbf{j}<\textbf{i}=(i_1,i_2,\cdots,i_{l_1+l_2},0,\cdots,0)$ with $|\textbf{j}|=|\textbf{i}|=k$ and $i_p\not=0$ for some $p\in\ove{2,l_1+l_2}$. Recall that $r_1>1$ (cf. (\ref{4.9})). So by the total order defined in (\ref{2.1}), we have
\begin{equation}
 i_{r_1}=0 \textrm{ or } \textbf{i}+1_{[1]}-1_{[r_1]}<\textbf{i}  .
\end{equation}
Thus Lemma \ref{le:3.9}, Lemma \ref{le:4.7} and the induction hypothesis give
\begin{eqnarray}
& &D_{p,q}(x^{\al}).v_{0,\textbf{i}}\nonumber\\
&=&\frac{1}{\rho_{r_1}(i_1+1)}D_{p,q}(x^{\al}).\big(x^{\rho}(\rho_{r_1}\ptl_1-\rho_1\ptl_{r_1}). v_{-\rho,\textbf{i}+1_{[1]}} + \rho_1 i_{r_1} v_{0,\textbf{i}+1_{[1]}-1_{[r_1]}}\big)\nonumber\\
&=&\frac{1}{\rho_{r_1}(i_1+1)}\Big(x^{\rho}(\rho_{r_1}\ptl_1-\rho_1\ptl_{r_1}).x^{\al}(\al_p \ptl_q -\al_q \ptl_p).
v_{-\rho,\textbf{i}+1_{[1]}}\nonumber\\
& &+x^{\al+\rho}\big((\al_p \rho_q -\al_q \rho_p)(\rho_{r_1}\ptl_1-\rho_1\ptl_{r_1})-(\rho_{r_1}\al_1-\rho_1\al_{r_1})(\al_p \ptl_q -\al_q \ptl_p)\big).v_{-\rho,\textbf{i}+1_{[1]}} \nonumber\\
& & + \rho_1 i_{r_1} x^{\al}(\al_p \ptl_q -\al_q \ptl_p).v_{0,\textbf{i}+1_{[1]}-1_{[r_1]}}\Big)\nonumber\\
&=&\al_p i_q v_{\al,\textbf{i}-1_{[q]}}-\al_q i_pv_{\al,\textbf{i}-1_{[p]}}
\end{eqnarray}
for all $\al\in\G\backslash\{0,\rho\}$, $p,q\in\ove{1,l}$ with $p\not=q$. So (\ref{1.101}) holds by induction on $\textbf{i}$, which completes the proof of the step.\vspace{0.2cm}

{\it Step 2.} For all $\textbf{i}\in\mbb{N}^{l_1+l_2}$ with $|\textbf{i}|=k$, $p,q\in\ove{1,l}$ with $p\not=q$, we have
\begin{equation}\label{1.102}
D_{p,q}(x^{\rho}).v_{0,\textbf{i}}=\rho_p i_q v_{\rho,\textbf{i}-1_{[q]}}-\rho_q i_pv_{\rho,\textbf{i}-1_{[p]}},
\end{equation}
where $\rho$ is the fixed element in (\ref{4.9}).\vspace{0.2cm}

Replacing $\rho$ by $-\rho$, and $\al$ by $\rho$ in Step 1, we can similarly prove (\ref{1.102}). Here we omit the details.

So this lemma follows from Steps 1 and 2.
 $\qquad\Box$ \vspace{0.3cm}

Under the condition that $\mu=0$, we get a set $\{v_{\be,\textbf{i}}\mid \be\in\G,\textbf{i}\in\mbb{N}^{l_1+l_2} \}$ from Definitions \ref{d:3.3} and \ref{d:4.4}.
In analogy with Lemma \ref{le:3.10}, it can be deduced from Lemma \ref{le:3.8}, Lemma \ref{le:4.5} and Definition \ref{d:1.5} that:\vspace{0.2cm}

\begin{lemma}\label{le:4.10}
The set $\{v_{\be,\textbf{i}}\mid \be\in\G,\textbf{i}\in\mbb{N}^{l_1+l_2} \}$ is an $\mbb{F}$-basis of $V$.
\end{lemma}\vspace{0.2cm}

So Proposition \ref{t:2.1}, Lemma \ref{le:3.6}, Corollary \ref{c:3.7}, Lemma \ref{le:3.8}, Lemma \ref{le:3.9}
and Lemma \ref{le:4.5}--\ref{le:4.10} give\vspace{0.2cm}

\begin{lemma}\label{le:4.11}
If $\mu\in\G$, then $V=V(\mu)\simeq A_{\mu}\simeq A_{0}$.
\end{lemma}\vspace{0.2cm}

In summary, Lemmas \ref{le:3.1}, \ref{le:3.11} and \ref{le:4.11} show that, Theorem \ref{t:1.3} holds.


{\small \begin{thebibliography}{99}

\bibitem{BP}Jeffrey Bergen and D. S. Passman, Simple Lie algebras of Special type, {\it J. Algebra} {\bf 227} (2000), no. 1, 45-67.

\bibitem{C}L. Chen, Multiplicity-One Representations of Divergence-Free Lie Algebras, {\it J. Algebra} (2012), doi:10.1016/j.jalgebra.2011.12.005.

\bibitem{DZ1}D. Dokovic and K. Zhao, Derivations, isomorphisms, and second cohomology of generalized Witt algebras,
{\it Trans. Amer. Math. Soc.} {\bf 350(2)} (1998), 643-664.

\bibitem{DZ3}D. Dokovic and K. Zhao, Generalized Cartan type S Lie algebras in characteristic 0,
{\it J. Algebra} {\bf 193} (1997), 144-179.

\bibitem{Ka1}V. G. Kac, A description of filtered Lie algebras whose associated
graded Lie algebras are of Cartan types, {\it Math.~USSR-Izvestija}
{\bf 8} (1974), 801-835.

\bibitem{Ka2}V. G. Kac, Lie superalgebras, {\it Adv. Math.}
{\bf 26} (1977), 8-96.

\bibitem{Ka3}V. G. Kac, Classification of infinite-dimensional simple linearly compact Lie superalgebras, {\it Adv. Math.}
{\bf 139} (1998), 1-55.

\bibitem{Ka4} V. G. Kac, Some problems on infinite dimensional Lie algebras and their representations,
{\it Lect. Notes Math., Lie algebras and related topics,
Berlin Heidelberg New York: Springer} {\bf  933} (1981), 117-126.

\bibitem{KR} V. G. Kac, A. K. Raina, Highest weights representations of infinite dimensional Lie
algebras, {\it Adv. Ser. Math. Phys.} {\bf 2} (1988).

\bibitem{Kap} I. Kaplansky, The Virasoro algebra,
{\it Commun. Math. Phys.} {\bf 86} (1982), 49-54.

\bibitem{KS} I. Kaplansky and L. J. Santharoubane, Harish-Chandra modules over the Virasoro algebras,
{\it MSRI Publ.} {\bf 4} (1987), 217-231.


\bibitem{M} O.~Mathieu, Classification of Harish-Chandra modules over the Virasoro Lie algebra, {\it Inv. Math.} {\bf 107} (1992), 225-234.


\bibitem{Ma2} V.~Mazorchuk, Verma modules over generalized Witt algebras,
{\it Composito Math.} {\bf 115(1)} (1999), 21-35.

\bibitem{O}J. M. Osborn, New simple infinite-dimensional Lie algebras of
characteristic 0, {\it J. Algebra} {\bf 185} (1996), 820-835.

\bibitem{PS}I. Penkov and V. Serganova, Weight representations of the polynomial Catan type Lie algebras $W_{n}$ and $\bar{S}_{n}$,
{\it Math. Res. Lett.} {\bf 6} (1999), 397-416.

\bibitem{R1}S. E. Rao, Representations of Witt algebras, {\it Publ.~Rims, Kyoto Univ.} {\bf 29} (1994), 191-201.

\bibitem{R2}S. E. Rao, Irreducible representations of the Lie algebra of the
diffeomorphisms of a d-dimensional torus, {\it J. Algebra} {\bf 182}
(1996), 401-421.

\bibitem{Sg1}G. Shen, Graded modules of graded Lie algebras of Cartan type
(I)--mixed product of modules, {\it Sci.~China Ser.~A} {\bf 29}
(1986), 570-581.

\bibitem{Sg2}G. Shen, Graded modules of graded Lie algebras of Cartan type
(II)--positive and negative graded modules, {\it Sci.~China Ser.~A}
{\bf 29} (1986), 1009-1019.

\bibitem{Sg3}G. Shen, Graded modules of graded Lie algebras of Cartan type (III)--irreducible modules,
{\it Chinese Ann.~Math.~Ser.}~{\bf B9} (1988), 404-417.

\bibitem{S1}Y. Su, Harish-Chandra modules of the intermediate series over the high rank Virasoro algebras and high rank super-Virasoro algebras, {\it J. Math. Phys.} {\bf 35} (1994), 2013-2023.

\bibitem{S2}Y. Su, Classification of Harish-Chandra modules over the super-Virasoro algebras, {\it Comm. Algebra} {\bf
23} (1995), 3653-3675.

\bibitem{SX}Y. Su and X. Xu, Structure of divergence-free Lie
algebras, {\it J. Algebra} {\bf 243} (2001), 557-595.


\bibitem{SZh}Y. Su and J. Zhou, Some representations of nongraded Lie algebras of generalized Witt type, {\it
J. Algebra} {\bf 246} (2001), 721-738.

\bibitem{SZk}Y. Su and K. Zhao, Generalized Virasoro and super-Virasoro algebras and modules of intermediate series,
{\it J. Algebra} {\bf 252} (2002), 1-19.

\bibitem{X2}X. Xu, New generalized simple Lie algebras of Cartan type over a
field with characteristic 0, {\it J. Algebra} {\bf 224} (2000),
23-58.

\bibitem{X3}X. Xu, Quadratic conformal superalgebras, {\it J. Algebra} {\bf
231} (2000), 1-38.

\bibitem{X4}X. Xu, Equivalence of conformal superalgebras to Hamiltonian
superoperators, {\it Algebra Colloq.} {\bf 8} (2001), 63-92.

\bibitem{Zk1}K. Zhao, Generalized Cartan type S Lie algebras in characteristic zero, II, {\it Pacific J. Math.} {\bf
192} (2000), 431-454.

\bibitem{Zk2}K. Zhao, Weight modules over generalized Witt algebras with 1-dimensional weight spaces, {\it Forum
Math.} {\bf 16} (2004), 725-748. Preprint, 2000.

\bibitem{Zy1}Y. Zhao, Irreducible representations of nongraded Witt type Lie algebras,
{\it J. Algebra} {\bf 298} (2006), 540-562.

\bibitem{Zy2} Y. Zhao,
Representations of nongraded Lie algebras of Block type, {\it
Manuscripta math.}
{\bf 119} (2006), 183-216.

\bibitem{Zy3} Y. Zhao, Composition series for a family of modules
of nongraded Hamiltonian type Lie algebras, {\it J. Lie Theory} {\bf
19} (2009), 1-27.




\end{thebibliography}}
\end{document}